\documentclass{article}

\usepackage[T1]{fontenc}
\usepackage[numbers]{natbib}
\usepackage{mathrsfs}
\usepackage{amsfonts, amsmath, amssymb, amsthm}
\usepackage{booktabs, tabularx, array}

\usepackage{graphicx}
\usepackage{amsmath, amsfonts}
\usepackage{amsthm}
\usepackage{bm}
\usepackage{tikz}
\usetikzlibrary{positioning, calc, decorations.pathreplacing, arrows.meta}
\usepackage{subfig}
\usepackage[ruled,vlined]{algorithm2e}
\newtheorem{theorem}{Theorem}
\newtheorem{lemma}{Lemma}
\newtheorem{corollary}{Corollary}
\newtheorem{remark}{Remark}
\DeclareMathOperator*{\argmin}{arg\,min}
\usepackage[margin=1.75cm]{geometry}

\begin{document}



\title{POD-CNF: reduced order modeling of parametrized stochastic PDEs via
normalizing flows and proper orthogonal decomposition}  



%

\author{Andrea Violante$^1$, Nicola Rares Franco$^1$, Andrea Manzoni$^1$}

\date{\small $^1$MOX, Department of Mathematics, Politecnico di Milano, Piazza Leonardo da Vinci 32, 20133, Milan, Italy}

















\maketitle

\begin{abstract}
Reduced-order models (ROMs) learned from trusted numerical simulations enable the real-time simulation of complex physical systems, making otherwise computationally prohibitive tasks—such as uncertainty quantification and optimal control—tractable. In particular, data-driven ROMs are opening new opportunities for developing powerful surrogate models that can learn from both numerical simulations and experimental data. However, integrating these sources of information is far from straightforward. Real-world measurements often involve multiple sources of uncertainty, including non-Gaussian noise, latent heterogeneities, unknown or hidden system parameters, and complex multiscale phenomena occurring below the resolution of the available data. Consequently, real-world data are more naturally described by parametrized stochastic partial differential equations (SPDEs), for which the development of effective ROMs remains relatively unexplored. In this setting, the system response is better characterized by a probability distribution than by a deterministic field, motivating the development of probabilistic surrogates. 
In this work, we introduce POD-CNF, a data-driven generative ROM for efficient uncertainty quantification and characterization of parametrized high-dimensional systems. The proposed approach combines projection-based model reduction with deep generative modeling. First, Proper Orthogonal Decomposition (POD) is employed to obtain a low-dimensional representation of the solution manifold. A Conditional Normalizing Flow (CNF) is then trained to approximate the conditional probability density of the reduced coefficients given the input parameters. This formulation allows the surrogate to capture complex and potentially multimodal output distributions, going beyond standard regression-based ROMs that typically focus on mean least-squares predictions. 
We assess the effectiveness of POD-CNF through a series of numerical experiments, demonstrating its ability to accurately reproduce output distributions while substantially reducing computational costs. Furthermore, a proof-of-concept application illustrates how the generative point of view provided by the normalizing flow facilitates the integration of POD-CNF into Bayesian inverse problems, highlighting promising directions for future research.
\end{abstract}




\newcommand{\nicolahere}{\begin{center}\;\\\textcolor{blue}{\textnormal{--- Nicola: \textit{On my side, all of the stuff above this line can be considered "almost" final}. ---}\;\\\;}\end{center}}
\newcommand{\guide}[1]{\noindent\textcolor{purple}{\textsf{[Nicola's proposal]\;\;}\emph{{#1}}}}
\newcommand{\bmu}{\boldsymbol{\mu}}
\newcommand{\btheta}{\boldsymbol{\theta}}
\newcommand{\bx}{\boldsymbol{x}}
\newcommand{\bz}{\boldsymbol{z}}
\newcommand{\wasserstein}{\mathcal{W}_1}
\newcommand{\ewass}{\tilde{\mathcal{W}}_1}
\newcommand{\pspace}{\mathcal{P}}
\newcommand{\prob}{\mathbb{P}}
\newcommand{\bu}{\mathbf{u}}
\newcommand{\espace}{\mathcal{E}}
\newcommand{\event}{\omega}
\newcommand{\bV}{\mathbf{V}}
\newcommand{\pdensity}{\xi}

\newcommand{\by}{\mathbf{y}}
\newcommand{\bb}{\mathbf{b}}
\newcommand{\bbf}{\mathbf{f}}
\newcommand{\bn}{\boldsymbol{n}} 
\newcommand{\bg}{\boldsymbol{g}}
\newcommand{\identity}{\mathbf{I}}
\newcommand{\indicator}{\mathbf{1}}
\newcommand{\bzero}{\mathbf{0}}
\newcommand{\normal}{\mathcal{N}}
\newcommand{\cspace}{\mathcal{C}}
\newcommand{\density}{\mathsf{p}}
\newcommand{\norm}[1]{\left\lVert#1\right\rVert}

\newcommand{\guideAndre}[1]{\noindent\textcolor{green}{\textsf{[Andrea's proposal]\;\;}\emph{{#1}}}}

\section{Introduction}\label{intro}

Partial Differential Equations (PDEs) constitute a standard mathematical framework for modeling physical phenomena in engineering and applied sciences \cite{quarteroni2017}. Given suitable boundary and initial conditions, their solutions can be approximated numerically through high-fidelity discretization techniques such as the Finite Element Method (FEM), Finite Volume Method (FVM), or spectral methods \cite{quarteroni2017}. In many practical applications, however, the governing equations depend on a set of input parameters describing material properties, geometric configurations, operating conditions, or boundary data. Consequently, one is not interested in solving a single PDE instance, but rather an entire family of parametrized problems. From an abstract point of view, these latter can be expressed, in the stationary case, as the solution of a parametrized PDE of the form
\begin{equation}
\label{eq:deterministic}
\mathcal{L}(u,\bmu) = f(\bmu)
\qquad \text{in}\;\Omega,
\end{equation}
complemented by suitable boundary conditions, where $\Omega \subset \mathbb{R}^d$ denotes a bounded spatial domain and $\bmu \in \mathbb{R}^p$ a vector of deterministic parameters. Here, $\mathcal{L}$ and $f$ denote a parametrized differential operator and source term, respectively. Under suitable assumptions, the problem is well posed, meaning that for every parameter value $\bmu$ there exists a unique solution field $u=u(\cdot,\bmu)$.

While high-fidelity numerical solvers provide accurate approximations of \eqref{eq:deterministic}, repeatedly solving a differential problem for many parameter configurations may become computationally prohibitive. This issue arises in several many-query scenarios, including uncertainty quantification, optimization, inverse problems, and real-time control. Such challenges have motivated the development of Reduced Order Modeling (ROM) techniques, whose objective is to construct low-dimensional surrogate models capable of rapidly approximating the solution manifold associated with a given parametrized PDE \cite{maday2002reduced, quarteroni2015reduced, hesthaven2022reduced, franco2023deep}.

Among ROM methodologies, projection-based techniques such as Proper Orthogonal Decomposition (POD) have proven particularly effective in extracting compact low-dimensional representations from high-fidelity simulations. More recently, advances in machine learning and the increasing availability of data have fostered the emergence of data-driven reduced-order models, where neural networks and statistical learning techniques are employed to approximate reduced dynamics or parameter-to-solution maps directly from data. These approaches are progressively opening the possibility of constructing surrogate models capable of learning not only from numerical simulations but also from experimental observations. 
However, the integration of experimental data into ROM frameworks is far from being straightforward. Classical reduced-order models are typically formulated under the assumption that the system response is uniquely determined by a finite-dimensional parameter vector. While this assumption is often appropriate in simulation-driven settings, it becomes restrictive when dealing with real-world measurements. In fact, experimental observations frequently contain variability that cannot be fully explained by the observable parameters alone. Such variability may arise from, e.g.,  measurement noise, unresolved multiscale phenomena, latent heterogeneities, unknown external forcings, or partially observed operating conditions, to mention a few relevant examples. 
As a consequence, two experiments performed under nominally identical observable conditions may yield significantly different system responses. 

In these settings, the deterministic parameter-to-solution paradigm underlying standard ROM methodologies becomes inadequate. Rather than learning a deterministic map
$
\bmu \mapsto u(\cdot,\bmu),
$
one should instead aim at characterizing a conditional probability distribution over possible system responses.
A natural mathematical framework for formalizing this uncertainty is provided by parametrized stochastic PDEs (SPDEs). Let $(\espace,\mathcal{F},\prob)$ be a probability space, and let $\event \in\espace$ represent a generic event in such space. The nonobservability of such events serves to model our partial knowledge of the system, thus accounting for, e.g., unknown sources of variability or unresolved scales. We then consider governing equations of the form
\begin{equation}
\label{eq:stochastic}
\tilde{\mathcal{L}}(u,\bmu,\event)
=
\tilde{f}(\bmu,\event)
\qquad \text{in}\;\Omega,
\end{equation}
where $\bmu$ represents observable deterministic parameters, while $\event$ models hidden stochastic effects. The solution field $u=u(\bx,\bmu,\event)$ is therefore no longer deterministic, but rather described by a random field conditioned on the observable parameters $\bmu$. The distinction between \eqref{eq:deterministic} and \eqref{eq:stochastic} is conceptually important. Deterministic parametrized PDEs can often be interpreted as simplified models of more complex stochastic physical phenomena. For example, in porous media flow, a deterministic nonlinear Darcy model of the form
$
-\nabla\cdot(u^{\mu-1}\nabla u)=f
$
may represent an idealized approximation, whereas a more realistic formulation accounting for unresolved material heterogeneities could read as 
$
-\nabla\cdot(a(\bx,\event)u^{\mu-1}\nabla u)=f,
$
where the random field $a(\bx,\cdot)$ models uncertain permeability variations and $\event$ is therefore unknown.

Within this stochastic setting, the objective of surrogate modeling changes fundamentally. Rather than approximating a deterministic parameter-to-solution map, one seeks to characterize the conditional probability distribution of the solution field, namely
$\bmu \mapsto \mathsf{p}(u|\bmu)$,
where $\mathsf{p}(u|\bmu)$ denotes the conditional law of $u$ given $\bmu$. Equivalently, $\mathsf{p}(u|\bmu)$ may be interpreted as the push-forward measure of $\prob$ through the map $\event \mapsto u(\cdot,\bmu,\event)$. After discretization, this corresponds to learning conditional probability distributions in very high-dimensional spaces, making the direct construction of probabilistic surrogates particularly challenging.

In this work, we propose tackling problems such as \eqref{eq:stochastic} by combing generative artificial intelligence (AI) and reduced order modeling. The former is meant to account for the stochasticity, while the latter can effectively handle the high dimensions into play. Our research is therefore part of a broader emerging trend concerning the development of probabilistic surrogates. We mention, for instance, the work by Coscia et al. \cite{coscia2024generative}, where generative adversarial neural networks are incorporated within ROM pipelines in order to quantify output uncertainties. There, however, no intrinstic stochasticity is present as part of the phenomenon under study: rather, all uncertainties emerge from the very definition of the ROM, essentially capturing a ROM-FOM closure in a probabilistic sense. More aligned with our setting are, instead, the recent works in \cite{bulte2024, molinaro2024generative, lim2025score, berman2026stochastic, kim2026stochastic}. These, in fact, explore the use of generative AI and machine learning for the construction of probabilistic neural operators and surrogate models, where stochasticity comes from partial knowledge of the system. Our work aims at pursuing a similar objective, however avoiding the use of overly complex architectures (such as diffusion models), in addition restoring a strong connection with the ROM literature. Indeed, it is now widely recognized that reduced order modeling and operator learning are closely interconnected, providing complementary perspectives with different advantages and trade-offs: we mention, for instance, the eloquent analogy between the POD-NN method \cite{hesthaven2018} and the so-called POD-DeepONet architecture \cite{lu2022comprehensive}, just to mention two widely used techniques appeared in the community in the last decade.

Before presenting our approach, it is important to clarify that, despite the similarity in terminology, the setting considered in this work is fundamentally different from that considered in several existing contributions. An example is provided by the \emph{stochastic reduced basis (RB) methods} \cite{nair2002stochastic}, and, more generally, by most applications of the RB method to parametric PDEs with stochastic coefficients: see, e.g.,  \cite{boyaval2010reduced}. There, in fact, stochasticity typically enters the governing equations in a fully specified manner, making the construction of the reduced-order model relatively straightforward. Similar approaches can be found, e.g., in \cite{chen2013weighted,chen2014comparison,chen2017reduced}. By contrast, our central assumption is that the occurring event $\event$, and more generally model \eqref{eq:stochastic}, are effectively unknown. This is meant to reflect the reality of many experimental settings, where the underlying physical mechanisms are only partially understood and the governing laws are not directly accessible. Our goal, in fact, is to lay the foundations for a unified framework that is sufficiently flexible to support learning from both computational models (FOMs) and empirical data. The objectives of the aforementioned works, instead, were intrinsically different.

Similarly, our definition of \emph{probabilistic reduced order model} is not to be confused with the one adopted by some other authors: see, e.g., \cite{billaud2024probabilistic}, where the word «probabilistic»  does not refer to the ROM per se, but rather relates to the algorithm underlying its construction, which the authors exploit to improve stability and generalization. In particular, these are works in which stochasticity is introduced only during the offline stage, but disappears at inference time.


\subsection{Our contribution}
In this work, we propose POD-CNF, a probabilistic reduced-order modeling framework for high-dimensional systems arising from parametrized stochastic PDEs. The proposed approach combines Proper Orthogonal Decomposition (POD) and Conditional Normalizing Flows (CNF). The former is employed to extract a low-dimensional representation of the solution manifold from a suitable set of high-fidelity snapshots (training data). The latter, instead, is trained over the same data to approximate the conditional probability distribution of the reduced coefficients given the input parameters $\bmu$. In this manner, the surrogate model learns a parameter-to-distribution map in terms of a set of reduced coordinates, enabling the generation of complex stochastic solution realizations while maintaining computational efficiency. From a theoretical standpoint, our construction is rigorously grounded on a suitable error decomposition formula, which explicitly decouples the inaccuracies introduced by the POD basis and those caused by the CNF module. This is obtained by combining standard results on the POD projection error with a novel proof about the universality of CNFs.
%

The remainder of the paper is organized as follows. In Section~\ref{sec:setting}, we introduce the mathematical setting and discuss the probabilistic formulation of reduced-order surrogate modeling for stochastic parametrized systems. Section~\ref{sec:method} presents the proposed POD-CNF methodology, including the projection-based reduction step and the conditional generative modeling framework, discussing its underlying theoretical rationale and its practical implementation. Numerical experiments assessing the accuracy and efficiency of the proposed approach are presented in Section~\ref{sec:exp}. In Section~\ref{sec:uq}, we demonstrate a proof-of-concept application to forward and inverse uncertainty quantification. Finally, conclusions and future research directions are discussed in Section~\ref{sec:conclusions}.

\section{Problem formulation}
\label{sec:setting}
Let $\Omega \subset \mathbb{R}^d$ be a bounded spatial domain. We consider a time-independent physical field $u$ whose behavior is governed by a parametrized stochastic PDE (SPDE), namely 
\begin{equation}
\label{eq:spde_formulation}
\tilde{\mathcal{L}}(u, \bmu, \event) = \tilde{f}(\bmu, \event) \qquad \text{in}\;\Omega,
\end{equation}
where $\bmu\in\pspace\subset\mathbb{R}^p$ is a suitable parameter vector, taking values in a given compact set $\pspace$, whereas $\event$ represents a stochastic event in some probability space $(\espace,\mathcal{F},\mathbb{P}).$ Here, $\tilde{\mathcal{L}}:V\times \pspace\times\espace\to\mathbb{R}$ is a map acting as a differential operator in the first variable, where $V\subset L^{2}(\Omega)$ is a  suitable function space ---e.g., $V=H^1(\Omega)$, or  $V=H^1_0(\Omega)$--- and the solution to \eqref{eq:spde_formulation} exists unique for each parametric instance $\bmu$ and each event $\event$. Similarly, $\tilde{f}$ represents a parametrized stochastic
source, incorporating external forces and boundary data.

We assume that the solutions to the SPDE can be either approximated numerically or measured experimentally over a suitable discretization of the spatial domain $\Omega$, independent of $\bmu$ and $\event$. 
Under both scenarios,
every parametric instance $\bmu\in\pspace$ can now be associated to a corresponding high-dimensional random vector, $\bu_{\bmu}\in\mathbb{R}^{N_h}$, 
where $N_h$ are the degrees of freedom associated with either the resolution of the data or the discretization of the numerical solver. In general, assuming a Lagrangian discretization in the second case, this translates into
\begin{equation*}\mathrm{u}^{(j)}_{\bmu}=u(\bx_j, \bmu, \cdot) + \epsilon(\bx_j,\bmu, \cdot),\qquad j=1,\dots,N_h,\end{equation*}
where $x_j$ is the $j$th dof in the discretization and the dependence on $\omega$ has been made implicit. Here, $\epsilon=\epsilon(x,\bmu,\omega)$ is a further source of uncertainty, potentially associated with, e.g., measurement noise or numerical errors.

We collect all sources of uncertainty within a suitable random vector of hidden latent variables, $\btheta$. Informally, $\btheta\cong[\omega,\epsilon]$. More precisely, $\btheta$ encapsulates all the variability in the discrete solution $\bu$ that cannot be explained in terms of the sole problem parameters $\bmu$. We will thus write $\bu_{\bmu}=\bu_{\bmu}(\btheta)$ whenever willing to emphasize such stochastic dependency. In the case of numerical simulations, the vector $\bu_{\bmu}$ can be characterized as the solution to a suitable (nonlinear) system of algebraic equations,
\begin{equation}
    \label{eq:discrete_spde}
    \mathbf{L}(\mathbf{u},\bmu,\btheta)=\mathbf{f}(\bmu,\btheta)
\end{equation}
obtained by discretizing \eqref{eq:spde_formulation} using, e.g., the FEM. With this notation, it is evident that $\bu$ can be seen as the solution to a parametric problem in the pair $(\bmu,\btheta)$. The issue is that $\btheta$ is effectively \emph{unknown}, or unobserved, which causes the solution to be stochastic. In general, the dimension of $\btheta$ is also unknown. Still, it can be upper bounded by $N_h$ using standard arguments from Optimal Transport theory (cf. Section \ref{subsec:surrogates}).

We emphasize that the presented framework is general enough to encompass many different scenarios. In particular, it can be used to address problems in which the stochasticity arises from either \emph{(i)} hidden physical parameters, \emph{(ii)}  numerical errors, \emph{(iii)} measurement uncertainty, or a combination of the three. In this work, we will focus exclusively on the first scenario, leaving interactions with the other research strands for future explorations. This choice is motivated by two main reasons. The first one, is that an exhaustive treatment of problems with hidden variables is a fundamental prerequisite for scaling the framework to challenging real-world settings. In fact, experimental settings are often characterized by partial knowledge of the system and its governing equations. Furthermore, the set of explanatory variables is often underestimated, meaning that, among all of the relevant variables influencing the system, only a few are effectively observed and recorded (here corresponding to the vector of parameters $\bmu$).

The second reason is that, in contrast to \emph{(i)}, the literature on \emph{(ii)} and \emph{(iii)} is already quite mature. For instance, the research on 
 \emph{Probabilistic Numerics} \cite{hennig2022probabilistic, cockayne2019bayesian} and \emph{Statistical Solutions} \cite{foias1973statistical, fjordholm2017} provide a solid foundation for addressing \emph{(ii)}. The former, in fact, allows to complement classical \emph{a priori} and \emph{a posteriori} estimates with rigorous uncertainty quantification on numerical errors. The latter, instead, provides a tailored framework for studying turbulent and chaotic systems, where even mild changes in the numerical solver setup can produce remarkably different results. 

Similarly, the broad literatures on \emph{Bayesian Statistics} \cite{gelman1995bayesian, cockayne2019bayesian} and \emph{Data Assimilation} \cite{law2015data} provide a well-established background for studying problems subject to measurement noise, as in \emph{(iii)}.
\\\\
In light of these considerations, the present work focuses on problems where stochasticity arises from hidden, unobserved, problem parameters. Precisely, our goal is to develop a computational pipeline for constructing a probabilistic surrogate capable of capturing the complexity of the parameter-to-random-solution map, $\bmu\mapsto\bu_{\bmu}$, associated with the problem under study. This idea is made rigorous in the next subsection.


\subsection{Probabilistic surrogates of parameter-to-random-solution maps}
\label{subsec:surrogates}
There are at least two equivalent ways of formulating our learning problem: one in terms of conditional probability distributions, and the other in terms of maps whose outputs are random variables. In this work, we adopt the latter perspective, which we find more convenient at inference time: see Remark~\ref{remark:2ways} for a deeper discussion on this matter.

In order to set our framework, we shall first introduce some notation. Let $\wasserstein$ denote the Wasserstein-1 distance. We recall that, given two probability measures $\mathbb{P}$ and $\mathbb{Q}$ over $\mathbb{R}^m$, the latter is defined as
\begin{equation}
    \label{eq:wasserstein}\wasserstein(\mathbb{P},\mathbb{Q}):=\inf_{\gamma\in\Gamma(\mathbb{P},\mathbb{Q})}\int_{\mathbb{R}^m\times\mathbb{R}^m}|x-y|\gamma(dxdy),
\end{equation}
where $\Gamma(\mathbb{P},\mathbb{Q})$ is the set of all probability measures over $\mathbb{R}^{m}\times\mathbb{R}^{m}$ having $\mathbb{P}$ and $\mathbb{Q}$ as marginals, respectively.
With little abuse of notation, given two random variables $X,Y\in\mathbb{R}^{m}$ with probability laws $\prob_X$ and $\prob_Y$, we write $\wasserstein(X,Y)$ to intend $\wasserstein(\prob_X, \prob_Y)$. Note that the random variables $X$ and $Y$ need not be defined over the same probability space for the discrepancy $\wasserstein(X,Y)$ to be well defined. Given any positive $q\in\mathbb{N}$, let $\mathcal{N}(0,\mathbf{I}_q)$ denote the standard $q$-dimensional Gaussian distribution. We recall that, by standard results in Optimal Transport theory, for every random variable $X\in\mathbb{R}^m$ there exists an index $1\le q\le m$ and a measurable map $T:\mathbb{R}^q\to\mathbb{R}^m$ such that, if $\bz\sim\mathcal{N}(0,\mathbf{I}_q)$, then
$X$ and $T(\bz)$ have the same law, meaning $\wasserstein(X,T(\bz))=0$: see, e.g., \cite[Chapter 1]{villani2009optimal}.
Typically, $q=m$. However, smaller values may at times be possible, for instance when the support of $X$ lies inside some submanifold of $\mathbb{R}^m.$
\\\\
With this notation, we define a probabilistic surrogate of $\bmu\mapsto\bu_{\bmu}$ to be any map of the form
\begin{equation}\label{eq:probabilistic-surrogate}
\Phi:\mathbb{R}^p\times\mathbb{R}^q\to\mathbb{R}^{N_h}\quad \text{with}\; 1\le q\le N_h,
\end{equation}
for which $\wasserstein(\bu_{\bmu},\;\Phi(\bmu,\bz))$ is sufficiently small, where $\bz\sim\mathcal{N}(0,\mathbf{I}_q)$ and $\bmu$ varies in $\pspace$. Here, $\bz$ is auxiliary random vector of latent variables, which is introduced to make the output of $\Phi=\Phi(\bmu,\bz)$ aleatoric. Intuitively, this construction aims at mimicking that of the ground truth model, $\bu=\bu_{\bmu}(\btheta)$. However, the vectors $\bz$ and $\btheta$ need not to be related one to another. What matters, is that they generate (approximately) the same stochasticity in the final output.
From a quantitative point of view, if the parameter space is equipped with a suitable density function $\pdensity=\pdensity(\bmu)$, the accuracy of a probabilistic surrogate can be measured via the Average Wasserstein Distance (AWD), namely
\begin{equation}
    \label{eq:awd}
\text{AWD}(\Phi):=\int_{\pspace}\wasserstein(\bu_{\bmu}, \Phi(\bmu, \bz))\pdensity(\bmu)d\bmu.
\end{equation}
At inference time, such probabilistic surrogates can be used to generate multiple \emph{plausible} solutions in real-time. For instance, given an operating condition of interest $\bmu_0\in\pspace$, one can produce multiple independent random solutions $\tilde{\bu}_{\bmu_{0},1},\dots,\tilde{\bu}_{\bmu_0,m_{\text{gen}}}$ by sampling $\bz_1,\dots,\bz_{m_{\text{gen}}}\sim\mathcal{N}(0,\mathbf{I}_q)$ and computing $\tilde{\bu}_{\bmu_0,k}=\Phi(\bmu_0,\bz_{k})$. These realizations can then be used to approximate the conditional law of $\bu_{\bmu}$ and to estimate relevant quantities of interest. If the samples are very similar to each other then one can conclude that, under $\bmu_0$, the effect of the hidden parameters is nearly negligible. In contrast, large deviations or multi-modal behaviors can be taken as evidence of intrinsic uncertainties affecting the system.


\begin{remark}[Learning generators vs learning distributions] 
\label{remark:2ways} Since we are only interested in the conditional law of the solution field $\bu_{\bmu}$, and not in the random variable \emph{per se}, an alternative approach is to define a probabilistic surrogate as a (nonlinear) operator mapping parameters onto probability distributions, $\pspace\mapsto \mathscr{P}(\mathbb{R}^{N_h}).$ However, since $N_h$ can be fairly large, constructing an object that returns explicit probability distributions over $\mathbb{R}^{N_h}$ can be highly nontrivial. Furthermore, at inference time, using such a tool can pose significant challenges. Consider for instance the problem of employing the surrogate to estimate a given statistics of the solution, e.g., $\mathbb{E}[Q(\bu_{\bmu_0})]$ for some $Q:\mathbb{R}^{N_h}\to\mathbb{R}$ and $\bmu_0\in\pspace$. In this case, the approximation would read \begin{equation}
    \mathbb{E}[Q(\bu_{\bmu_0})]\approx \int_{\mathbb{R}^{N_h}}Q(\bu)\hat{\mathbb{P}}_{\bmu_0}(d\bu),
\end{equation} where $\hat{\mathbb{P}}_{\bmu_0}\in\mathscr{P}(\mathbb{R}^{N_h})$ is the probability law predicted by the surrogate. However, even if $\hat{\mathbb{P}}_{\bmu_0}$ were available in closed form, for example through its density, evaluating the above integral would remain computationally challenging: numerical quadrature methods do not scale sufficiently well to be practical in  $\mathbb{R}^{N_h}$, while Monte Carlo methods can easily become troublesome in high-dimensions, even when the underlying density is known explicitly: see, e.g., Markov-Chain Monte Carlo samplers and similar tools \cite[Chapters 10-12]{gelman1995bayesian}.

In contrast, learning a probabilistic surrogate that operates as a generative model greatly simplifies the problem. In fact, the generative model directly acts as an approximate conditional sampler, allowing one to estimate expectations via Monte Carlo sampling, e.g.,
\begin{equation}
\label{eq:conditional-monte-carlo}
\mathbb{E}[Q(\bu_{\bmu_0})]\approx\frac{1}{m_{\text{gen}}}\sum_{k=1}^{m_{\text{gen}}}Q(\Phi(\bmu_0,\bz_k)).\end{equation}
In general, the accuracy of the above will depend on the quality of the surrogate and on the number of Monte Carlo samples. However, we notice that this second source of inaccuracy can be made arbitrarily small at a negligible computational cost: in fact, the Monte Carlo error scales as $O(m_{\text{gen}}^{-1/2})$, but because the surrogate is computationally cheap to evaluate, $m_{\text{gen}}$ can indeed be taken as large as required.    
\end{remark}

\section{The POD-CNF approach}
\label{sec:method}

We assume to have access to a dataset   $\{\bmu_i, \bu_{i}\}_{i=1}^{N}$ consisting of different parametric configurations and corresponding random solutions. These may either consist of numerical simulations, obtained by means of a given stochastic FOM, or derive from suitable real-world experiments. We stress that, in both cases, the system under study is assumed to be \emph{time-independent}, but potentially subject to complex physical laws that entail nonlinearities.
We also remark that $\bu_i=\bu_{\bmu_i}(\btheta_i)$, meaning that the $i$-th solution in the dataset corresponds to the $i$-th parametric configuration ($\bmu_i$, which is known) and the $i$-th stochastic event ($\btheta_i$, which is unknown). Our goal is to leverage this data in order to construct a surrogate $\Phi$ of the parameter-to-random-solution map $\bmu\mapsto\bu_{\bmu}.$ To this end, we propose the following ansatz,
\begin{equation}
    \label{eq:ansatz}
    \Phi(\bmu,\bz):=\mathbf{V}\phi(\bmu,\bz)
\end{equation}
which is based on a low-rank representation of the solution field. More precisely, for a suitable reduced dimension $n \ll N_h$, we assume that:
\begin{itemize}
    \item $\mathbf{V}\in\mathbb{R}^{N_h\times n}$ is a deterministic orthonormal matrix, whose purpose is to reduce the dimensionality of the problem. The idea, in fact, is that although the solutions $\bu$ live in a high dimensional space, $\mathbb{R}^{N_h},$ there may exist a much smaller space, $\mathcal{V}=\text{span}(\mathbf{V})\cong\mathbb{R}^n$, in which the solutions approximately lie. Informally, if we denote by $\text{supp}(\bu_{\bmu})$ the support of the random solution $\bu_{\bmu}$, this corresponds to saying that $\cup_{\bmu\in\pspace}\text{supp}(\bu_{\bmu})$ is approximately contained in $\mathcal{V}$. From a computational standpoint, if the parameter space is endowed with a suitable density $\pdensity=\pdensity(\bmu)$, such degree of approximation can be quantified as
    \begin{equation}
        \label{eq:projerror}
        \int_{\pspace}\mathbb{E}\|\bu_{\bmu}-\mathbf{V}\mathbf{V}^\top\mathbf{u}_{\bmu}\|\pdensity(\bmu)d\bmu,
    \end{equation}
    i.e. by the average expected projection error; 
    \item $\phi:\mathbb{R}^{p}\times\mathbb{R}^{n}\to\mathbb{R}^{n}$ is a latent generative model -- corresponding to $q=n$ in Eq. \eqref{eq:probabilistic-surrogate}. The latter aims at capturing both the variability induced by the problem parameters, $\bmu$, as well as the one produced by the hidden stochastic variables, $\btheta$. Intuitively, one aims at constructing $\phi$ in such a way that, for any given $\bmu\in\pspace$, the distribution of the generative model, $\phi(\bmu,\bz)$ with $\bz\sim\mathcal{N}(0,\mathbf{I}_n)$, is close to that of the random projection coefficients, $\mathbf{V}^\top\bu_{\bmu}$. In analogy to Eq. \eqref{eq:awd}, the quality of such approximation can be measured using, e.g., the AWD metric.
\end{itemize}

From a practical standpoint, both the construction of the matrix $\mathbf{V}$ and that of the latent generative model $\phi$ will be based on the data at hand. Specifically, in this work we propose relying on the Principal Orthogonal Decomposition (POD) algorithm and on Conditional Normalizing Flow (CNF) architectures in order to construct $\mathbf{V}$ and $\phi$, respectively, reason for which we term our approach POD-CNF. 

We mention that, in spirit, our proposal can also be understood as a stochastic generalization of the already well-established POD-NN method \cite{hesthaven2018}. The latter, in fact, is a data-driven ROM strategy, developed for deterministic problems, in which the solution is approximated as $\mathbf{u}_{\bmu}\approx \mathbf{V}\tilde{\phi}(\bmu)$, with $\mathbf{V}$ the POD matrix and $\tilde{\phi}:\mathbb{R}^p\to\mathbb{R}^n$ a standard multilayer perceptron (MLP) network trained to approximate the parameter-to-coefficient map $\bmu\mapsto\mathbf{V}^\top\mathbf{u}_{\bmu}$. Since its introduction, the POD-NN method has had a significant impact on the literature, with several subsequent approaches building upon or tracing their origins back to this pioneering idea. For instance, the DeepONet architecture \cite{lu2022comprehensive} can be seen as a continuous generalization of the POD-NN; similarly, alternative designs and implementations have led to the development of other prominent strategies such as POD-DL-ROM \cite{fresca2022}, Conservative POD-NN \cite{boon2025deep}, DOD-NN \cite{franco2026deep, kotowski2026new} and many more \cite{wang2019non, pichi2023artificial}. In this sense, the POD-NN approach can be regarded as the forerunner of a much wider class of methods, which our contribution aims to further extend.


The remainder of this Section is devoted to establishing suitable theoretical results supporting our construction (\S\ref{sec:theory}), as well as to presenting its practical implementation from an algorithmic standpoint (\S\ref{sec:practice}).

\subsection{Theoretical rationale}
\label{sec:theory}

There are several reasons for using POD and CNFs. To start, POD is a powerful and well-established algorithm, which is known to enjoy remarkable optimality properties. Furthermore, its implementation is relatively straightforward, as it boils down to a truncated Singular Value Decomposition (SVD) of the snapshots matrix $\mathbf{S}=[\bu_1,\dots,\bu_N]\in\mathbb{R}^{N_h\times N}$. 
On top of this, POD allows for the construction of \textit{stable} ROMs, in that small errors in the reduced space are propagate isometrically to the full order space via the map $\mathbb{R}^n\ni\mathbf{c}\mapsto \mathbf{V}\mathbf{c}\in\mathbb{R}^{N_h}.$ Finally, by construction, POD tends to filter out high frequency modes, which can be helpful in removing external white noise-type fluctuations in the data. 

The use of CNFs, instead, is mainly motivated by their expressivity and their computational efficiency. Indeed, under fairly mild assumptions, normalizing flows are known to be universal approximators \cite{koehler2021representational}, which makes them good candidates for learning the latent generative model $\bmu\mapsto \mathbf{V}^\top\bu_{\mu}.$ At the same time, they are also: \emph{(i)} computationally much less expensive to train than modern score-based generative models, such as Diffusion Models \cite{rombach2022high}, \emph{(ii)} far more efficient than other flow-based architectures, such as Flow Matching \cite{lipman2023flow, albergo2025stochastic} or Stochastic Interpolants \cite{albergo2025stochastic, mucke2026physics}, when employed at inference time, which makes them better suited for many query applications.

By combining the universality of CNFs with classical results on POD truncation errors, it is then relatively straightforward to derive corresponding theoretical insights for POD-CNF. This connection is formalized in the next two paragraphs, where we state several foundational results along this direction. Right after that, a brief discussion on the limitations of our approach is also provided. In order to keep the paper self-contained, technical proofs are deferred to the Appendix~\ref{appendix:proofs}.


\paragraph{Proper Orthogonal Decomposition (POD): idea and motivation.}
Proper Orthogonal Decomposition (POD) is a powerful data-driven technique which commonly used for dimensionality reduction and reduced order modeling in general. Its primary objective is to extract the most dominant characteristics (often called "modes") from a high-dimensional dataset. By identifying these modes, POD allows for a low-dimensional approximation that captures the essential behavior of the original system while drastically reducing its complexity. It is best suited in those frameworks in which \emph{linear} reduction is effective, cf. \cite{quarteroni2015reduced,cohen2021linear}.

In the simplest scenario, the implementation of the POD consists of computing an SVD of the  snapshots matrix $\mathbf{S}=[\bu_1,\dots,\bu_N]\in\mathbb{R}^{N_h\times N}$, namely
%
\begin{equation*}\mathbf{S}=\mathbf{U}\mathbf{\Sigma}\mathbf{W}^\top,\end{equation*}
where $\mathbf{U} \in \mathbb{R}^{N_h \times N_h}$ and $\mathbf{W} \in \mathbb{R}^{N \times N}$ are orthogonal matrices containing the left and right singular vectors, respectively, while $\Sigma$ is a diagonal matrix containing the singular values, $\sigma_1 \ge \sigma_2 \ge \dotsi \ge 0$, sorted in decreasing order. The POD matrix $\mathbf{V}\in\mathbb{R}^{N_h \times n}$ is then computed by retaining the first $n$ columns of $\mathbf{U}$.
It is well-known that this procedure is \emph{empirically optimal}, in that
\begin{equation}
    \label{eq:pod-energy}
    \sum_{i=1}^{N}\|\bu_i-\mathbf{V}\mathbf{V}^\top\bu_i\|^2=\min_{\mathbf{Q}\in\mathbb{R}^{N_h\times n}}\sum_{i=1}^{N}\|\bu_i-\mathbf{Q}\mathbf{Q}^\top\bu_i\|^2=\sum_{i>n}\sigma_i^2,
\end{equation}
see, e.g., \cite{quarteroni2015reduced}. In other words, the POD basis is optimal in the mean-square sense, minimizing the reconstruction error of the computed snapshots among all subspaces of a prescribed dimension. If the parameter instances $\bmu_i\in\mathcal{P}$ are drawn independently according to $f$, and the events $\btheta_i$ are independent as well, then, in the infinite data limit, $\mathbf{V}$ is known to converge to the \emph{ideal} POD matrix $\mathbf{V}_\infty$, solution to the minimization problem below
\begin{equation}
    \min_{\mathbf{Q}\in\mathbb{R}^{N_h\times n}}\int_{\mathcal{P}}\mathbb{E}\|\bu_{\bmu}-\mathbf{Q}\mathbf{Q}^\top\bu_{\bmu}\|^2\pdensity(\bmu)d\bmu,
\end{equation}
see, e.g., \cite[Section 2.2]{brivio2024pod}.  Theoretically, much like $\mathbf{V}$ is built by computing an SVD of $\mathbf{S}$, the ideal POD matrix can be computed by considering an SVD of the Hilbert-Schmidt operator $\mathscr{S}:L^2(\mathcal{P}\times\mathbb{R}^{n})\ni g\mapsto \int_{\mathcal{P}}\mathbb{E}[\bu_{\bmu}g]\pdensity(\bmu)d\bmu\in\mathbb{R}^{N_h}.$ Equivalently, $\mathbf{V}_{\infty}$ can also be derived by diagonalizing the uncentered covariance matrix, \begin{equation}
\label{eq:uncentered-covariance}
\mathbf{C}:=\int_{\pspace}\mathbb{E}[\bu_{\bmu}\otimes \bu_{\bmu}]\pdensity(\bmu)d\bmu\in\mathbb{R}^{N_h\times N_h},\end{equation}
and then collecting the $n$ eigenvectors associated with the $n$ largest eigenvalues $\lambda_i$. 

In general, if $\sigma_{i,\infty}$ are the singular values of $\mathscr{S}$ in the infinite data-limit, and $\sigma_{i,N}$ are those of $\mathbf{S}$ obtained with $N$ randomly drawn snapshots, then one has $$\lim_{N\to+\infty}\frac{1}{N}\sigma_{i,N}^2=\sigma_{i,\infty}^2=\lambda_i,$$
almost surely: see, for instance, \cite[Section 2.2]{brivio2024pod} and \cite[Chapter 3]{hesthaven2015certified}.
The whole machinery can also be generalized very easily to account for other inner-product norms, different from the plain Euclidean one, cf. \cite{quarteroni2015reduced, hesthaven2015certified}.
%

In our setting, the use of a projection-based dimensionality reduction technique is motivated by the following \emph{decoupling lemma}, which shows how the error a generative ROM built according to our ansatz can be split into two contributions: {\em (i)} a projection error, and {\em(ii)} a latent generative error.

\begin{lemma}
\label{lemma:decomposition}
Let $\{\bu_{\bmu}\}_{\bmu\in\pspace}$ be a parametric family of random variables in $\mathbb{R}^{N_h}$, with $\pspace\subseteq\mathbb{R}^{p}$, and let $\bV\in\mathbb{R}^{N_h\times n}$ be an orthonormal matrix. Let $\phi:\mathbb{R}^p\times\mathbb{R}^{n}\to\mathbb{R}^n$ be given, and let $\bz\in\mathbb{R}^n$ be a random vector with law $\bz\sim\mathbb{Q}$. 
Then,
\begin{equation}
\label{eq:main inequality}
\wasserstein(\bu_{\bmu}, \;\bV\phi(\bmu,\bz)) \leq \mathbb{E}\| \bu_{\bmu} - \bV\bV^\top \bu_{\bmu}\| + \wasserstein(\bV^\top \bu_{\bmu}, \;\phi(\bmu, \bz))\qquad\forall\bmu\in\pspace.
\end{equation}
\end{lemma}
\begin{proof}
    See Appendix~\ref{appendix:proofs}.
\end{proof}

Combining the result above with the aforementioned optimality of POD immediately yields the following corollary.

\begin{corollary}
\label{corollary:upperbound}
Let $\{\bu_{\bmu}\}_{\bmu\in\pspace}$ be a parametric family of random variables in $\mathbb{R}^{N_h}$, with $\pspace\subseteq\mathbb{R}^{p}$. Let $f$ be a probability density over $\pspace.$ Denote by $\mathbf{C}$ the $N_h\times N_h$ positive semi-definite matrix in \eqref{eq:uncentered-covariance}. Let $\lambda_1\ge\dots\ge\lambda_{N_h}\ge0$ be the eigenvalues of $\mathbf{C}$, listed in decreasing order, and let $\mathbf{v}_1,\dots,\mathbf{v}_{N_h}$ be the corresponding orthonormal eigenvectors. 
Given any truncation index $1\le n\le N_h$, assemble the orthonormal matrix $\bV=[\mathbf{v}_1,\dots,\mathbf{v}_{n}]\in\mathbb{R}^{N_h\times n}.$
For any $\phi:\mathbb{R}^p\times\mathbb{R}^n\to\mathbb{R}^n$ and any $n$-dimensional random vector $\bz$, one has
\begin{equation}
\label{eq:corollary}
\int_{\pspace}\wasserstein(\bu_{\bmu}, \;\bV\phi(\bmu,\bz))\pdensity(\bmu)d\bmu \leq \left[\sum_{j=n+1}^{N_h}\lambda_j\right]^{1/2} + \int_{\pspace}\wasserstein(\bV^\top \bu_{\bmu}, \;\phi(\bmu, \bz))\pdensity(\bmu)d\bmu.
\end{equation}
\end{corollary}
\begin{proof}
    See Appendix~\ref{appendix:proofs}.
\end{proof}

Basically, Corollary~\ref{corollary:upperbound} ensures that if the eigenvalues $\lambda_i$ decay sufficiently fast, then one can construct an effective probabilistic ROM by doing a POD and transferring all of the stochasticity at the latent level, in analogy to what is commonly done in the deterministic case.

\paragraph{Conditional Normalizing Flows (CNF): idea and motivation.}

Conditional Normalizing Flows (CNFs) are a class of generative models that aim to approximate conditional densities by leveraging optimal transport theory and the change of variables formula. Simply put, in order to approximate the conditional density of a given random variable, $\by|\bmu$, CNFs progressively transform Gaussian noise into the desired distribution by applying a sequence of $\bmu$-dependent transformations, $T_{\bmu,1},\dots, T_{\bmu,K}$. Such transformations are designed to be invertible in their main argument. This allows making the density of the transformed vector, $(T_{\bmu,K}\circ\dots\circ T_{\bmu,1})(\bz)$, fully explicit; in fact, thanks to the change of variables formula, one has
\begin{equation}
\label{eq:change-of-variables}
f_{(T_{\bmu,K}\circ\dots\circ T_{\bmu,1})(\bz)}(\by)=f_{\bz}((T_{\bmu,1}^{-1}\circ\dots\circ T_{\bmu,K}^{-1})(\by))\cdot\prod_{k=1}^K\left|\text{det}\;J_{T_{\bmu,k}}^{-1} ((T_{\bmu,k}^{-1}\circ\dots\circ T_{\bmu,K}^{-1})(\by))\right|,\end{equation}
with $f_{\bz}$ the density of the standard normal distribution in $\mathbb{R}^n.$

From a practical perspective, the transformations $T_{\bmu,1},\dots,T_{\bmu,K}$ constitute the trainable components of the CNF architecture, and they are usually designed by cleverly combining different units based on multi-layer perceptrons (MLPs). Such design is crucial, as it prescribes the expressivity of the whole architecture, but it also entails some nontrivial challenges. In fact, in order to make \eqref{eq:change-of-variables} practical, one must design the maps $T_{\bmu, k}$ so that they are easily invertible and so that Jacobian determinants can be computed efficiently.

In this respect, an approach which has proven successful, and which we adopt in this work, is to construct the transformations $T_{\bmu,k}$ by relying on \emph{affine coupling layers}. The latter were first introduced by Dinh et al. in \cite{dinh2016density} as part of a normalizing flow architecture thereon referred to as \emph{RealNVP}. Here, we shall consider a conditional variant of such architecture, which can be described as follows. For every $k\in\{1,\dots,K\}$, let $\bz_{k-1}$ denote a generic input for the $k$-th layer. 
%
%
%
The vector is split into two halves, $\bz_{k-1,A_k}$ and $\bz_{k-1,B_k}$, with $A_k\cup B_k=\{1,\dots,n\}$ being a suitable partition of the vector indices, $A_k\cap B_k=\emptyset$. Such partition is typically chosen at random and varies from layer to layer. If $n$ is even, a standard choice is to set $A_k=\{1,\dots,n/2\}$ and $B_k=\{n/2+1,\dots,n\}$ when $k$ is even, and $A_k=\{n/2+1,\dots,n\}$ and $B_k=\{1,\dots,n/2\}$ when $k$ is odd. In this case, the RealNVP is said to consist of \emph{alternating affine coupling layers}.

With this convention, the action of the $k$th layer $T_{\bmu,k}$ reads
%
\begin{equation}
\label{eq:couplinglayer}
T_{\bmu,k}(\bz_{k-1})=
\begin{cases}
    \bz_{k,A_k}=\bz_{k-1,A_k} \\
    \bz_{k,B_k}=\bz_{k-1,B_k} \odot \exp(s_k(\bz_{k-1,A_k},\bmu)) + t_k(\bz_{k-1,A_k}, \bmu)
\end{cases}
\end{equation}
where $s_k,t_k:\mathbb{R}^{|A_k|+p}\to\mathbb{R}^{|B_k|}$ are trainable MLPs, referred to as the \emph{scale} and \emph{translation} parameter, respectively, whereas $\odot$ denotes the Hadamard product. The full CNF architecture is then given by
\begin{equation}
    \phi(\bmu,\bz):=(T_{\bmu, K}\circ\dots\circ T_{\bmu,1})(\bz).
\end{equation}
We will say that $\phi$ is a ReLU CNF with alternating affine coupling layers if: \emph{(i)} the transformations $T_{\bmu,k}$ are defined as in \eqref{eq:couplinglayer}, with $A_k$ and $B_k$ given by the canonical splitting; \emph{(ii)} the MLPs $s_k$ and $t_k$ use the ReLU activation in all of their hidden layers; \emph{(iii)} the translation MLPs, $t_k$, have no terminal activation; \emph{(iv)} the scale MLPs, $s_k$, are equipped with a terminal continuous activation of choice, $\rho_k:\mathbb{R}\to\mathbb{R}$. For numerical stability, such activations are often chosen to be bounded: in our experiments, for instance, we use the hyperbolic tangent activation for all $k=1,\dots,K.$ Still, when discussing theoretical results we will always assume that $-\infty\in\overline{\rho_k(\mathbb{R})}$, thus preventing boundness from below. This is to guarantee the universality of the resulting CNF architecture under fairly general conditions: we refer the interested reader to Remark~\ref{remark:universality} in the Appendix.

The peculiar construction in \eqref{eq:couplinglayer}, where each layer leaves one half of the input vector unchanged, while scaling and translating the other half, is what makes RealNVPs particularly effective in practice. This, in fact, guarantees that each $T_{\bmu,k}$, is invertible ---regardless of the learned MLPs $s_k$ and $t_k$--- and that the Jacobian of $T_{\bmu,k}^{-1}$ is triangular, allowing for a fast computation of its determinant. At the same time, this design is general enough to make RealNVPs universal approximators: see, e.g., \cite{koehler2021representational} for results in the unconditional case, or \cite{lyu2022cflows} and Lemma~\ref{lemma:universality} in Appendix~\ref{appendix:universality} for suitable extensions to parametric problems.
\\\\
For all such reasons, we propose relying on CNFs in order to learn the latent generative model in  Eq.~\eqref{eq:ansatz}. With this choice, the aforementioned universality property, together with Corollary~\ref{corollary:upperbound}, immediately yields the following result, which provides the theoretical foundation for the POD-CNF approach.

\begin{theorem}
\label{theorem:universality}
Let $\bu_{\bmu}$ be a parametrized random vector in $\mathbb{R}^{N_h}$, where $\bmu\in\pspace\subset\mathbb{R}^p$ varies over a suitable parameter space, with $\pspace$ compact. Let $f$ be a probability density defined over $\pspace$. Assume that the law of $\bu_{\bmu}$ varies continuously in the $\wasserstein$ metric with respect to $\bmu$, and that there exists some $R>0$ such that $\|\bu_{\bmu}\|\le R$ almost surely for all $\bmu\in\pspace.$
Given any truncation index $n\ge q$, let $\bV$ be as in Corollary~\ref{corollary:upperbound}. 
Let $\prob_{n,\bmu}$ be the probability distribution of the  reduced random vector $\bV^\top \bu_{\bmu}$. Assume that, for every $\bmu\in\pspace$, the distribution $\prob_{n,\bmu}$ is
absolutely continuous.
Then, for every $\epsilon>0$, there exists a ReLU conditional normalizing flow $\phi_\epsilon$ with alternating affine coupling layers such that
\begin{equation*}\int_{\pspace}\wasserstein(\bu_{\bmu}, \;\bV\phi_\epsilon(\bmu,\bz))\pdensity(\bmu)d\bmu \leq \left[\sum_{j=n+1}^{N_h}\lambda_j\right]^{1/2} + \epsilon.\end{equation*}
\end{theorem}

\begin{proof}
See Appendix~\ref{appendix:proofs}.
\end{proof}

\paragraph{Intrinsic limitations of the POD-CNF approach.}
While the POD-CNF framework provides a principled and computationally efficient strategy for learning parametric solution distributions, it also relies on a number of assumptions that may limit its applicability in challenging regimes.

First, the effectiveness of the approach is intrinsically tied to the compressibility of the solution manifold by a linear subspace. In particular, the quality of the POD representation depends on the decay of the eigenvalues of the correlation operator. This is closely related to the behavior of the Kolmogorov $n$-width of the \emph{supporting} solution set 
$\mathcal{S}=\cup_{\bmu\in\pspace}\text{supp}(\bu_{\bmu})=\cup_{\bmu\in\pspace}\{\mathbf{v}\in\mathbb{R}^{N_h}\;:\;\forall\delta>0\;\;\mathbb{P}(|\bu_{\bmu}-\mathbf{v}|<\delta)>0)\}$: when this quantity decays slowly, a large number of POD modes may be required to achieve a satisfactory approximation accuracy. In such cases, the output dimension of the CNF increases, hindering the computational advantages of the reduced representation and potentially making the training of the generative model more challenging.

A related limitation is that POD is a linear dimensionality reduction technique and, by construction, preferentially retains energetic modes while filtering out small-scale and high-frequency components. Although this property can be beneficial when the discarded fluctuations are mainly associated with noise, it may become detrimental for problems where localized phenomena, sharp interfaces, or fine-scale structures play a relevant role in the solution dynamics. In such settings, nonlinear dimensionality reduction techniques, such as autoencoder-based representations, may provide more expressive alternatives and constitute a natural direction for future investigation.

A second class of limitations concerns the choice of CNFs as generative models. Although such models provide an attractive combination of expressivity and computational efficiency, their construction is constrained by invertibility assumption. For instance, architectures based on affine coupling transformations, such as the RealNVP model adopted in this work, achieve tractable likelihood evaluation at the expense of restricting the class of admissible transformations. More flexible alternatives, including continuous normalizing flows and flow-matching approaches, alleviate some of these architectural constraints, but typically require the numerical solution of an ordinary differential equation during inference, which can become computationally demanding in multi-query applications.

Finally, the scalability of CNFs with respect to their ambient dimension remains an important practical consideration. Even after POD truncation, the reduced random variables may exhibit highly anisotropic or nearly singular distributions, for instance when the data are concentrated around lower-dimensional structures embedded in the latent space. Since standard CNFs are constructed by transporting a full-dimensional reference distribution, accurately representing such distributions may require increasingly complex transformations and can lead to optimization difficulties. Developing dimensionality-adaptive generative models or combining POD with nonlinear manifold-learning techniques may provide effective strategies to address these challenges, representing an interesting direction for future research.

\subsection{Practical implementation}
\label{sec:practice}
From an operational standpoint, we construct a POD-CNF surrogate through a three-step \emph{offline} procedure, whose workflow is further summarized in Algorithm~\ref{algo:pod-cnf}.

\paragraph{Step 1: data collection and splitting.}
To start, we first sample $S$ parametric configurations $\bmu_1,\dots,\bmu_S\in\pspace$ by drawing independent realizations from the sampling density $\xi$. For each $\bmu_i$, we then run the FOM (resp., the lab experiment) to obtain a random solution $\bu_i=\bu_{\bmu_i}(\btheta_i).$ We then split the dataset $$\{\bmu_i,\bu_i\}_{i=1}^{S}=\{\bmu_i,\bu_i\}_{i=1}^{N}\cup\{\bmu_i^{\text{valid}},\bu_i^{\text{valid}}\}_{i=1}^{N_{\text{valid}}}\cup\{\bmu_i^{\text{test}},\bu_i^{\text{test}}\}_{i=1}^{N_{\text{test}}},$$ into \emph{training}, \emph{validation}, and \emph{test} sets, generally following a 75:20:5 ratio.

\paragraph{Step 2: POD basis extraction.}
We leverage the training data to assemble the snapshots matrix $\mathbf{S}=[\bu_{1},\dots,\bu_{N}]$ and compute its SVD. In light of Eq. \eqref{eq:pod-energy}, we perform POD by truncating the SVD to the first \emph{even} index $n\in 2\mathbb{N}$ for which the relative residual energy,
$$\left[\sum_{i>n}\sigma_i^2\right]/\left[\sum_{i=1}^{\min{N_h,N}}\sigma_i^2\right],$$
is below a given threshold $\tau\in(0,1).$ We require $n$ to be an even number in order to standardize the subsequent construction of the CNF module via alternating coupling layers, cf. \S\ref{sec:theory}.


\paragraph{Step 3: CNF implementation and training.}
To construct the CNF module, we start by choosing a depth $K$ and designing the architectures of the trainable MLPs in the different coupling layers. Let us denote by $\boldsymbol{\eta}\in\mathbb{R}^{D}$ the vector collecting all the trainable parameters appearing in the definition of the CNF. The trainable architecture can be then regarded as a map $\phi_{\boldsymbol{\eta}}:\mathbb{R}^{p}\times\mathbb{R}^n\to\mathbb{R}^n$, where $\boldsymbol{\eta}$ is to be optimized.

To train the CNF, i.e., to estimate the optimal value for $\boldsymbol{\eta}$, we first assemble the reduced dataset $\{\bmu_i,\mathbf{c}_i\}_{i=1}^{N}$ where $\mathbf{c}_i:=\mathbf{V}^\top\bu_i$ are the POD-coefficients of the training snapshots. 
The model is then trained by minimizing the \emph{negative log-likelihood} (NLL),
\begin{equation}
    \label{eq:nll}
    \textsf{NLL}(\boldsymbol{\eta}):=-\frac{1}{N}\sum_{i=1}^{N}\log \hat{f}_{\boldsymbol{\eta},\bmu_i}(\mathbf{c}_i)
\end{equation}
where $\hat{f}_{\boldsymbol{\eta},\bmu}$ denotes the density of $\phi_{\boldsymbol{\eta}}(\bmu,\bz)$ with $\bz\sim\mathcal{N}(\mathbf{0},\mathbf{I}_n).$ If we denote by $T_{\bmu,k}^{\boldsymbol{\eta}}$ the $k$th coupling layer in $\phi_{\boldsymbol{\eta}}$, using Eq. \eqref{eq:change-of-variables}, we may rewrite Eq. \eqref{eq:nll} as
\begin{equation}
    \label{eq:nll2}
    \textsf{NLL}(\boldsymbol{\eta})=\frac{1}{N}\sum_{i=1}^{N}\left[\frac{1}{2}\|[\phi_{\boldsymbol{\eta}}(\bmu_i,\cdot)]^{-1}(\mathbf{V}^\top\bu_i)\|^2 - \sum_{k=1}^{K}\log\left|\text{det}\;J_{T_{\bmu_i,k}}^{-1} ((T_{\bmu_i,K}^{\boldsymbol{\eta}}\circ\dots\circ T_{\bmu_i,k}^{\boldsymbol{\eta}})^{-1}(\mathbf{V}^\top\bu_i))\right|\right]+\text{constant}.
\end{equation}
Minimizing (numerically) the above function with respect to $\boldsymbol{\eta}$ yields the trained architecture $\phi=\phi_{\boldsymbol{\eta}_*}.$

We mention that, in the literature, this is the standard approach for training normalizing flows and CNFs in general \cite{papamakarios2021normalizing}. From a theoretical standpoint, there are several reasons that make NLL minimization a natural choice. One of such reasons is that, up to a Monte Carlo approximation, it can be shown that minimizing \eqref{eq:nll} is equivalent to minimizing the expected Kullback-Leibler divergence (EKL) between the true conditional density of $\mathbf{V}^\top\bu_{\bmu}$, denoted as $f_{\bmu}$, and the conditional density predicted by the CNF module, that is
\begin{equation}
\label{eq:KL}
\textsf{EKL}(\boldsymbol{\eta})=\int_{\pspace}D_{\text{KL}}(f_{\bmu},\hat{f}_{\bmu,\boldsymbol{\eta}})\xi(\bmu)d\bmu.\end{equation}
We refer the interested reader to, e.g., \cite[Chapter 5.5]{Goodfellow-et-al-2016} for a deeper discussion on such equivalence. In our framework, this correspondence becomes particularly valuable as it allows us to connect the CNF training with our theoretical analysis in \S\ref{sec:theory}. In fact, if we assume the two densities in Eq. \eqref{eq:KL} to be compactly supported, it is easy to see that, by Pinsker's inequality,
$$\wasserstein(\mathbf{V}^\top\bu_{\bmu}, \phi(\bmu,\bz))\le c\|f_{\bmu}-\hat{f}_{\bmu,\boldsymbol{\eta}}\|_{L^1}(\mathbb{R}^n)\le c\sqrt{2D_{\text{KL}}(f_{\bmu},\hat{f}_{\bmu,\boldsymbol{\eta}})},$$
where $c>0$ is a suitable constant depending only the sizes of the supports of the two densities.  
This draws a connection between the $\textsf{EKL}$ and the $\textsf{AWD}$ (which can be made more explicit by applying Jensen's inequality), thus linking our theoretical analysis with the minimization of \eqref{eq:nll}.

\begin{remark}
    In practice, the design of the CNF architecture will depend on the problem under consideration. Following standard practice in machine and deep learning, we account for this implicit dependency by tuning the model and training hyperparameters—including the number of coupling layers, the complexity of the MLPs in the architecture, the learning rate, and the number of training epochs—based on NLL performance on the validation set. More precisely, given a finite collection of candidate architectures, $\phi_{\boldsymbol{\eta}}^{(1)},\dots,\phi_{\boldsymbol{\eta}}^{(r)}$, and corresponding training routines, we train all the models separately, obtaining $\phi_{\boldsymbol{\eta}_*}^{(1)},\dots,\phi_{\boldsymbol{\eta}_*}^{(r)}$, and then set $\phi=\phi_{\boldsymbol{\eta}_*}^{j_*}$ where
$$j_*:=\argmin_{j\in\{1,\dots,r\}}\;-\frac{1}{N_{\text{valid}}}\sum_{i=1}^{N_{\text{valid}}}\log \hat{f}_{\boldsymbol{\eta}_*,\bmu_i^{\text{valid}}}^{(j)}(\mathbf{V}^\top\bu_{j}^{\text{valid}})$$
with $\hat{f}_{\boldsymbol{\eta}_*,\bmu}^{(j)}$ being the conditional density associated to the $j$th model $\phi_{\boldsymbol{\eta}_*}^{(j)}.$
\end{remark}

\subsubsection{Deployment at inference time} 
After training, the POD-CNF becomes fully operational. Given any query parameter $\bmu_0\in\pspace$, the model can efficiently generate an arbitrarily large amount of independent samples $\tilde{\bu}_{\bmu_0,1},\dots,\tilde{\bu}_{\bmu_0,J}$ by sampling the latent variables $\bz_1,\dots,\bz_J\sim\mathcal{N}(0,\mathbf{I}_n)$, computing the latent coefficients $\phi(\bmu_0,\bz_j)$ and expanding them over the POD basis, $\tilde{\bu}_{\bmu_0,j}:=\mathbf{V}\phi(\bmu_0,\bz_j).$ The generated ensemble then allows for a complete characterization of the uncertainties and a proper description of the conditional law of $\bu_{\bmu_0}$.

The accuracy of the model can be certified by computing suitable performance metrics over the test set, such as the NLL in the latent space and the POD projection error. Estimating the AWD in \eqref{eq:awd} is also an option,
but this generally only possible in ideal scenarios. In order to approximate \eqref{eq:awd} numerically, in fact, one needs access to a test set with a \emph{nested structure}, namely $\{\bmu_{i}^{\text{test}},\bu_{\bmu_i^\text{test}}(\boldsymbol{\theta}_{i,j})\}$ with ${j=1,\dots,N_{\text{test}}}$ and ${i=1,\dots,\ell}$. While this is generally not always possible, in our experiments we will force this assumption, effectively overwriting our previous definition of "test set". This will be essential to ensure a rigorous and comprehensive assessment of the probabilistic surrogates, which we find particularly relevant given the exploratory and pioneering nature of the present work.

\begin{algorithm}[t]
\SetAlgoLined
\caption{Offline construction of the POD-CNF surrogate}
\label{algo:pod-cnf}
\KwIn{Sampling density $\xi$, total number of snapshots $S$, number of validation samples $N_{\text{valid}}$, number of test samples $N_{\text{test}}$, POD tolerance $\delta$, trainable CNF architectures $\{\phi_{\boldsymbol{\eta}}^{(j)}\}_{j=1}^r$, optimization routines $\{\mathscr{R}_j\}_{j=1}^r$.}\vspace{0.25cm}
\KwOut{POD-CNF $\Phi$}
\vspace{0.25cm}
\textit{Step 1. Data collection}

$[\bmu_1,\dots,\bmu_S]\leftarrow$ sample from $\xi$

\vspace{0.05cm}

\For{$i=1,\ldots,S$}{
    $\mathbf{u}_i\leftarrow$ sample random event $\btheta_i$ and compute $\bu_{\bmu_i}(\btheta_i)$ 
}

\vspace{0.05cm}

$N\leftarrow S-(N_{\text{valid}}+N_{\text{test}})$

\vspace{0.05cm}

\textbf{for} $i=1,\ldots,N_{\textnormal {valid}}$ \textbf{set} $\bu_i^{\text{valid}}\leftarrow \bu_{N+i}$

\vspace{0.05cm}

\textbf{for} $i=1,\ldots,N_{\textnormal {test}}\hspace{0.25cm}$\textbf{set} $\bu_i^{\text{test}}\hspace{0.15cm}\leftarrow \bu_{N+N_{\textnormal{valid}}+i}$

\vspace{0.25cm}

\textit{Step 2. POD basis extraction}

$\mathbf{S}\leftarrow[\bu_1,\ldots,\bu_N]$ assemble snapshots matrix

\vspace{0.05cm}

$\mathbf{U},\mathbf{\Sigma},\mathbf{W}\leftarrow$ compute SVD so that $\mathbf{S}=\mathbf{U}\mathbf{\Sigma}\mathbf{W}^\top$

\vspace{0.05cm}

$\sigma_1,\dots,\sigma_{\min\{N, N_h\}}\leftarrow \text{diag}(\mathbf{\Sigma})$ extract singular values in decreasing order

\vspace{0.05cm}

$n\leftarrow 2\min \{m\;:\;[\sum_{i>2m}\sigma_i^2]/[\sum_{i=1}^{\min\{N_h,N\}}\sigma_i^2]\leq\delta\}$ select even POD dimension

\vspace{0.05cm}

$\mathbf{V}\leftarrow\mathbf{U}[:,1:n]$ extract POD basis via truncation

\vspace{0.25cm}

\textit{Step 3. Train and select CNF module}

\For{$j=1,\ldots,r$}{
$\boldsymbol{\eta}_*^{(j)}\leftarrow$ minimize $ 
\textsf{NLL}^{(j)}(\boldsymbol{\eta})
=-\frac{1}{N}\sum_{i=1}^{N}
\log\hat{f}_{\boldsymbol{\eta},\bmu_i}^{(j)}(\mathbf{V}^\top\bu_i)$ numerically using optimization routine $\mathscr{R}_j$
}

\vspace{0.05cm}

$
j_*\leftarrow\argmin_{j\in\{1,\ldots,r\}}
-\frac{1}{N_{\mathrm{valid}}}
\sum_{i=1}^{N_{\mathrm{valid}}}
\log\hat{f}_{\boldsymbol{\eta}_*}^{(j)}
\left(\mathbf{V}^{\top}\bu_i^{\mathrm{valid}}\right),$
select best-in-validation CNF module

\vspace{0.35cm}
\textit{Assemble the trained POD-CNF}

\vspace{0.02cm}

$\Phi\leftarrow\mathbf{V}\phi^{(j_*)}_{\boldsymbol{\eta}_*}$, optionally print test NLL

\vspace{0.25cm}

\Return $\Phi$
\end{algorithm}



\section{Numerical experiments}
\label{sec:exp}
In this Section, we assess the capabilities of the POD-CNF framework in a controlled simulation setting. Specifically, we investigate whether the proposed approach can approximate the solution of a parametric problem whilst capturing additional stochasticities induced by suitable hidden parameters. This setting thus serves as a proof of concept for the proposed methodology, allowing us to assess its ability to recover the relevant stochastic structures in a well-defined environment. As such, the present analysis represents a preliminary yet necessary step towards more challenging scenarios, including applications to real experimental data, which are beyond the scope of this work.

We consider two benchmark problems: an advection-diffusion process in a partially controlled Stokes flow (Section~\ref{sec:scalar_case}) and a linear elasticity problem involving a biphasic material with a stochastic interface (Section~\ref{sec:bivariate_case}). The two test cases serve to demonstrate the capabilities of POD-CNF in both scalar and vector-valued settings.
In both cases, the goal is to approximate the conditional distribution of the solution field given deterministic parameters, $\bmu$, whilst accounting for the uncertainty induced by unobserved hidden variables, $\btheta$. 
The analysis is based on different performance metrics and benchmarks which are briefly outlined in Section~\ref{sec:exp_setup}.

We anticipate that, while this Section focuses on assessing the generative capabilities of POD-CNF, a preliminary investigation of its application to uncertainty quantification, including forward and inverse problems, is discussed immediately after: see Section \ref{sec:uq}.

\subsection{Numerical setup}
\label{sec:exp_setup} 

\paragraph{Performance metrics}
Let $\bu_{\bmu}$ denote the stochastic ground truth solution, conditional on $\bmu\in\pspace$. Conversely, let $\tilde{\bu}_{\bmu}$ denote the output of a given (trained) probabilistic surrogate. 

As a first metric, we consider the reduced negative-loglikelihood over the test set, computed as
\begin{equation}
    \label{eq:nll_test}
    \textsf{NLL}_{\text{test}}:=-\frac{1}{N}\sum_{i=1}^{N}\log \hat{f}_{\bmu_i^{\text{test}}}(\mathbf{V}^\top\bu_{i}^{\text{test}})
\end{equation}
where $\hat{f}_{\bmu_i}$ is the density of the reduced random vector $\mathbf{V}^\top\tilde{\bu}_{\bmu}$, whereas each $\bu_{i}^{\text{test}}$ is a single independent realization of $\bu_{\bmu_i^\text{test}}$. As discussed in~\ref{sec:practice}, this metric provides a global indication of the goodness of fit; furthermore, it has the advantage of requiring a single FOM evaluation for each $\bmu_i^{\text{test}}.$ In general, this metric is specifically designed for POD-based probabilistic surrogates that admit an absolutely continuous density at the latent level.

As a second set of metrics, we consider more general indicators based on approximated Wasserstein distances. In general, these require repeated conditional evaluations of both the FOM and the surrogate and are thus much more expensive to evaluate. The idea goes as follows. Given any $\bmu_0\in\pspace$, we exploit the FOM and the surrogate to generate $m_{\text{gen}}$ independent realizations of both models, $\bu_{\bmu_0}^{(1)},\dots,\bu_{\bmu_0}^{(m_{\text{gen}})}$ and $\tilde{\bu}_{\bmu_0}^{(1)},\dots,\tilde{\bu}_{\bmu_0}^{(m_{\text{gen}})}$, respectively. The two clouds of points in $\mathbb{R}^{N_h}$ are then compared to each other using a discrete version of the Wasserstein distance $\ewass$, defined as
\begin{multline}
    \label{eq:empirical_wasserstein}
    \ewass(\{\bu_{\bmu_0}^{(1)},\dots,\bu_{\bmu_0}^{(m_{\text{gen}})}\},\;\{\tilde{\bu}_{\bmu_0}^{(1)},\dots,\tilde{\bu}_{\bmu_0}^{(m_{\text{gen}})}\})=\\=\min\left\{\sum_{i,j}\Gamma_{i,j}\|\mathbf{u}_{\bmu_0}^{(i)}-\bu_{\bmu_0}^{(j)}\|\;\text{where}\;\Gamma\in\mathbb{R}^{m\times q}\;:\;\sum_{j}\Gamma_{i_0,j}=\sum_{i}\Gamma_{i,j_0}=1/m_{\text{gen}}\;\forall i_0,j_0\right\}.
\end{multline}
The above is used as proxy for the true Wasserstein distance between the FOM and generated distributions, that is
\begin{equation}
     \label{eq:wasserstein-approx}\wasserstein(\bu_{\bmu_0},\tilde{\bu}_{\bmu_0})\approx \ewass(\{\bu_{\bmu_0}^{(1)},\dots,\bu_{\bmu_0}^{(m_{\text{gen}})}\},\;\{\tilde{\bu}_{\bmu_0}^{(1)},\dots,\tilde{\bu}_{\bmu_0}^{(m_{\text{gen}})}\}).
\end{equation}
In our implementation, we compute \eqref{eq:empirical_wasserstein} using the \textsf{pot} Python library, with a computational overload of $O(m_{\text{gen}}^3),$ excluding generation costs.

This local analysis is exploited in two ways. First, to us to assess the performances of the model over a specific parametric scenario $\bmu_0\in\pspace$. Second, to investigate the accuracy of the model throughout $\pspace$; this is done by computing \eqref{eq:empirical_wasserstein} for multiple, independent, parametric configurations, randomly drawn from $\pspace$. We generally set $m_{\text{gen}}=1000$ in the first case and $m_{\text{gen}}=200$ in the second case to mitigate the computational overload. 

\paragraph{Benchmarks}
We consider two benchmarks against which we compare our proposed strategy, in two different test cases. The first benchmark 
represents an idealized comparison with purely \emph{deterministic} surrogates. Rather than implementing and evaluating individual deterministic methods, this analysis is carried out by considering an ideal reference predictor, which serves as a baseline for all deterministic models. To this end, we notice that for any $\bmu_0\in\pspace$, all deterministic surrogates must return a single, nonaleatoric prediction in $\mathbb{R}^{N_h}$. In mean-square terms, it is well-known that the optimal such prediction is given by the conditional expectation $\mathbb{E}[\bu_{\bmu}\mid\bmu=\bmu_0]=\mathbb{E}[\bu_{\bmu_0}]$, see, e.g., \citep[Chapter 9]{williams1991probability}. We therefore take $\bmu\mapsto\mathbb{E}[\bu_{\bmu}]$ as our ideal deterministic baseline: no deterministic surrogate can improve upon its mean-square accuracy, regardless of its architecture or training procedure. In practice, for any given $\bmu_0\in\pspace$, we approximate the ideal prediction by the empirical mean of a FOM-generated population, $\overline{\bu}_{\bmu_0}:=m_{\text{gen}}^{-1}\sum_{i=1}^{m_{\text{gen}}}\bu_{\bmu_0}^{(i)}$.
From a probabilistic perspective, because the output of a deterministic surrogate is intrinsically degenerate, its accuracy cannot be assessed reliable using the negative log-likehood. For this reason, when analyzing the deterministic baseline, we focus on Wasserstein-based metrics only. In this case \eqref{eq:wasserstein-approx} becomes
\begin{equation}
\ewass(\{\bu_{\bmu_0}^{(i)}\}_{i=1}^{m_{\text{gen}}},\;\{\overline{\bu}_{\bmu_0}\})
=\frac{1}{m_{\text{gen}}}\sum_{i=1}^{m_{\text{gen}}}\|\bu_{\bmu_0}^{(i)}-\overline{\bu}_{\bmu_0}\|
\approx\mathbb{E}\|\bu_{\bmu_0}-\mathbb{E}[\bu_{\bmu_0}]\|
=\wasserstein(\bu_{\bmu_0},\mathbb{E}[\bu_{\bmu_0}]).
\end{equation}

As a second benchmark, instead, we consider POD-GP \cite{guo2018reduced}.
The latter is an established approach for probabilistic reduced-order modeling that is based on the integration of POD and Gaussian Process (GP) regression. From a generative standpoint, POD-GP can be seen as a model of the form
\begin{equation}
    \Phi_{\text{POD-GP}}(\bmu,\bz)=\sum_{j=1}^{n}\left[m_j(\bmu)+z_js_j(\bmu)\right]\mathbf{v}_j,
\end{equation}
where $z_j$ and $\mathbf{v}_j$ are the $j$th component and the $j$th column of $\bz\sim\mathcal{N}(0,\mathbf{I}_n)$ and $\mathbf{V}$, respectively, whereas $m_j$ and $\sigma_j$ are the predicted mean and standard deviation returned by the $j$th GP. In the POD-GP approach, in fact, a sequence of $n$ independent GPs is fitted to approximate the parameter-to-reduced state map, where each GP focuses on a different POD coefficient.
In our experiments, we fit POD-GP by leveraging the same data used to design and train POD-CNF. Specifically, when implementing POD-GP, we exploit the validation set to estimate the hyperparameters of each GP, and use the training set to fit the models and obtain the final surrogate. 
All GPs are implemented using a squared-exponential (RBF) covariance kernel for the prior, scaled by a trainable constant, together with an additive white-noise kernel. The signal variance, length scale, and noise variance are optimized in the ranges $[10^{-6},10^{2}]$, $[10^{-2},10^{1}]$, and $[10^{-8},10^0]$, respectively. To ensure that such ranges are meaningful, all inputs and outputs are suitably normalized before training.

An important feature of POD-GP is that its predictive distributions are always Gaussian. Precisely, if we denote by $\mathbf{m}_{\text{GP}}(\bmu):=[m_1(\bmu),\dots m_n(\bmu)]$ and $\mathbf{\Sigma}_{\text{GP}}(\bmu):=\text{diag}\left(s_1(\bmu),\dots,s_n(\bmu)\right)$, then $\Phi_{\text{POD-GP}}(\bmu,\cdot)\sim \mathcal{N}(\mathbf{V}\mathbf{m}(\bmu), \mathbf{V}\mathbf{\Sigma}(\bmu)\mathbf{V}^\top)$. Similarly, at the latent level, $\hat{f}_{\bmu}$ coincides with the density of a multivariate normal $\mathcal{N}(\mathbf{m}(\bmu), \mathbf{\Sigma}(\bmu)).$

\subsection{Advection-Diffusion in a Stokes Flow}
\label{sec:scalar_case}

\begin{figure}
 \centering
   \includegraphics[width=0.37\textwidth]{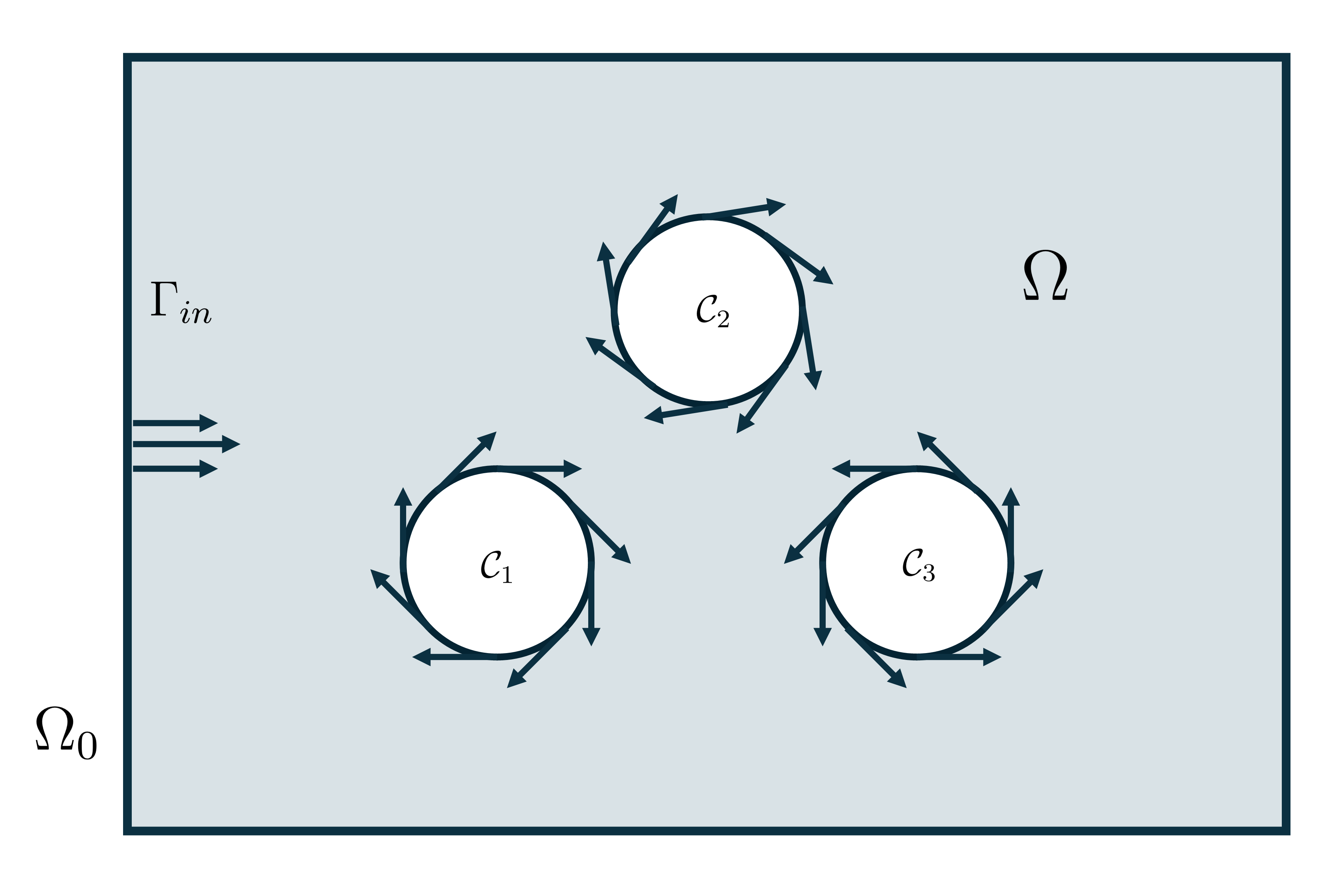}
   \includegraphics[width=0.37\textwidth]{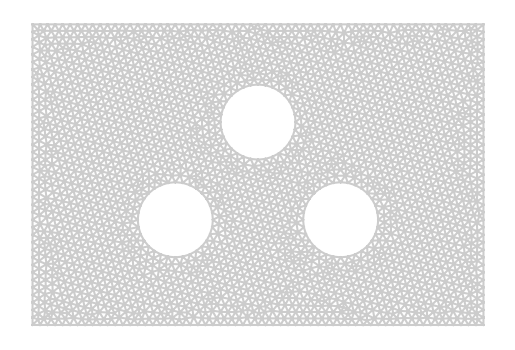}
   \caption{Problem domain for the advection-diffusion problem (left) and corresponding discretization (right).}\label{fig:mesh_stokes}
\end{figure}

As a first test case, we consider the passive transport of a chemical in a rectangular chamber equipped with three rotating cylinders, cf. Fig.~\ref{fig:mesh_stokes}. Here, $\Omega=\Omega_0\setminus(\mathcal{C}_1\cup \mathcal{C}_2\cup \mathcal{C}_3)$ with $\Omega_0=(0,3)\times(0,2)$. This setup simulates a physical mixer where the cylinders rotate at specific velocities, driving the advection field. The objective is to construct a surrogate that can predict the  concentration field of the chemical, $u=u(\bx)$, given the configuration of the mixer, while simultaneously accounting for potential uncertainties.

To model the ground truth mechanism, and produce the desired variability, we use the following sets of parameters.

\begin{figure}[tb]
 \centering
   \includegraphics[width=0.8\textwidth]{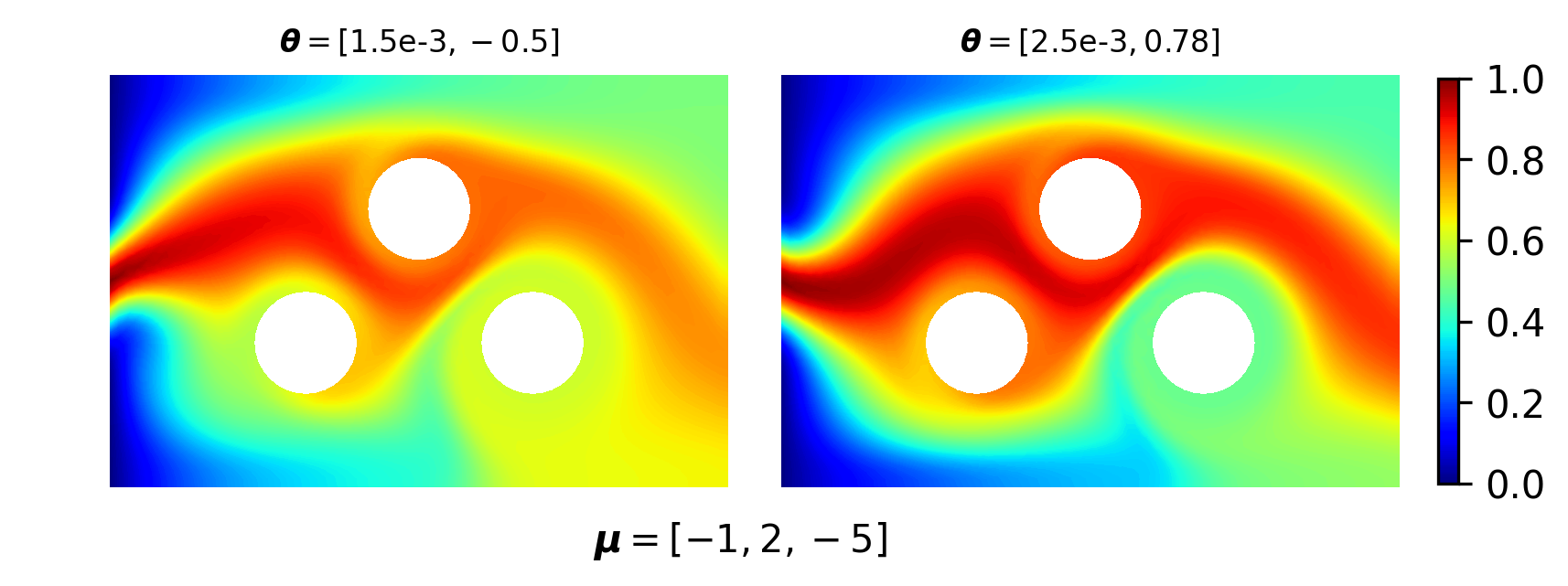}
   \caption{Advection diffusion problem. Two random solutions associated with the same operating conditions $\bmu=[-1,2,-5]$. The variability is induced by the hidden parameters, $\btheta$, whose values are assumed to be unknown (to the user, and to the surrogate). \label{fig:stokes_setup}}
\end{figure}

\begin{itemize}
    \item \textbf{Conditioning Parameters $\bmu=[c_1,c_2,c_3]$.} Here, $c_i$ is the rotational velocity of the $i$th cylinder. Such parameters encode \emph{what is known} about the configuration of the system, serving as controllable inputs for the mixer. We allow each $c_i$ to vary in a pre-defined range, $c_i\in[-5,5].$
    
    \item \textbf{Hidden Parameters $\btheta=[\kappa,\alpha]$.} We use $\kappa$ to model the diffusion coefficient of the chemical, and $\alpha$ for the inflow angle of the fluid. We exploit such parameters to introduce an additional, formally \emph{unknown}, source of variability, which we obtain by letting $\kappa\sim\mathcal{U}[0.0025, 0.005]$ and $\alpha\sim\mathcal{U}[-\pi/4, \pi/4]$. The first hidden parameter accounts for the fact that, although the identity of the diffusing chemical is known, its effective transport properties may vary across experiments as a result of changes in, e.g., temperature, solvent composition, or other experimental conditions that can alter key properties such as the molecular diffusivity. These alternations are assumed to be out-of-control, making the effective value of $\kappa$ formally unknown form the outside. Similarly, the stochasticity in the inflow condition can be though as a way to model imperfections in the injection mechanism or other perturbations taking place outside of the mixing chamber.
\end{itemize}
Due to the stochastic nature of $\kappa$ and $\alpha$, a fixed configuration of cylinder velocities $\bmu$ does not yield a unique solution but rather a distribution of physically consistent concentration fields $\mathsf{p}(u | \bmu)$: see, e.g., Fig.~\ref{fig:stokes_setup}.

\paragraph{Mathematical Formulation and Discretization}

As anticipated, we model the concentration field $u(\bx; \bmu,\btheta)$ as the solution to the following advection-diffusion equation,
\begin{equation}
\label{eq:AD}
\begin{cases}
-\kappa\Delta u(\bx; \bmu,\btheta) + \bb(\bx; \bmu, \btheta)\cdot\nabla u(\bx; \bmu,\btheta) = 0 & \text{in } \Omega, \\
\nabla u(\bx; \bmu,\btheta) \cdot \bn = 0 & \text{on } \partial\Omega\setminus\Gamma_{in},\\
u(\bx; \bmu,\btheta) = g(\bx) & \text{on } \Gamma_{in},
\end{cases}
\end{equation}
where $\Gamma_{in} = \{0\}\times[0,2]$, $\bn$ is the external unit normal, $g(\bx) = g(x_1,x_2) = \exp{(-25(x_2-1)^2)}$ is a fixed boundary condition, while $\bb$ is the underlying velocity field. The latter, which depends on $\alpha, c_1, c_2, c_3$, is determined by the stationary Stokes equations for an incompressible viscous fluid, neglecting inertial terms:
\begin{equation}
\label{eq:stokes_strong}
\begin{cases}
-\Delta \bb(\bx; \bmu, \btheta) + \nabla p(\bx; \bmu, \btheta) = \bzero & \text{in } \Omega, \\
\nabla \cdot \bb(\bx; \bmu, \btheta) = 0 & \text{in } \Omega,
\end{cases}
\end{equation}
where $p$ represents the fluid pressure. The above is complemented by the following boundary conditions,
\begin{equation}
\label{eq:stokes_bc}
\begin{cases}
\bb(\bx; \bmu, \btheta) = [\cos\alpha, \sin\alpha]^\top g(\bx) & \text{on } \Gamma_{in}, \\
\bb(\bx; \bmu, \btheta) = c_i [-(x_2-x_2^i), x_1-x_1^i]^\top & \text{on } \partial\mathcal{C}_i, \quad i \in \{1,2,3\}, \\
\bb(\bx; \bmu, \btheta) = \bzero & \text{on } \partial\Omega\setminus(\Gamma_{in} \cup \partial\mathcal{C}_1 \cup \partial\mathcal{C}_2 \cup\partial \mathcal{C}_3),
\end{cases}
\end{equation}
where $(x_1^i,x_2^i)$ is the center of the  $i$th cylinder $\mathcal{C}_i$.
The physics of the problem is deliberately kept simple, as our primary interest lies in addressing the parametric variability of the solution with respect to $\bmu$ and $\btheta$.


As reference solver, we consider a FOM that approximates the solution to \eqref{eq:AD}–\eqref{eq:stokes_bc} using the Finite Element Method (FEM). To this end, we discretize the spatial domain $\Omega$ using an unstructured triangular mesh comprising $N_h=2763$ vertices; see Fig.~\ref{fig:mesh_stokes}.
We use P1 continuous elements for the advection-diffusion equation, whereas, in order to satisfy the inf-sup condition, we use the stable Mini-Element pair to solve for the Stokes' flow equation (continuous P1 elements for the pressure, and P1-bubble for the velocity).


\paragraph{Dataset and POD-CNF configuration}
For the deterministic parameter, we consider a uniform sampling density over the parameter space $\pspace=[-5,5]^3$, so that $\xi\sim\mathcal{U}(\pspace)$. We run the FOM repeatedly in order to produce a total of $S=6400$ independent simulations, which we split into training $(N=4800)$, validation ($N_\text{val}=1280$), and test ($N_\text{test}=320$) instances. For the POD, we set an energy threshold of $\tau=10^{-4}$. For the CNF module, instead, we select the best-in-validation design among architectures with  $K\in\{12, 16, 24, 32\}$ alternating affine coupling layers, scaling and translation MLPs of depth $\ell\in\{1,2\}$ and width $w\in\{64, 128, 256, 512\}$ for each coupling layer, and training routines based on the Adam optimizer \cite{kingma2015adam} with learning rate $\textsf{lr}\in\{10^{-3}, 5 \times 10^{-4}\}$, a batch size of $64$, early stopping with a patience of $10$ epochs, weight decay set to $10^{-5}$, and a total of $50$ epochs per configuration.





\begin{figure}
    \centering
    \includegraphics[width=0.55\textwidth]{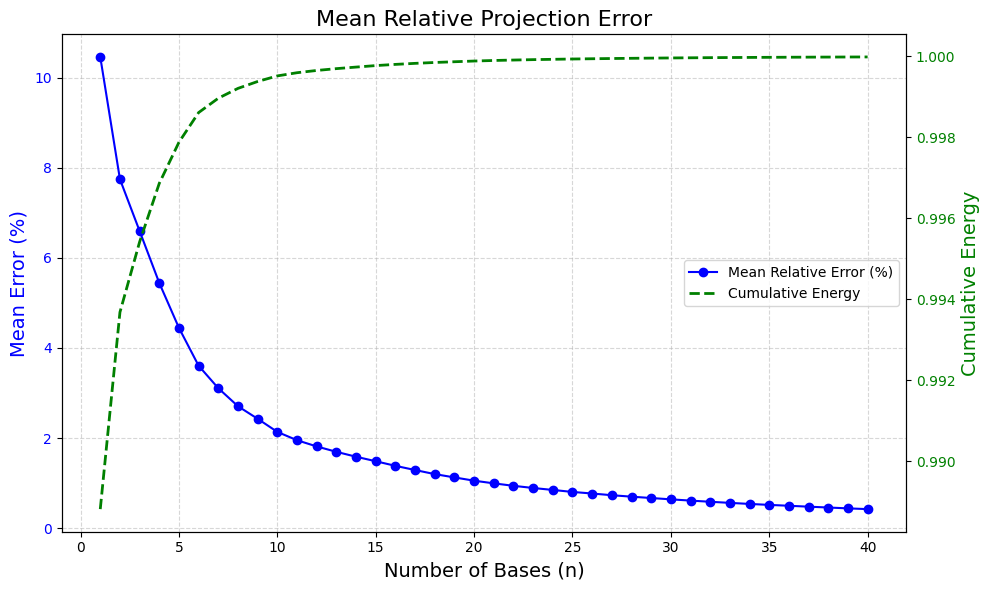}
    \caption{Advection-diffusion problem, \S\ref{sec:scalar_case}. Mean relative reconstruction error (solid blue line) and cumulative energy (dashed black line) versus the number of POD bases $n$. Truncating at $n=20$ guarantees a projection error in euclidean norm equal to $1.05\%$.}
    \label{fig:svd_stokes}
\end{figure}

\paragraph{Results I: final design of the POD-CNF surrogate.}
Following the methodology described in Section~\ref{sec:method}, we computed the SVD of the training snapshot matrix. 
The singular value and truncation analysis (Fig.~\ref{fig:svd_stokes}) revealed a distinct---but not sharp---energy drop-off. Still, $n=20$ modes eventually suffice to meet the chosen energy threshold. On the validation set, the resulting POD basis $\mathbf{V}$ striked a Mean Relative Projection Error (MRPE) of $1.05\%$,
\begin{equation}
    \label{eq:mrpe}\mathsf{MRPE}_{\text{valid}}:=\frac{1}{N_{\text{valid}}}\sum_{i=1}^{N_{\text{valid}}}\left[\frac{\|\bu_{i}^{\text{valid}}-\mathbf{V}\mathbf{V}^\top\mathbf{u}_{i}^{\text{valid}}\|}{\|\mathbf{u}_{i}^{\text{valid}}\|}\right],
\end{equation}
confirming the accuracy of the reduced representation.

The CNF was thus configured to operate in a $20$-dimensional space under the conditioning of the physical parameters $\bmu\in\mathbb{R}^3.$ Following the results of hyperparameter tuning procedure, a CNF architecture with $K=24$ coupling layers, and internal MLP units of depth $\ell=2$ and width $w=128$, was selected. The latter was then trained for 100 epochs using the Adam optimizer with $\textsf{lr}=10^{-3}$ and a weight decay of $10^{-5}$, to obtain the final model. The offline training procedure proved to be highly efficient, requiring approximately $10$ minutes on a single NVIDIA T4 GPU, cf. Table~\ref{tab:mu0-stokes}.


\begin{figure}
  \centering
  \parbox{0.04\textwidth}{\textbf{(a)}}%
  \parbox{0.94\textwidth}{\includegraphics[width=\linewidth]{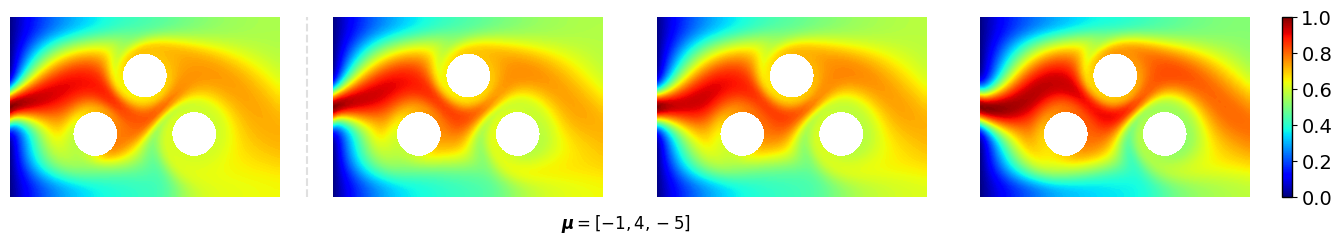}}\\[0.8em]

  \parbox{0.04\textwidth}{\textbf{(b)}}%
  \parbox{0.94\textwidth}{\includegraphics[width=\linewidth]{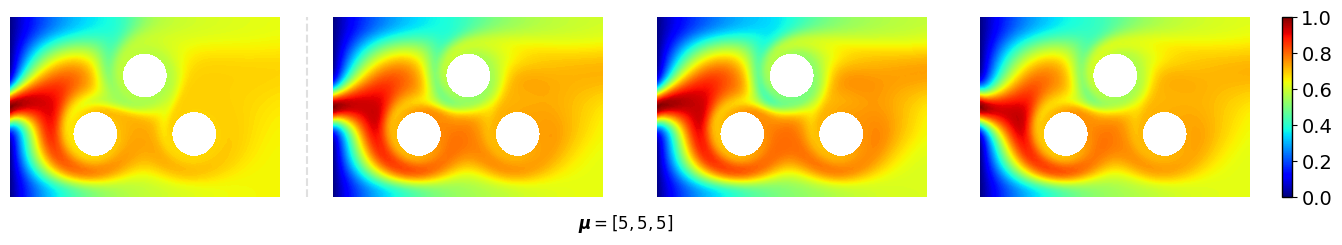}}\\[0.8em]

  \parbox{0.04\textwidth}{\textbf{(c)}}%
  \parbox{0.94\textwidth}{\includegraphics[width=\linewidth]{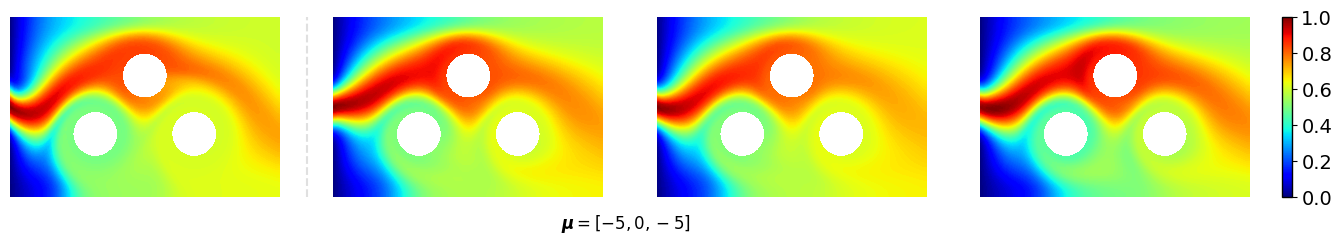}}

  \caption{Advection-diffusion problem, \S\ref{sec:scalar_case}. Solutions fields across different parametric configurations: 
    (a)~$\bm{\mu}_0 = [-1, 4, -5]$, 
    (b)~$\bm{\mu}_1 = [5, 5, 5]$, and 
    (c)~$\bm{\mu}_2 = [-5, 0, -5]$. For each row, the leftmost panel shows a ground-truth snapshot, computed with the numerical solver, while the three panels to its right show independent realizations sampled from the POD-CNF surrogate model. Configuration~(a) corresponds to the parameter $\bmu_0\in\pspace$ chosen for the local analysis.}
  \label{fig:recon_stokes}
\end{figure}

\begin{figure}[h]
    \centering
    \includegraphics[width=0.85\textwidth]{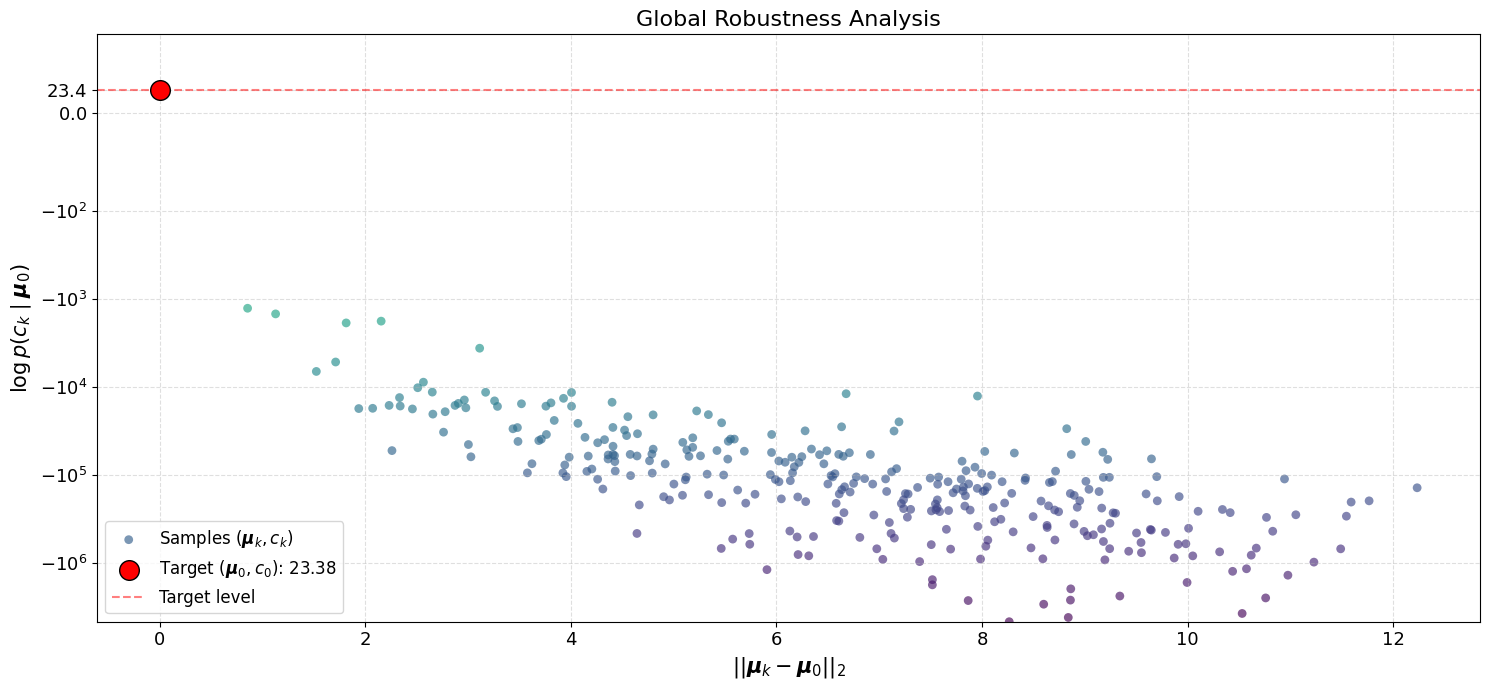}
    \caption{Reduced log-likelihood sensitivity of the POD-CNF surrogate evaluated at $\bmu_0$ against solutions belonging to different parameter instances $\bmu_k$.}
    \label{fig:loglike}
\end{figure}

\paragraph{Results II: model performance on a selected $\bmu_0\in\pspace$.} Evaluating the performance of a probabilistic surrogate is, by itself, a fairly complex task. For this reason, we start by assessing the accuracy of the trained POD-CNF on a representative test configuration, obtained by fixing the values of the
physical parameters to $\bmu_0 := [-1, 4, -5]$. In particular, we aim at understanding if the generative model is actually able to capture the conditional law of $\bu\mid\bmu=\bmu_0$.

Fig.~\ref{fig:recon_stokes}a compares a random solution obtained using the FOM, with $\bmu=\bmu_0$ and $\btheta\sim\mathcal{U}[0.0025, 0.005]\otimes\mathcal{U}[-\pi/4, \pi/4]$ drawn at random, against three \emph{i.i.d.} realizations sampled from the trained POD-CNF model. As we can see, all generated solutions exhibit advection patterns that are consistent with the selected configuration of the system (recall that $\bmu_0$ defines the speed of the rotating cylinders). The variability in the outputs of the POD-CNF shows how the latter is effectively attempting at capturing an intrinsic stochasticity that is not imputable to $\bmu$ alone. Most importantly, such variability is not uncontrolled, but rather presents a precise structure that makes the predicted outputs effectively compatible with the FOM simulation.

Further evidence of this fact can be found in Fig.~\ref{fig:loglike}. There, we compare the POD-CNF log-likelihood (in reduced space) of different FOM simulations conditionally on $\bmu_0$. Each of these simulations has been extracted from the test set and corresponds to a different $\bmu\in\pspace$.
The reduced log-likelihood shows a clear trend, with FOM simulations becoming less and less likely as their physical configuration, $\bmu$, departs from the one assumed by the POD-CNF, $\bmu_0$. This confirms that the proposed surrogate is able to discriminate in between different scenarios, effectively assigning a higher likelihood to those physical solutions that are consistent with the prescribed input condition.

\begin{figure}
    \centering
    \includegraphics[width=0.5\textwidth]{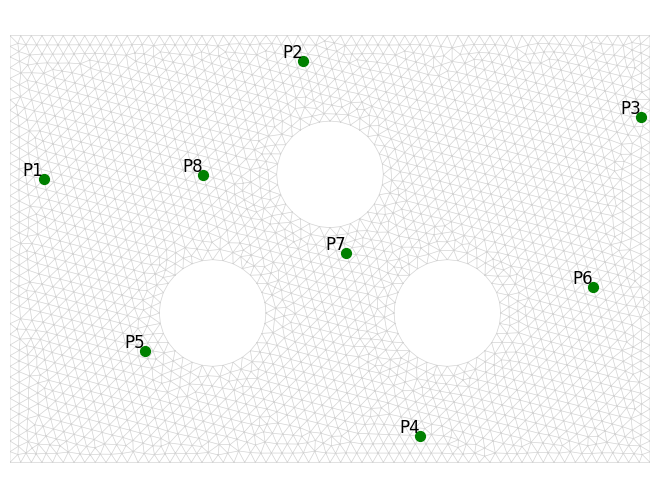}
    \includegraphics[width=\textwidth]{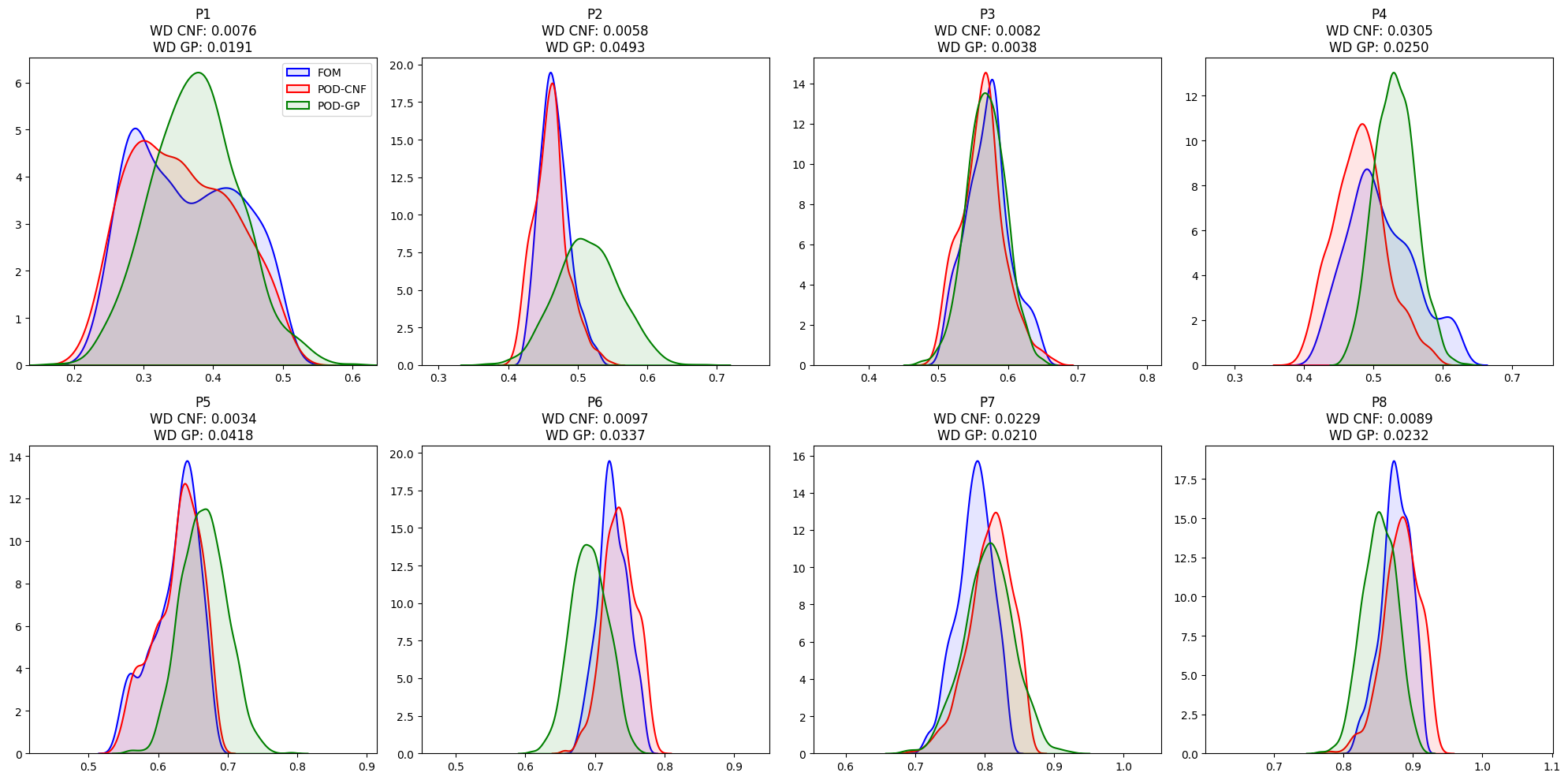}
    \caption{Advection-diffusion problem, \S\ref{sec:scalar_case}. Comparison between FOM and POD-CNF marginal densities (bottom panel) at eight selected monitoring points (top panel).}
    \label{fig:distribution_comparison_6179}
\end{figure}

To corroborate these qualitative findings, a complementary statistical analysis, comparing true and predicted distributions, has been carried out. To this end, we sampled $m_{\text{gen}}=1000$ new simulations using both the FOM and the POD-CNF, all conditional on $\bmu=\bmu_0.$ The two populations have then been compared using the (discrete) Wasserstein distance $\ewass$, computed as a proxy for $\wasserstein(\bu_{\bmu_0}, \mathbf{V}\phi(\bmu_0))$: cf. Section~\ref{sec:exp_setup}. Results are in Table~\ref{tab:mu0-stokes}. Notably, POD-CNF scores particularly well, with an estimated Wasserstein distance of 0.6735. In particular, the method results 56\% more accurate that \emph{any} deterministic surrogate potentially conceivable, demonstrating that its capabilities extend far beyond the standard regression to the mean. 
In contrast, when repeating the same analysis for POD-GP, we note a significant drop in accuracy.
This latter can fact can be explained in at least two ways. The first one, is that the current analysis is entirely focused on a limit scenario, as the parametric configuration under study falls \emph{at the boundary} of the parameter space, $\mu_0\in\partial\pspace.$ As such, even with the prescribed training data, the GP remains highly uncertain about its prediction: this results in an inflation of the conditional variance $\text{Var}(\bu|\bmu_0)$, and a consequent degradation in the approximation of the conditional law (see also Fig.~\ref{fig:distribution_comparison_6179}). The second explanation concerns the fact that the predictive distribution of a GP is necessarily Gaussian. In our setting, where stochasticity arises from hidden variables rather than from external white noise, this simplification ends up posing major limitations.

Our conclusion are further supported by the results shown in Fig.~\ref{fig:distribution_comparison_6179}. By focusing on eight monitoring points spread across the spatial domain, the latter compares the marginal distributions of the true FOM population against the ones obtained with POD-CNF and POD-GP, respectively. 
Because 1-dimensional distributions are much easier to inspect compared to high-dimensional ones, this allows us to make both quantitative and qualitative considerations.
For instance, we can see how POD-CNF correctly identifies regions with higher volatility, matching both the range and distribution shape of the values attained by the solution field (point P1). Highly peaked distributions (points P5–P6) are also accurately represented, even when they exhibit asymmetry. This capability is particularly valuable if we consider that conventional probabilistic surrogates, such as POD-GP, generally fail at reproducing these features.
Clearly, the model also shows some limitations in certain points of the domain (P7-P8), although the performances remain satisfactory.

\begin{figure}
    \centering
    \includegraphics[width=0.6\textwidth]{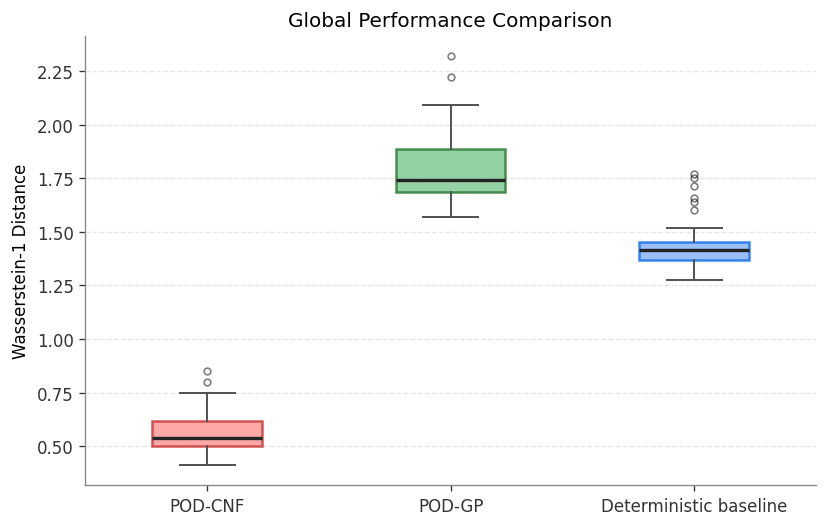}
    \caption{Advection-diffusion problem, \S\ref{sec:scalar_case}. Global analysis over 50 unseen test configurations. The boxplots show the distribution of the Wasserstein-1 error over the considered test set, comparing FOM and surrogate outputs. From left to right: POD-CNF (red, mean $ = 0.5673$), POD-GP (green, mean $ = 1.7970$) and deterministic baseline (blue, mean $ = 1.4369$).} 
    \label{fig:global_w1_stokes}
\end{figure}

\paragraph{Results III: global performance.}
To rigorously evaluate the generalization capabilities of POD-CNF, the previous analysis has been extended
to 50 unseen parameter values, randomly drawn from the test set. In this case, in order to mitigate the computational overload, the discrepancy in $\wasserstein$ metric has been estimated using only $m_{\text{gen}}=200$ samples per parametric configuration.~\ref{fig:global_w1_stokes} shows the distribution of such errors across the 50 configurations under study for POD-CNF, POD-GP and the deterministic baselines. Averages and standard deviations are reported in Table~\ref{tab:mu0-stokes}. Notably, POD-CNF consistently outperforms its competitors, 
demonstrating remarkable generalization capabilities across all parametric regimes. 
The same conclusion can be drawn by looking at the negative log-likelihood, Table~\ref{tab:mu0-stokes} (note: being computationally cheaper to evaluate, the latter as been computed on the whole test set). Finally, Fig.~\ref{fig:recon_stokes}a-c shows how POD-CNF is able to produce meaningful solutions conditionally on the desired parametric configuration. Of note, we can see that the POD-CNF generated samples show a clear variability in the inflow condition at the left boundary. This is particularly interesting if we consider that our model was never given information about the existence of an inflow angle $\alpha$. In this sense, we can see how POD-CNF is effectively capturing the stochasticity produced by the latent, hidden, parameters. This is fundamentally different from what POD-GP normally does. The latter, in fact, is built based on a Bayesian paradigm, which makes the model sensitive to \emph{epistemic} uncertainty. While this feature can be particularly valuable in certain applications, here it can obscure the \emph{aleatoric} component of uncertainty, which is our primary focus.

\renewcommand{\arraystretch}{1.5}
\begin{table}\centering
    \centering
    \begin{tabular}{lllrrr}
        \hline
        
        & $\ewass$ at $\bmu_0$
        & Average $\ewass$ on $\mathcal{P}$
        & $\textsf{NLL}_{\text{test}}$ & Training time & Inference time\\
        \hline
        \textbf{Deterministic baseline}  & 1.5284 & 1.4369 $\pm$ 0.1103& --  & -- & --\\
        \textbf{POD-GP}                 &  2.3819 &    1.7970 $\pm$ 0.1605
        &  4.5315 $\pm$ 0.9112 & 54s  & 0.29s\\
        \textbf{POD-CNF}                 & 0.6735 & 0.5691 $\pm$ 0.0997& -20.7059 $\pm$ 7.6416 & 10m 33s & 0.06s\\
        \hline
    \end{tabular}
    \caption{Advection-diffusion problem, \S\ref{sec:scalar_case}. Performance metrics and training times for POD-CNF, POD-GP and the deterministic baseline (cf. Section~\ref{sec:exp_setup}). First column: surrogate under study. Second column: empirical Wasserstein distance between FOM-generated and surrogate-generated samples ($m_{\text{gen}}=1000$). Third column: average Wasserstein distance, computed across $50$ random parametric configurations ($m_{\text{gen}}=200)$, and corresponding standard deviation., 
    Fourth column: negative loglikelihood, averaged across the whole test set (320 joint samples in $\bmu,\btheta$), reported with a 95\% confidence interval. Fifth column: fitting time, excluding hyperparameter tuning. Last column: \emph{sequential} time required to generate a single scenario for a given parametric configuration.}
    \label{tab:mu0-stokes}
\end{table}

\paragraph{Results IV: sample size sensitivity analysis.}
We conclude by presenting the results of a complementary analysis, devoted at assessing the sensitivity of the CNF module to the number of training data. 
To this end, we re-trained the CNF from scratch, keeping the same POD basis ($n = 20$) and the identified hyperparameters, but progressively reducing the amount of training samples, from 100\% (4800 snapshots) down to 5\% (240 snapshots).~\ref{stokes_sens} shows the impact of the latter on the performances of POD-CNF. 
%

\begin{table}\centering
    \begin{tabular}{lllll}
    \toprule
    \textbf{\% Data} & $N_{\text{train}}$ & $N_{\text{val}}$ & \textbf{Average} $\ewass$ \textbf{on} $\pspace$ & $\ewass$ \textbf{at} $\bmu_0$ \\
    \midrule
    100\% & 4800 & 1280 & 0.5691  $\pm$ 0.0997 & 0.6735 \\
    75\%  & 3600 & 960  & 0.5102  $\pm$ 0.0872 & 0.6675 \\
    50\%  & 2400 & 640  & 0.5487  $\pm$ 0.0943 & 0.8187 \\
    25\%  & 1200 & 320  & 0.6607  $\pm$ 0.1287 & 0.9352 \\
    10\%  & 480  & 128  & 1.5019  $\pm$ 0.4024 & 1.8266 \\
    5\%   & 240  & 64   & 2.0273  $\pm$ 0.4023 & 2.3104 \\
    \bottomrule
    \end{tabular}\caption{Advection-diffusion problem, \S\ref{sec:scalar_case}. Sensitivity of POD-CNF performance to training dataset size. Evaluation metrics read as in Table~\ref{tab:mu0-stokes}.} 
\label{stokes_sens}
\end{table}

Interestingly, the model remains quite stable even when halving the amount of training data. This suggests a convergence in the performance: in order to make the surrogate more accurate, instead of adding more data, one would need to either increase the expressivity of the CNF or the number of POD basis. 
At 25\% sample size, the performance starts to deteriorate slightly. However, the most notable drop happens at 10\% (and below), where POD-CNF loses most of its value: there, in fact, its performances essentially become on par with the deterministic baseline and POD-GP, cf. Table~\ref{tab:mu0-stokes}.

\newpage

\subsection{Linear Elasticity}
\label{sec:bivariate_case}

\begin{figure}
  \centering
  \includegraphics[width=0.3\textwidth]{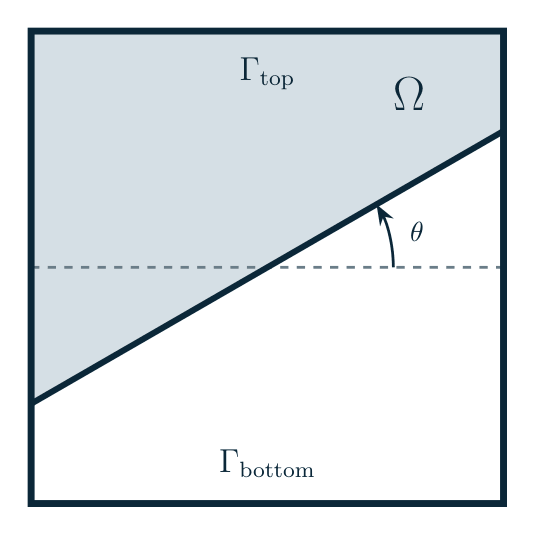}\label{fig:pb_setup_mat}
  \includegraphics[width=0.3\textwidth]{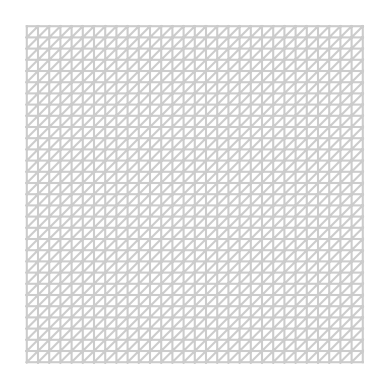}
  
  \caption{Linear elasticity problem, \S\ref{sec:bivariate_case}. Physical domain (left) and corresponding discretization (right).}\label{fig:pb_setup}
\end{figure}

The second numerical experiment extends the proposed framework to a vector-valued setting, addressing the mechanical response of a square block $\Omega = (0,1)^2$ composed of a biphasic material. The domain features two distinct phases—a "soft" and a "hard" material—separated by a linear interface with an unobserved orientation angle $\theta$ (Fig.~\ref{fig:pb_setup}). Deformation is induced by 
a localized vertical force, acting at the top boundary. The system is governed by the following parameters.

\begin{figure}
    \centering
    \includegraphics[width=0.7\linewidth]{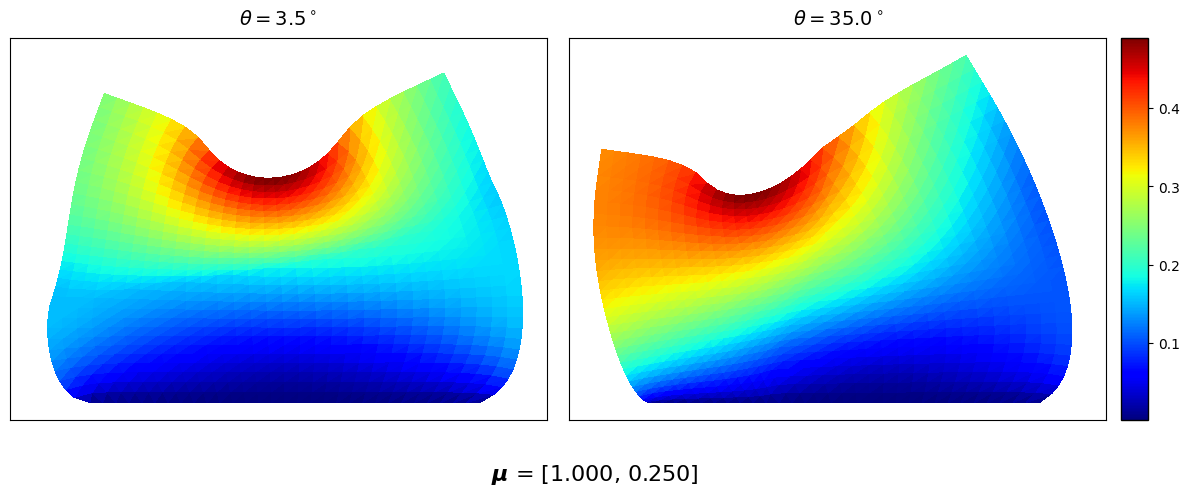}
    \caption{Linear elasticity problem, \S\ref{sec:bivariate_case}. Two random solutions associated with the same operating conditions $\bmu=[1.0, 0.25]$. The variability is induced by the hidden parameters, $\btheta$, whose values are assumed to be unknown (to the user, and to the surrogate).}
    \label{fig:elastic_varia}
\end{figure}

\begin{itemize}
    \item \textbf{Conditioning Parameters $\bmu=[m, \delta]$.} We use these controllable parameters to model the magnitude $(m)$ and the locality $(\delta)$ of the vertical force acting on the box. 
    We allow $m$ and $\delta$ to vary uniformly, and independently, in the ranges $[1, 2]$ and $[0.15, 0.25]$, respectively, so that $\pspace:=[1,2]\times[0.15,0.25]$.
    
    \item \textbf{Latent Parameter $\theta$.} We consider a single unobserved parameter, $\theta \in [0, 2\pi)$, which we use to prescribe the orientation of the linear interface separating the two material phases in the physical domain, cf Fig.~\ref{fig:pb_setup}. 
   Intuitively, this setup can be thought as a simplified model for an indentation test performed on a heterogeneous, biphasic material block whose internal interface orientation cannot be determined non-destructively prior to loading.
\end{itemize}
We let $\theta\sim\mathcal{U}(0,2\pi)$.
Consequently, for a fixed loading configuration $\bmu\in\pspace$, the displacement field $\boldsymbol{u} = (u_x, u_y)$ is not unique but rather follows a conditional distribution $\mathsf{p}(\boldsymbol{u} | \bmu)$ dependent on the random realization of the material interface: see Fig.~\ref{fig:elastic_varia}.

\paragraph{Mathematical Formulation and Discretization.}

We model the system using the equations of linear static elasticity:
\begin{equation}
    \label{eq:elasticity}
    \begin{cases}
        -\nabla\cdot\left(\lambda_\theta(\bx)(\nabla\cdot\boldsymbol{u} )\identity + \nu\left(\nabla\boldsymbol{u} +\nabla\bu^{\top}\right)\right) = - \bg & \textnormal{in}\;\Omega,\\
        \boldsymbol{u} \equiv\bzero & \textnormal{on}\;\Gamma_{bottom},\\
        \lambda_\theta(\bx)(\nabla\cdot \boldsymbol{u} )\bn+\nu\left(\nabla\bu+\nabla\boldsymbol{u} ^{\top}\right)\bn\equiv\mathbf{f}_{m, \delta}(\bx) & \textnormal{on}\;\Gamma_{top},\\
        \lambda_\theta(\bx)(\nabla\cdot \boldsymbol{u} )\bn+\nu\left(\nabla\bu+\nabla\boldsymbol{u} ^{\top}\right)\bn\equiv\bzero & \textnormal{on}\;\partial\Omega\setminus(\Gamma_{top}\cup\Gamma_{bottom}).
    \end{cases}
\end{equation}
The first equation represents the static force balance within the domain, where the body force $-\bg$ accounts for the effect of gravity, and is assumed to be fixed. This equilibrium is subject to mixed boundary conditions. On the lower boundary $\Gamma_{bottom}$, a homogeneous Dirichlet condition $\boldsymbol{u}  \equiv \bzero$ is imposed, meaning that the block is rigidly clamped and experiences no deformation at the bottom. On the upper boundary $\Gamma_{top}$, the system is subjected to a localized traction force $\mathbf{f}$ exerted by the external loading. This Neumann condition is modeled as a downward vertical force with parametrized magnitude and position, here defined as $$\bbf_{m,\delta}(\bx) = \bbf_{m,\delta}(x_1,x_2) = [0, -m]^\top \indicator_{[\frac{1}{2}-\delta, \frac{1}{2}+\delta]}(x_1).$$ The remaining lateral boundaries $\partial\Omega\setminus(\Gamma_{top}\cup\Gamma_{bottom})$, instead, are assumed to be stress-free (homogeneous Neumann condition).
The biphasic nature of the material is made explicit in the definition of the  Lamé parameters. Indeed, while the shear modulus $\nu$ is assumed to be constant across the entire domain, the first Lamé parameter $\lambda_\theta$ is modeled as a piecewise constant function, whose definition depends on the stochastic interface angle $\theta$. More precisely, we let
\begin{equation}
    \lambda_{\theta}(\bx) = 
    \begin{cases}
        0.1 & \text{if } \bx \in \mathcal{B}_{\theta}^1, \\
        7.1 & \text{if } \bx \in \Omega \setminus \mathcal{B}_{\theta}^1,
    \end{cases}
\end{equation}
where $\mathcal{B}_{\theta}^1$ denotes the subdomain of the first material phase, effectively introducing a structural discontinuity within the physical domain.

We discretize the problem using the FEM on a structured triangular mesh ($30 \times 30$ elements - Fig.~\ref{fig:pb_setup}). The displacement field is approximated with continuous piecewise linear functions ($P1$), while material properties are modeled with piecewise constant functions ($P0$). The resulting system has $N_h = 1922$ degrees of freedom, combining $x$ and $y$ displacements.

We emphasize that, in order favor visual inspections and qualitative assessment, the units considered in 
Eq.~\eqref{eq:elasticity} have been magnified on purpose, far beyond the conventional ranges in which the equations of linear elasticity are commonly assumed to provide a reasonable approximation of real-world deformations.

\paragraph{Dataset and POD-CNF configuration.}
We leverage the FOM to construct a dataset consisting of $N=6400$ independent simulations, partitioned into training ($N_{\text{train}} = 4800$), validation ($N_{\text{val}} = 1280$) and test ($N_{\text{test}} = 320$) sets, following the same split used for the advection-diffusion problem. We sample the conditioning parameters uniformly across the parameter space, $\bmu\sim\mathcal{U}(\pspace)$, and proceed analogously for the unobserved interface angle, $\theta \sim \mathcal{U}(0,2\pi)$.
To handle the bivariate nature of the solution field, the $u_x$ and $u_y$ components are flattened and concatenated into a single vector $\bu \in \mathbb{R}^{N_h}$. For the POD, instead of fixing 
an energy threshold \emph{a priori}, we directly selected the same number of modes used for the previous test case and check their accuracy \emph{a posteriori} on the validation set. This is to guarantee an easier comparison between the two numerical experiments. Similarly, for the CNF module, we initialize the hyperparameter tuning procedure using the same setup adopted for the advection-diffusion problem. 

\begin{figure}
    \centering
    \includegraphics[width=0.6\textwidth]{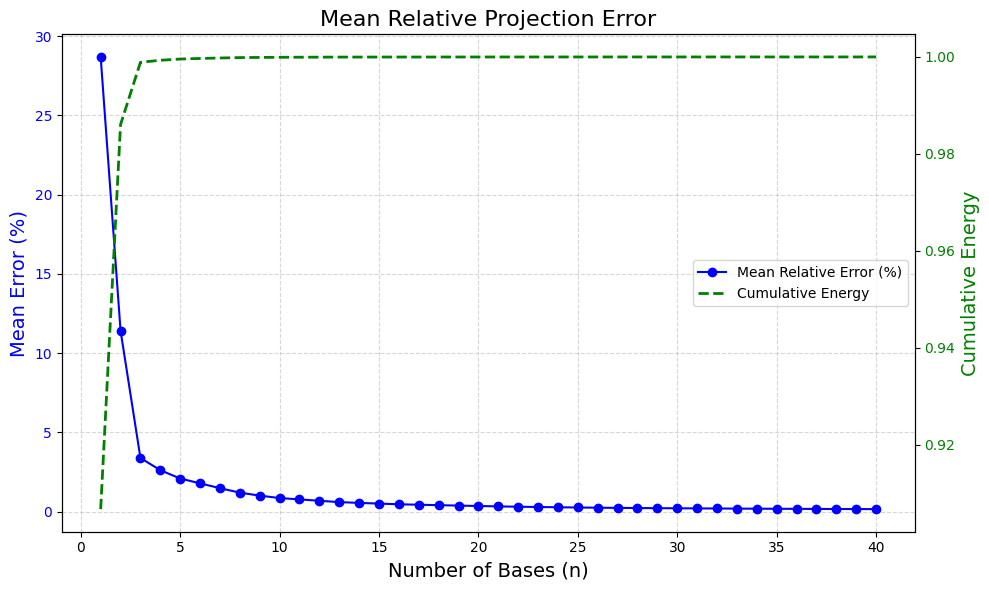}
    \caption{Linear elasticity problem, \S\ref{sec:bivariate_case}. POD truncation analysis: mean relative reconstruction error on the validation set (blue solid line) and cumulative retained energy (dashed black line) as a function of the number of POD modes $n$.}
    \label{fig:svd_elastic}
\end{figure}

\paragraph{Results I: final design of the POD-CNF surrogate.}
The SVD analysis of the training snapshots is illustrated in Fig.~\ref{fig:svd_elastic}. Retiring $n=20$ basis captures 97.15\% of the cumulative energy and achieves an average projection error of 0.36\% on the validation set, thus yielding the desired accuracy. 
Following the results of the hyperparameter tuning procedure, instead, a CNF architecture with $K=16$ coupling layers was selected, with each layer featuring conditioner MLPs of depth $\ell=2$ and width $w=256$.

%

\paragraph{Results II: model performance on a selected $\bmu_0\in\pspace$.}

Similarly to the scalar case, we first focus on a representative test configuration,
$\bm{\mu}_0 = [1.5, 0.2]$, falling right at the center of the parameter space. 

\begin{figure}
  \centering
  \parbox{0.04\textwidth}{\textbf{(a)}}%
  \parbox{0.94\textwidth}{\includegraphics[width=\linewidth]{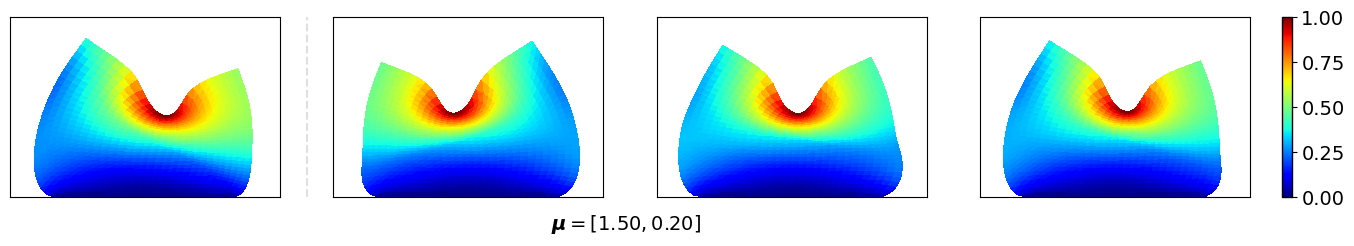}}\\[0.8em]

  \parbox{0.04\textwidth}{\textbf{(b)}}%
  \parbox{0.94\textwidth}{\includegraphics[width=\linewidth]{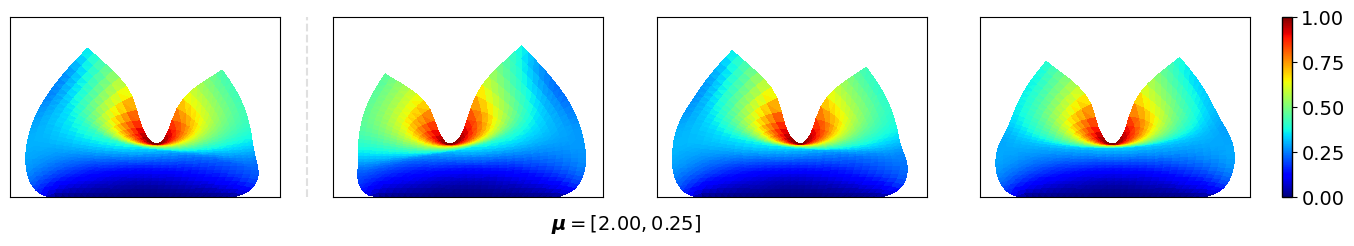}}\\[0.8em]

  \parbox{0.04\textwidth}{\textbf{(c)}}%
  \parbox{0.94\textwidth}{\includegraphics[width=\linewidth]{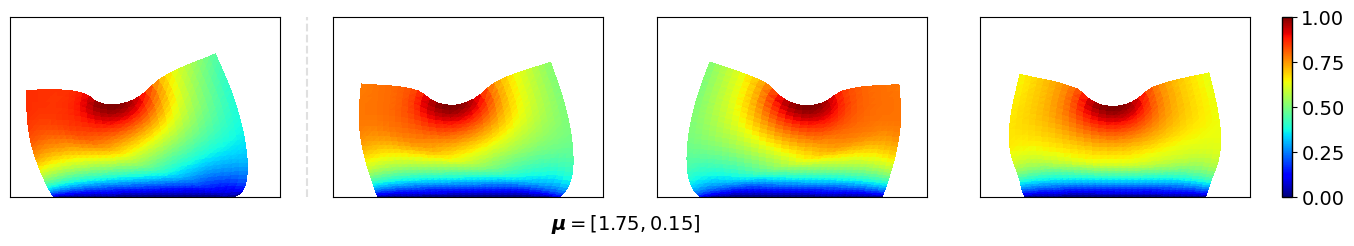}}

  \caption{Linear elasticity problem, \S\ref{sec:bivariate_case}. Solutions fields across different parametric configurations: 
    (a)~$\bm{\mu}_0 = [1.5, 0.20]$, 
    (b)~$\bm{\mu}_1 = [2.0, 0.25]$, and 
    (c)~$\bm{\mu}_2 = [2.0, 0.25]$. For each row, the leftmost panel shows a ground-truth snapshot, computed with the numerical solver, while the three panels to its right show independent realizations sampled from the POD-CNF surrogate model. Configuration~(a) corresponds to the parameter $\bmu_0\in\pspace$ chosen for the local analysis.}
  \label{fig:recon_elastic}
\end{figure}

Fig. \ref{fig:recon_elastic}a compares a random FOM solution with $\bmu=\bmu_0$ and $\theta \sim \mathcal{U}(0,2\pi)$ against three \textit{i.i.d.} realizations sampled from the trained POD-CNF model. The displacements predicted by the surrogate appear consistent with the prescribed operating condition, $\bmu_0$, capturing most of the qualitative features observed in the FOM simulation. 
At the same time, each sample is different. In particular, we can see how the model is able to reproduce the variability induced by the interface angle $\theta$, despite not being explicitly aware of it.


\begin{figure}
    \centering
    \includegraphics[width=0.8\textwidth]{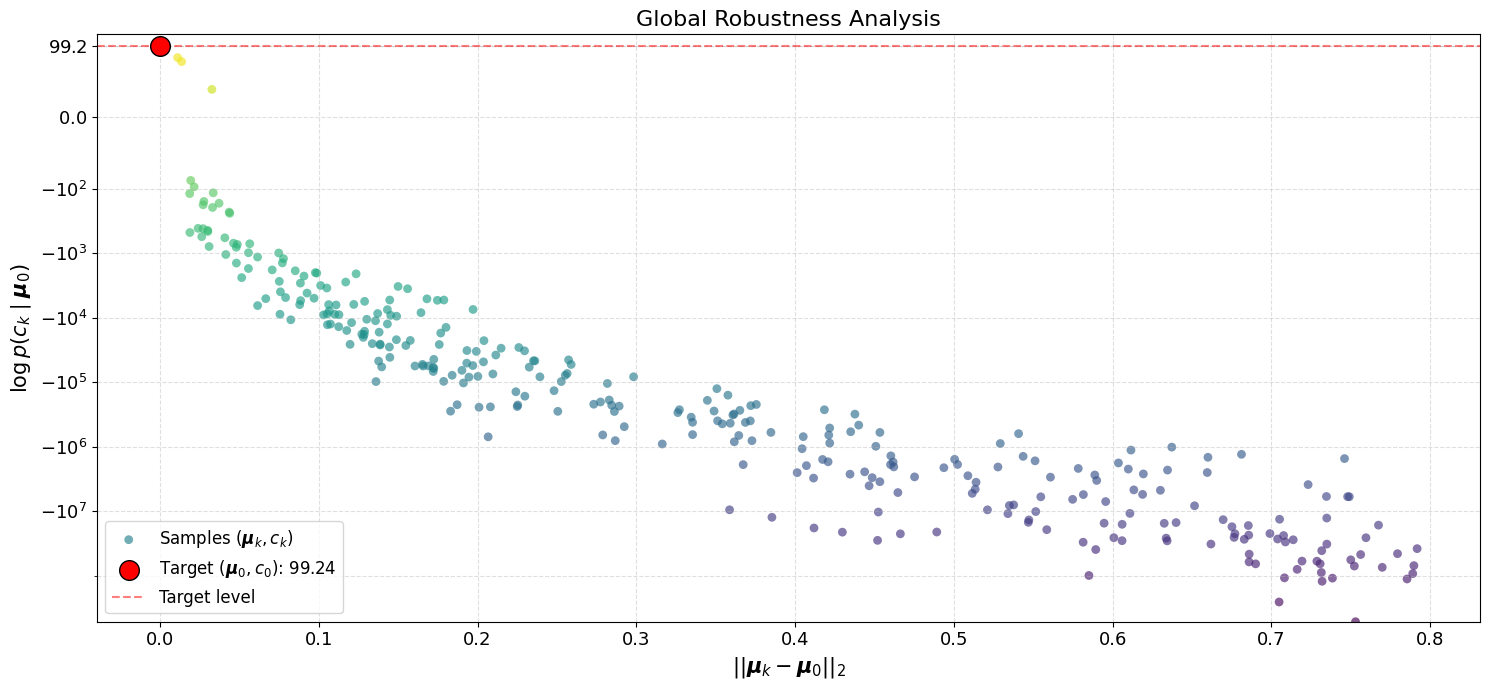}
    \caption{Reduced log-likelihood sensitivity of the POD-CNF surrogate evaluated at $\bmu_0$ against solutions belonging to different parameter instances $\bmu_k$.}
    \label{fig:glob_rob_el}
\end{figure}

Fig. \ref{fig:glob_rob_el} displays the (reduced) log-likelihood, conditioned on the chosen $\bmu_0$ and evaluated on the test set. The behavior is analogous to the one seen in the advection-diffusion problem, with the log-likelihood presenting a clearly decreasing trend in $\|\bmu-\bmu_0\|$. 
This further confirms that the stochasticity introduced by POD-CNF is highly structured, allowing the model to distinguish amongst solutions coming from different operating conditions.

\begin{figure}
    \centering
    \begin{minipage}[l]{0.49\textwidth}
        \vspace{0.1cm}\hspace{0.2cm}\includegraphics[width=0.65\linewidth]{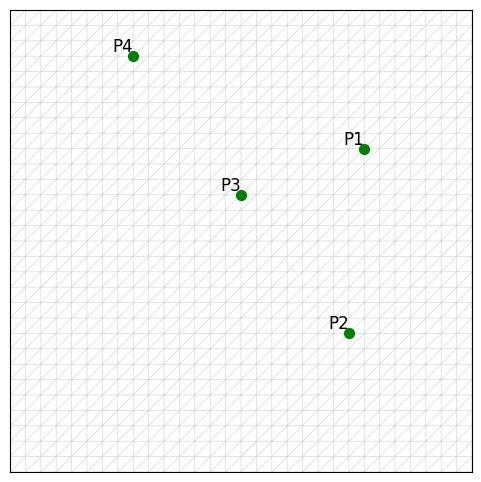}\vspace{0.22cm}\\
        \includegraphics[width=\linewidth]{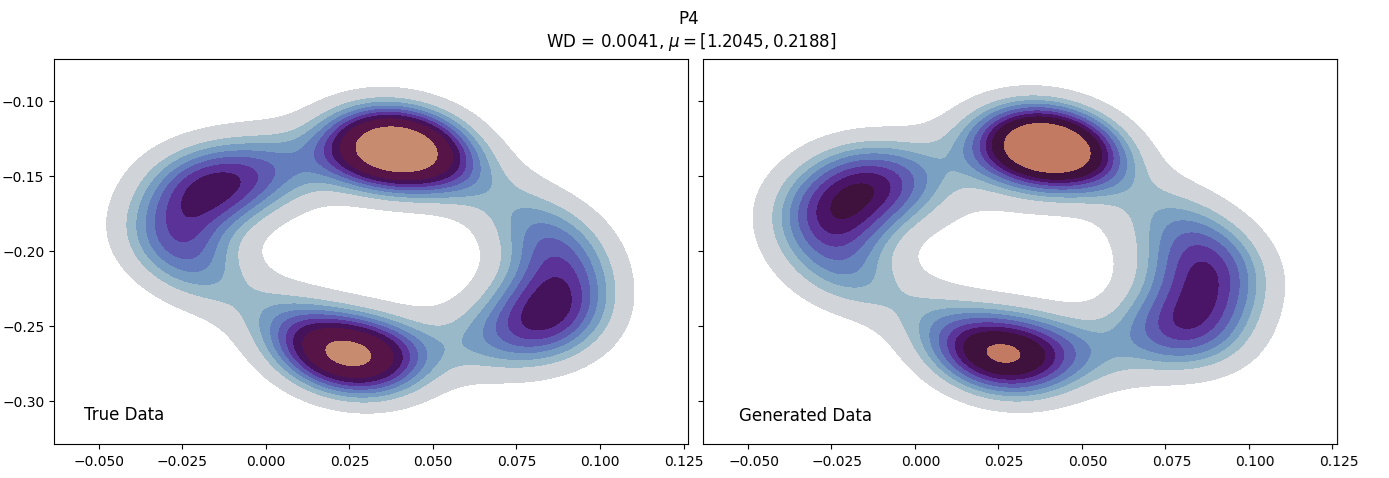}
    \end{minipage}
    \hfill
    \begin{minipage}[r]{0.49\textwidth}
        \centering
        \includegraphics[width=\linewidth]{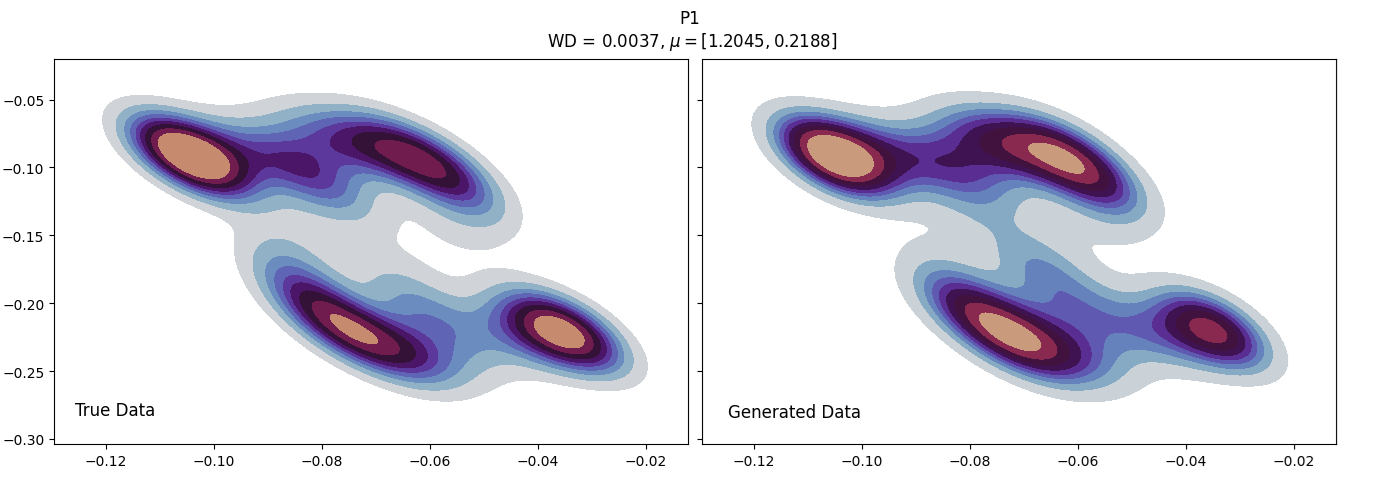}\\
        \includegraphics[width=\linewidth]{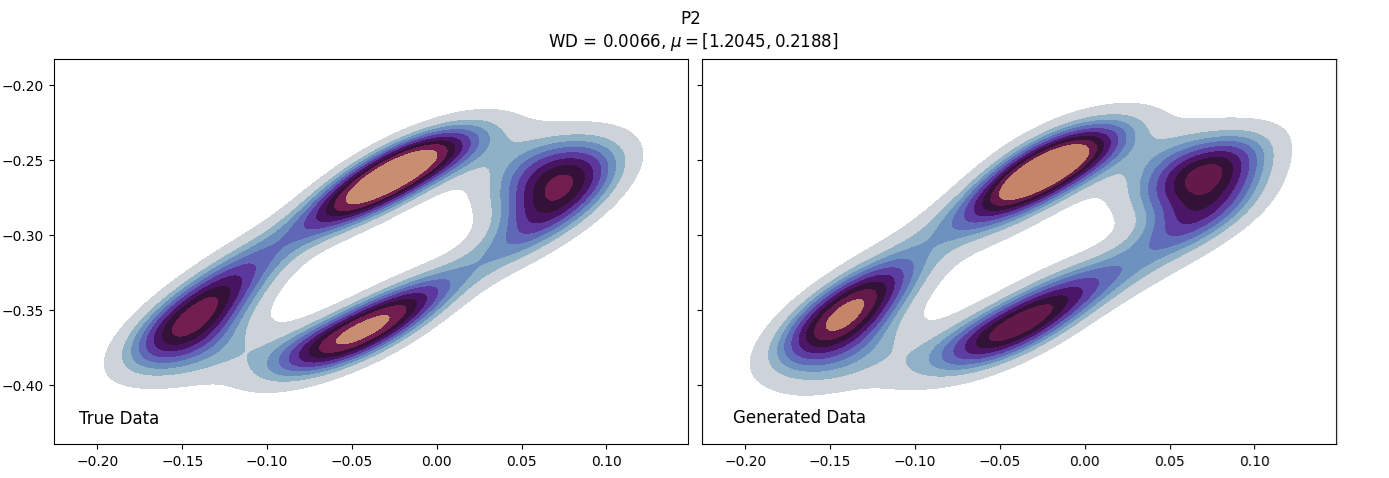}\\
        \includegraphics[width=\linewidth]{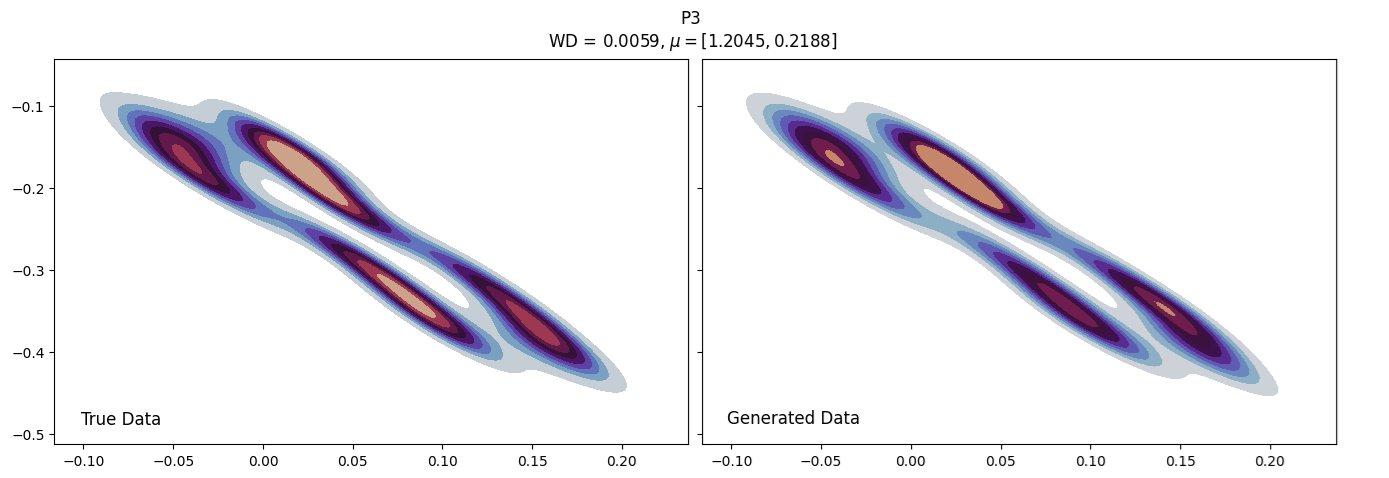}
    \end{minipage}
    \caption{Linear elasticity problem, \S\ref{sec:bivariate_case}. Top left: computational mesh and
    locations of the four monitoring points (P1--P4). Bottom and right: joint
    probability density contours of the horizontal and vertical
    displacements $(u_x, u_y)$ at P1--P4 for $\bmu_0$, comparing the
    FOM reference with the POD-CNF generated distribution
    (left and right of each subplot, respectively).}
    \label{fig:mesh_dist_el}
\end{figure}

The results of the quantitative analysis, comparing POD-CNF, POD-GP, deterministic baseline and FOM over $m_{\text{gen}} = 1000$ generated samples, can be found in Table \ref{tab:mu0-elastic}. This time,
%
%
%
For the parametric configuration under study, the proposed approach achieves a 
a Wasserstein-1 error of 0.4561, resulting nearly four times more accurate than POD-GP and improving upon the deterministic baseline by a factor of 83.76\%. 
To verify that the model correctly captures componentwise cross-correlations in the whole physical space, we monitor the 2D joint displacement ($u_x,u_u$) at four locations
across the domain. Fig. \ref{fig:mesh_dist_el} compares the resulting joint probability density contours against the FOM reference. Remarkably, POD-CNF is to able to capture several nontrivial features seen in the FOM distribution, including the characteristic multi-modal behavior induced by the interface angle $\theta$. Once again, this is particularly valuable considering that standard Gaussian approximations would generally be unable to replicate such level of complexity.

\begin{figure}
    \centering
    \includegraphics[width=0.6\textwidth]{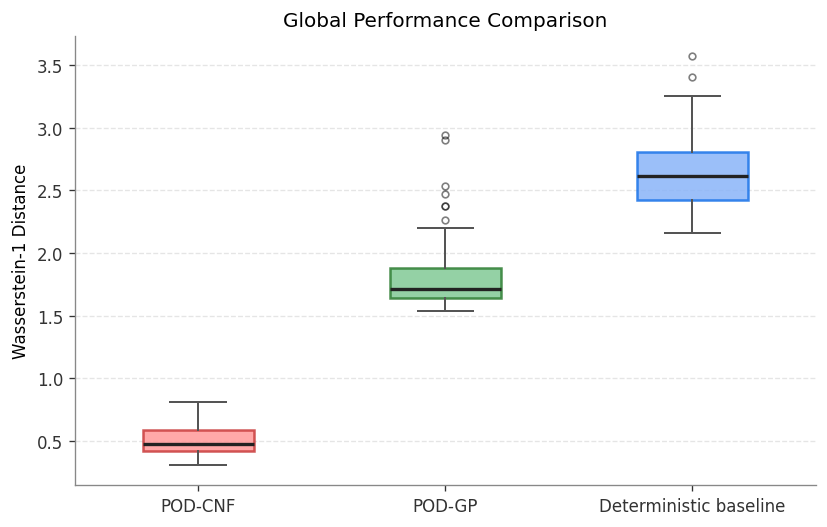}
    \caption{Linear elasticity problem, \S\ref{sec:bivariate_case}. Global analysis over 50 unseen test configurations. The boxplots show the distribution of the Wasserstein-1 error over the considered test set, comparing FOM and surrogate outputs. From left to right: POD-CNF (red, mean $ = 0.5117$), POD-GP (green, mean $ = 1.8402$) and deterministic baseline (blue, mean $ = 2.6443$).} 
    \label{fig:global_w1_elastic}
\end{figure}

\paragraph{Results III: global performance.}

We assess the generalization capabilities of POD-CNF by repeating the previous analysis on 50 unseen loading configurations extracted the test set, this time using $m_{gen}=200$ samples per instance to reduce computational cost. As shown in Figure \ref{fig:global_w1_elastic}, the Wasserstein-1 error remains well below the system's intrinsic variability (deterministic baseline). Across all test cases, POD-CNF achieves an average $\ewass$ error of $0.5517$, approximately $80\%$ lower than the baseline variability ($2.6443$). The model also performs significantly better than POD-GP (see also the negative log-likelihood in Table \ref{tab:mu0-elastic}), further demonstrating that the two approaches serve entirely different purposes, and that POD-CNF is a more natural choice in this framework.

\renewcommand{\arraystretch}{1.5}
\begin{table}\centering
    \centering
    \begin{tabular}{lllrrr}
        \hline
        & $\ewass$ at $\bmu_0$
        & Average $\ewass$ on $\mathcal{P}$
        & $\textsf{NLL}_{\text{test}}$ & Training time & Inference time\\
        \hline
        \textbf{Deterministic baseline}  & 2.8098 & 2.6443 $\pm$ 0.3308 & --  & -- & --\\
        \textbf{POD-GP}                 &  1.6707 &    1.8402 $\pm$ 0.3343
        &  \textcolor{white}{0}-27.7385 $\pm$ 0.4076 & 55.57s  & 0.29s\\
        \textbf{POD-CNF}                 & 0.4561 & 0.5117 $\pm$ 0.1259 
        & -101.9474 $\pm$ 2.3082 & 12m 24s & 0.08s\\
        \hline
    \end{tabular}
    \caption{Linear elasticity problem, \S\ref{sec:bivariate_case}. Quantitative analysis comparing POD-CNF, POD-GP and the deterministic baseline, cf. \S\ref{sec:exp_setup}. Table entries read as in Table~\ref{tab:mu0-stokes}.}
    \label{tab:mu0-elastic}
\end{table}

\paragraph{Results IV: sample size sensitivity analysis.}

We investigate the data efficiency of the architecture by retraining the CNF on progressively smaller datasets, all while keeping the POD basis ($n = 20$) and the model hyperparameters fixed (Table~\ref{tab:sens_el}). Consistently with the results seen in the previous test case, the performance remains remarkably stable down to $25\%$ of the original sample size ($N_{\text{train}} = 1200$, average $\wasserstein = 0.6768$). Reducing the training size further leads to a sharp drop in accuracy, suggesting that $N_{\text{train}} \approx 1200$ is the critical data threshold for this problem.

\begin{table}[t]\centering
    \begin{tabular}{lllll}
    \toprule
    \textbf{\% Data} & $N_{\text{train}}$ & $N_{\text{val}}$ & \textbf{Average} $\ewass$ \textbf{on} $\pspace$ & $\ewass$ \textbf{at} $\bmu_0$ \\
    \midrule
    100\% & 4800 & 1280 & 0.5117  $\pm$ 0.1259 & 0.4561 \\
    75\%  & 3600 & 960  & 0.5110  $\pm$ 0.1624 & 0.4135 \\
    50\%  & 2400 & 640  & 0.5062  $\pm$ 0.1327 & 0.4936 \\
    25\%  & 1200 & 320  & 0.6761  $\pm$ 0.1799 & 0.6768 \\
    10\%  & 480  & 128  & 1.2932  $\pm$ 0.2403 & 1.4152 \\
    5\%   & 240  & 64   & 2.3952  $\pm$ 0.2498 & 1.3998 \\
    \bottomrule
    \end{tabular}\caption{Linear elasticity problem, \S\ref{sec:bivariate_case}. Sensitivity of POD-CNF accuracy to training sample size. Table entries read as in Table \ref{stokes_sens}.} 
\label{tab:sens_el}
\end{table}

\section{Application to Forward and Inverse Uncertainty Quantification (UQ)}
\label{sec:uq}

The numerical experiments presented in Section~\ref{sec:exp} exemplify scenarios where the physical model is governed by a PDE subject to intrinsic uncertainties, encoded by latent random parameters. Having assessed the capabilities of POD-CNF in this framework, we now aim at illustrating how the latter can be deployed to tackle Uncertainty Quantification (UQ) tasks. We shall consider both forward and inverse problems, alternating between the two test cases discussed in Section~\ref{sec:exp}. 

\subsection{Forward UQ}
\label{sec:forward}
%
Traditionally, Forward UQ aims at characterizing how uncertainties about the configuration of the system propagate to a given Quantity of Interest (QoI). For instance, given a QoI $Q$ that can be expressed as a Lipschitz continuous functional of the state, $Q:\mathbb{R}^{N_h}\to\mathbb{R}$, a classical problem in Forward UQ concerns the estimation of the conditional expectation
%
\begin{equation}
\label{eq:conditional-exp}
\mathbb{E}[Q(\bu_{\bmu})] = \int Q(\bu_{\bmu})\mathsf{p}(\bu|\bmu)d\bu,
\end{equation}
often approximated via Monte Carlo, or related techniques, as
$\mathbb{E}[Q(\bu_{\bmu})] \approx \sum_{j=1}^{N}Q(\bu^{(j)}_{\bmu})/N,$ where the $\bu^{(j)}_{\bmu}$'s are independent FOM simulations, all conditional on $\bmu$. 
In our setting, we can accelerate this process, by exploiting the generative capabilities of the probabilistic surrogate, effectively bypassing the PDE solver. For instance, a natural approach is to replace the conventional estimator with
\begin{equation} \label{eq:mc_nf}
\mathbb{E}[Q(\bu_{\bmu})] \approx \frac{1}{N}\sum_{j=1}^{N}Q(\bV\phi(\bmu, \bz^{(j)})).
\end{equation}
where $\bz^{(1)},\dots\bz^{(N)}\sim\mathcal{N}(0,\mathbf{I}_n)$ are independent Gaussian random vectors.
Indeed, if we denote by $\tilde{\bu}_{\bmu}:=\mathbf{V}\phi(\bmu,\bz)$ the output of the POD-CNF, then, as a direct consequence of the Kantorovich-Rubinstein duality,
%
\begin{equation}
\label{eq:rubstein}
\left| \mathbb{E}[Q(\bu_{\bmu})] - \mathbb{E}[Q(\tilde{\bu}_{\bmu})] \right| \leq \text{Lip}(Q) \wasserstein(\bu_{\bmu}, \tilde{\bu}_{\bmu}),
\end{equation}
where $\text{Lip}(Q)$ is the Lipschitz constant of the QoI, and the $\wasserstein$-term can be bounded as discussed in Section~\ref{sec:theory}. 
Because the right-hand side of \eqref{eq:mc_nf} is nothing but the Monte Carlo approximation of $\mathbb{E}[Q(\tilde{\bu}_{\bmu})]$, this readily motivates the definition of the proposed estimator.
%
\\\\
To demonstrate how the idea unfolds in practice, we present a simple application featuring the advection-diffusion problem discussed in Section \ref{sec:scalar_case}. We consider five different QoIs, defined as in Table \ref{tab:quantities-of-interest}: each of them relates to a different feature of the solution field. For instance, $Q_3$ quantifies the amount of solute  accumulated around the bottom-right cylinder, in comparison with the other two. Conversely, $Q_5$ identifies potential concentration spikes in the wake region of the mixer.

\begin{figure}
 \centering
   \includegraphics[width=\textwidth]{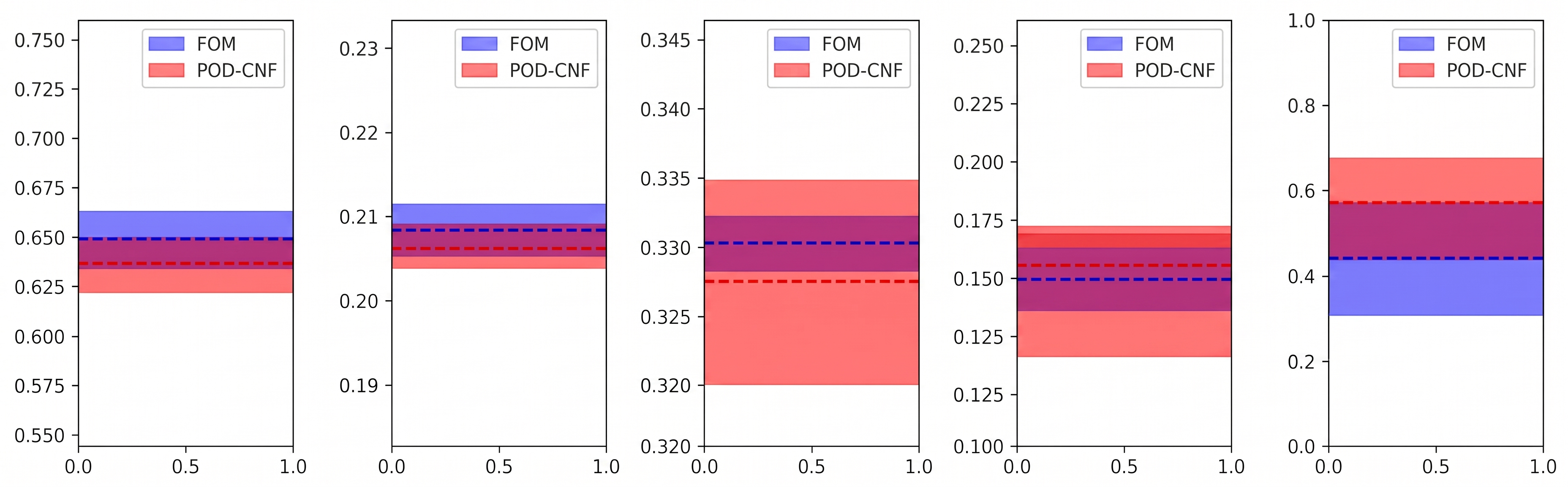}
   \caption{Comparison of the $99\%$ confidence intervals obtained using the true FOM solver and the proposed POD-CNF surrogate for the evaluated Quantities of Interest (QoIs).}\label{fig:forward_res}
\end{figure}

Figure \ref{fig:forward_res} shows the results obtained for the different QoIs when approximating \eqref{eq:conditional-exp} using $N=100$ Monte Carlo samples. In particular, the picture reports the estimate mean values, together with their 99\% confidence interval. Notably, the estimates obtained with POD-CNF perfectly agree with the ones obtained using the FOM: the two approaches, in fact, systematically yield overlapping confidence intervals. Moreover, the width of the latter is often comparable between the two approaches, further confirming the ability of POD-CNF in replicating variability of the FOM. Still, there are also major differences amongst the different QoIs. For instance, POD-CNF clearly overestimates the uncertainties in $Q_3$ and $Q_4$. We speculate that this is related to the very definition of the QoIs: some of them, in fact, correspond to nonlinear functionals that are not continuous in the $L^2$-norm. This makes \eqref{eq:rubstein} formally invalid (at least in the limit $N_h\to+\infty$), thus hindering the advantages of probabilistic surrogates  in general.

In principle, because the proposed surrogate is extremely efficient at inference time, we could consider improving our estimate by increasing the number of Monte Carlo samples in \eqref{eq:mc_nf}. However, we notice that this is only suitable up to a point, which is when the statistical error---related with the amplitude of the confidence interval---becomes comparable with the magnitude of the model error in \eqref{eq:rubstein}. Circumventing this problem requires the development of suitable multi-fidelity strategies, specifically designed to  balance these two sources of inaccuracy: see, e.g., \cite{peherstorfer2016optimal, farcaș2023context}. These methods, however, have been mostly developed for deterministic surrogates: extending the idea to probabilistic approaches is in order, but beyond the scope of this work. 

Nonetheless, it can be insightful to estimate the potential reduction in computational cost. Table \ref{tab:forward-results} shows the hypothetical running costs required to approximate \eqref{eq:conditional-exp} with a relative uncertainty, at 99\% confidence, below 1\%. The latter have been estimated by leveraging the previous pilot run and inverting the confidence interval formula; this made it possible to estimate the required sample size $N^*$ and the corresponding running costs. According to such projections, POD-CNF has the potential to reduce computational times by three orders of magnitude. This motivates developing the idea further, extending beyond this preliminary exploration.

\begin{table}\centering
\centering\small
\renewcommand{\arraystretch}{2.5}
\begin{tabular}{@{}clll@{}}
\toprule
\textbf{Label} & \textbf{QoI} & \textbf{Expression}
& \textbf{Continuous counterpart} \\
\midrule

$Q_1$
& Average Concentration Gradient
& $
  \frac{1}{N_c}\sum_{j=1}^{N_c}
  \|\nabla\mathrm{u}^{(j)}_{\bmu}\|_2$
& $
  \frac{1}{|\Omega|}\int_{\Omega}
  \|\nabla u_{\bmu}\|\,d\bx$
\\

$Q_2$
& Spatial Standard Deviation
& $
  \left[
  \frac{1}{N_h}\sum_{j=1}^{N_h}
  |\mathrm{u}^{(j)}_{\bmu}-\overline{u}_{\bmu}|^2
  \right]^{1/2}$
& $
  \left[
  \frac{1}{|\Omega|}\int_{\Omega}
  |u_{\bmu}-\overline{u}_{\bmu}|^2\,d\bx
  \right]^{1/2}$
\\

$Q_3$
& Relative Boundary Concentration
& $
  \left[\sum_{j\in I_{c_3}}|\mathrm{u}^{(j)}_{\bmu}|\right]/\left[\sum_{i=1}^{3}\sum_{j\in I_{c_i}}|\mathrm{u}^{(j)}_{\bmu}|\right]$
& $  \left[\int_{\Gamma_{c_3}}|u_{\bmu}|\,ds\right]/\left[\sum_{i=1}^{3}\int_{\Gamma_{c_i}}|u_{\bmu}|\,ds\right]$
\\

$Q_4$
& High-Concentration Area Fraction
& $
  \frac{1}{N_h}\sum_{j=1}^{N_h}
  \indicator_{(0.8,+\infty)}(\mathrm{u}^{(j)}_{\bmu})$
& $
  \frac{1}{|\Omega|}\int_{\Omega}
  \indicator_{(0.8,+\infty)}(u_{\bmu})\,d\bx$
\\

$Q_5$
& Localized High-Concentration Event
& $
  \indicator_{(0.8,+\infty)}\big(
  \max_{j\in I_{\mathrm{right}}}\mathrm{u}^{(j)}_{\bmu}
  \big)$
& $
  \indicator_{(0.8,+\infty)}\big(
  \sup_{\bx\in\Omega_{\mathrm{right}}}u_{\bmu}(\bx)
  \big)$
\\

\bottomrule
\end{tabular}%
\caption{QoIs for the Forward UQ analysis. $|\Omega| = $ area of the spatial domain. $N_h = $ number of vertices in the mesh. $N_c = $ number of cells in the mesh. $\mathrm{u}^{(j)} = $ $j$th component of $\bu_{\bmu}$. $\nabla\mathrm{u}^{(j)} = $ value of the gradient of the FEM interpolant at the barycenter of the $j$th cell.  $\overline{u}_{\bmu} = $ the nodal average (resp. integral) average of the concentration field (third and fourth columns, respectively). $\Gamma_{c_i}$ and $I_{c_i} = $ continuous boundary and discrete degrees of freedom associated with the $i$-th cylinder obstacle, respectively. $\Omega_{\text{right}} = \{ \bx \in \Omega : x_1 > 2.25 \}$. $I_{\text{right}} = $ degrees of freedom in the mesh belonging to $\Omega_{\text{right}}$.}
\label{tab:quantities-of-interest}
\end{table}

\begin{table*}
    \begin{minipage}[t]{0.48\textwidth}
        \begin{center}
            \textbf{POD-CNF}
        \end{center}
        \begin{tabular*}{\linewidth}{@{\extracolsep{\fill}}l c c c c@{}}
        \toprule
        \textbf{QoI} & \textbf{Est. Mean} & \textbf{99\% CI} & \textbf{$N^*$} & \textbf{Time}$^*$ \textbf{(s)} \\ 
        \midrule
        $Q_1$ & 0.6360 & 0.0138 & 472   & 0.8349 \\
        $Q_2$ & 0.2061 & 0.0040 & 376   & 0.6651 \\
        $Q_3$ & 0.3275 & 0.0073 & 501   & 0.8862 \\
        $Q_4$ & 0.1553 & 0.0248 & 25500 & 45.1065 \\
        $Q_5$ & 0.5700 & 0.1310 & 52564 & 92.9796 \\
        \bottomrule
        \end{tabular*}
    \end{minipage}
    \hfill
    \begin{minipage}[t]{0.48\textwidth}
    \begin{center}
            \textbf{FOM}
        \end{center}
        \begin{tabular*}{\linewidth}{@{\extracolsep{\fill}}l c c c c@{}}
        \toprule
        \textbf{QoI} & \textbf{Est. Mean} & \textbf{99\% CI} & \textbf{$N^*$} & \textbf{Time}$^*$ \textbf{(s)} \\ 
        \midrule
        $Q_1$ & 0.6486 & 0.0140 & 464   & 478.5703 \\
        $Q_2$ & 0.2083 & 0.0030 & 213   & 219.6885 \\
        $Q_3$ & 0.3302 & 0.0020 & 36    & 37.1305 \\
        $Q_4$ & 0.1495 & 0.0132 & 7771  & 8015.0204 \\
        $Q_5$ & 0.4400 & 0.0248 & 88680 & 91464.6777 \\
        \bottomrule
        \end{tabular*}
    \end{minipage}
    \caption{Summary of the Forward UQ analysis for POD-CNF (left) and FOM (right). Est. Mean = value of the Monte Carlo estimator with $N=100$ samples. 99\% CI = semi-amplitude of the 99\% confidence interval with $N=100$ Monte Carlo samples. $N^* = $ estimated sample size required to reduce the relative uncertainty below 1\%. Time$^*$ = running time (in seconds) required by the model to generate $N^*$ samples.}\label{tab:forward-results}
\end{table*}

\subsection{Inverse UQ}
\label{sec:inverse}

In many applications, one is not only interested at predicting the response of a physical system for a prescribed problem configuration (\emph{forward problem}), but rather aims at inferring the values of the problem parameters given, incomplete, experimental observations of the system itself (\emph{inverse problem}). Because inverse problems are typically ill-posed, a common strategy is to tackle them using probabilistic methods based on Bayesian inference \cite{gelman1995bayesian, stuart2010}. In essence, rather than providing a pointwise estimate, these approaches leverage system measurements to produce a full probability distribution over the parameter space (the \emph{posterior distribution}), thus quantifying uncertainties and providing additional insights.

In the case of intrinsic, structured stochasticities produced by hidden parameters, this task becomes particularly challenging. In fact, traditional approaches generally require modeling aleatoric uncertainty explicitly. This is perfectly doable when such uncertainty is purely associated with measurements noise, which can be reasonably described using additive models based on Gaussian assumptions. However, when uncertainties arise from \emph{within} the system, the situation changes substantially, adding an additional layer of complexity.
\\\\
To illustrate the idea, we shall present a simple example based on the linear elasticity problem discussed in Section \ref{sec:bivariate_case}. Recall that, there, $\bu_{\bmu}$ represented the vector field associated with the local displacement of a solid body subject to an external load. The magnitude and the spatial extension of such load was determined by two parameters, $m$ and $\delta$, respectively, collected as $\bmu=[m,\delta]$. The intrinsic uncertainty was generated by a hidden variable, $\theta\sim\mathcal{U}(0,2\pi)$, controlling the interface angle between two material layers constituting the body.

We consider the problem of inferring the values of $m$ and $\delta$ given sensor measurements of the displacement vector field at 31 equally spaced sensors located at the top surface of the body. We shall write $\bu_{\text{sens}}\in\mathbb{R}^{62}$ to the denote the data vector collecting such measurements. Assuming that all stochasticites arise from the hidden variable $\theta$, whose value is unknown, we shall neglect potential measurement noise affecting the sensors. This is depart from the conventional setting and focus entirely on the new challenges posed by the considered framework. 

Following a Bayesian paradigm, we equip the parameter space with a prior probability measure $\density(\bmu)$, chosen to be the uniform distribution over $\pspace.$ According to Bayes' rule, the posterior distribution $\density(\bmu \mid \bu_{\text{sens}})$ satisfies
\begin{equation}
    \label{eq:bayes}
    \density(\bmu \mid \bu_{\text{sens}}) \propto \density(\bu_{\text{sens}} \mid \bmu) \density (\bmu).
\end{equation}
In standard settings, the likelihood function $\density(\bu_{\text{sens}} \mid \bmu) $ can be readily derived from the model used to described measurement uncertainties. In this case, however, we cannot express the likelihood function explicitly, as the latter depends intrinstically on the PDE solver and its interaction with the hidden variable $\theta\sim\mathcal{U}(0,2\pi)$.

However, we notice that, because our trained POD-CNF  is able to approximate the conditional law of $\bu\mid\bmu$, we can leverage this knowledge to surrogate $\density(\bu_{\text{sens}} \mid \bmu).$ To see this, let $\mathbf{P}_{\text{sens}}\in\mathbb{R}^{62\times N_h}$ the matrix mapping full displacement vectors onto sensor measurements. Then
\begin{equation}
\label{eq:sensors}
\bu_{\text{sens}}\mid\bmu\;\;=\;\;\left(\mathbf{P}_{\text{sens}}\bu\right)\mid\bmu\;\;\overset{d}{\approx}\;\;\mathbf{P}_{\text{sens}}\mathbf{V}\phi(\bmu,\bz)=:\mathbf{V}_{\text{sens}}\phi(\bmu,\bz),\end{equation}
where $\overset{d}{\approx}$ denotes approximation in distribution and $\mathbf{V}_{\text{sens}}:=\mathbf{P}_{\text{sens}}\mathbf{V}.$ While \eqref{eq:sensors} does not suffice to approximate $\density(\bu_{\text{sens}} \mid \bmu)$ explicitly, it can be used to estimate the distribution numerically. In fact, because $\phi$ is computationally inexpensive to evaluate, the law of $\mathbf{V}_{\text{sens}}\phi(\bmu,\bz)$, and thus that of $\bu_{\text{sens}}\mid\bmu$, can be estimated empirically by drawing $N_{gen}$ \emph{i.i.d.} samples 
and fitting a Kernel Density Estimation (KDE) model. This procedure can be then integrated within a suitable Markov chain Monte Carlo (MCMC) method to approximate the posterior distribution using \eqref{eq:bayes}.

\begin{figure}
    \centering
    \setlength{\tabcolsep}{2pt}
    \begin{tabular}{c c}
        \multicolumn{2}{c}{\small (a) $m=1.7384, \delta=0.2295$}\\[2pt]
        \includegraphics[height=0.168\textheight]{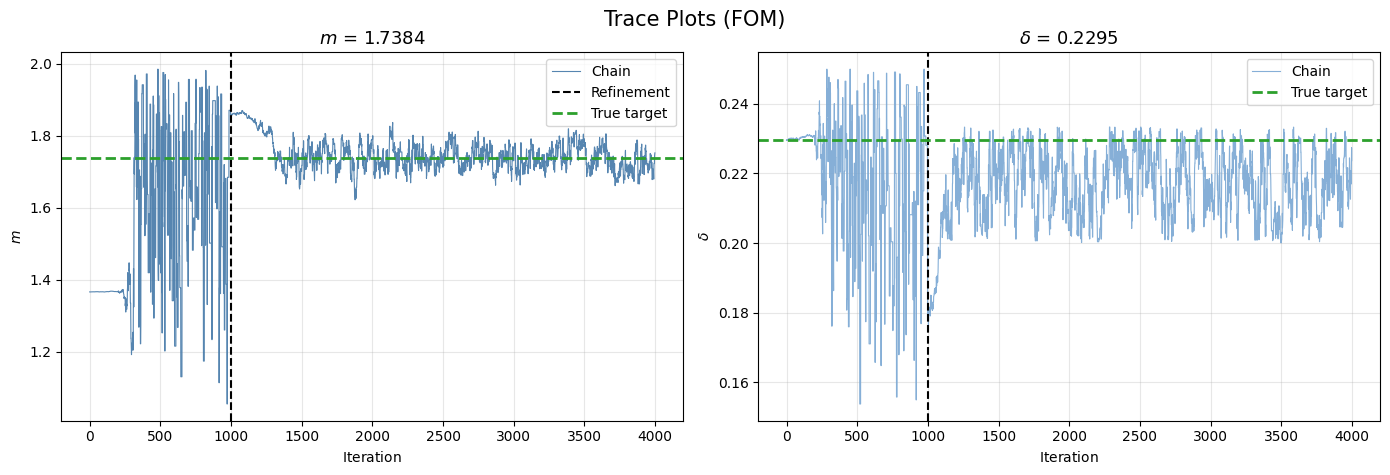} &
        \includegraphics[height=0.16\textheight]{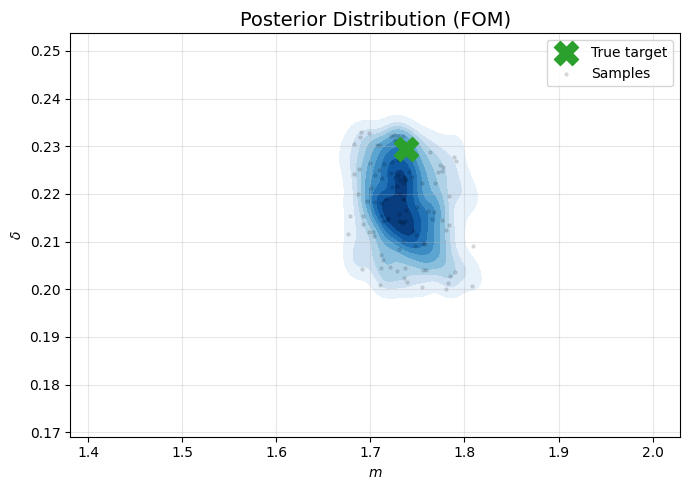} \\
        \includegraphics[height=0.168\textheight]{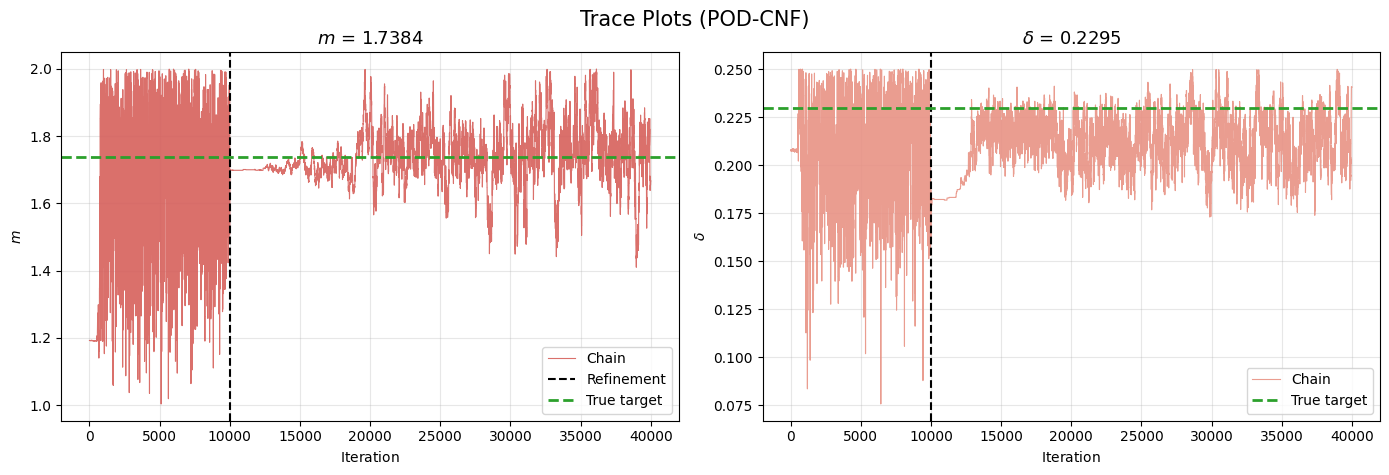} &
        \includegraphics[height=0.16\textheight]{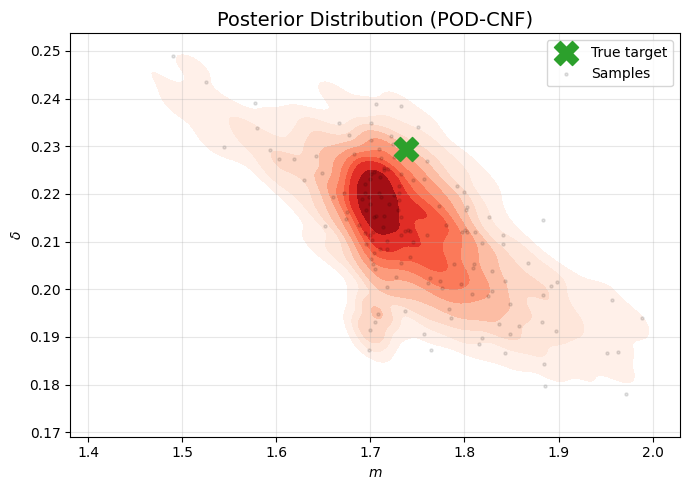} \\[6pt]
        \multicolumn{2}{c}{\small (b) $m=1.1656, \delta=0.1905$}\\[2pt]
        \includegraphics[height=0.168\textheight]{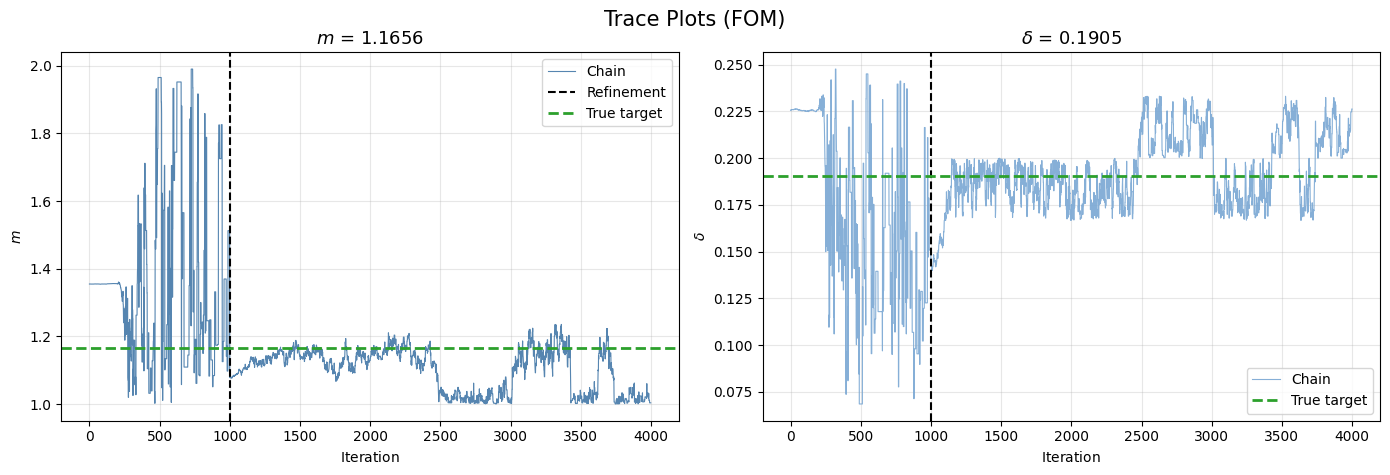} &
        \includegraphics[height=0.16\textheight]{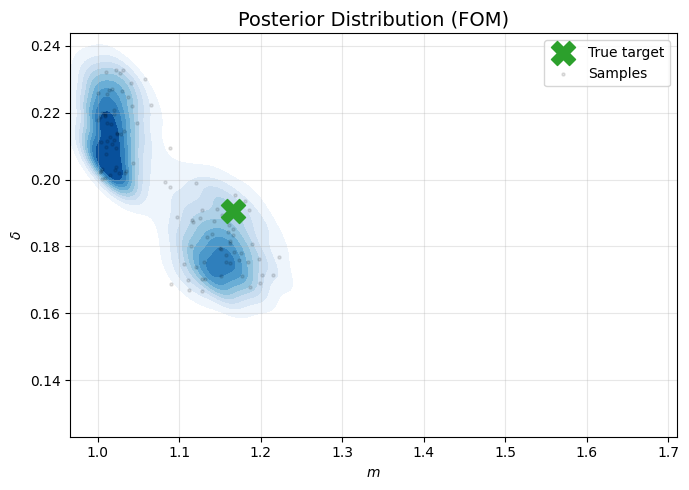} \\
        \includegraphics[height=0.168\textheight]{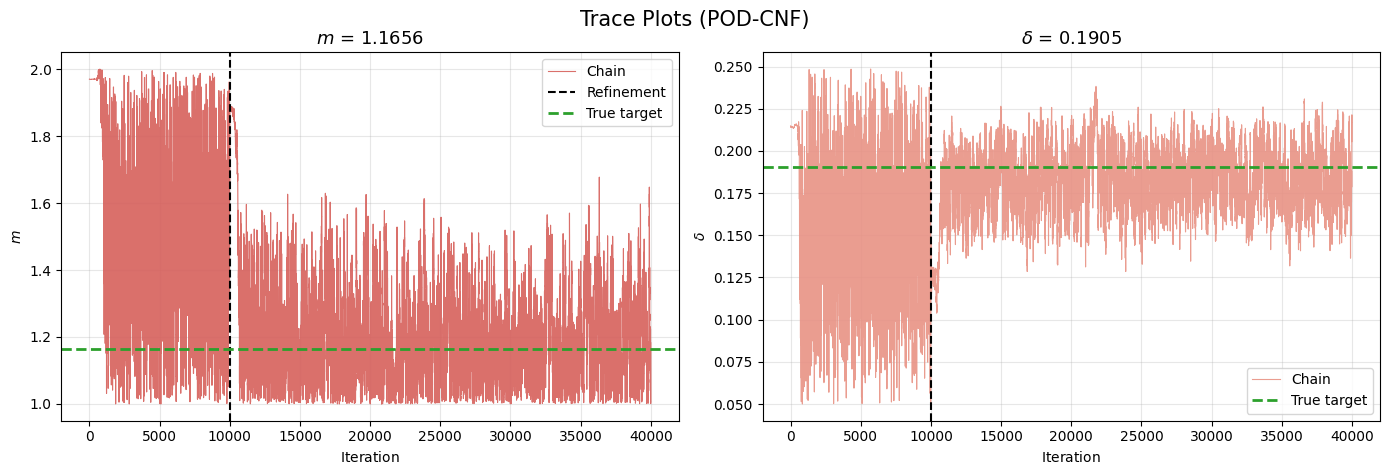} &
        \includegraphics[height=0.16\textheight]{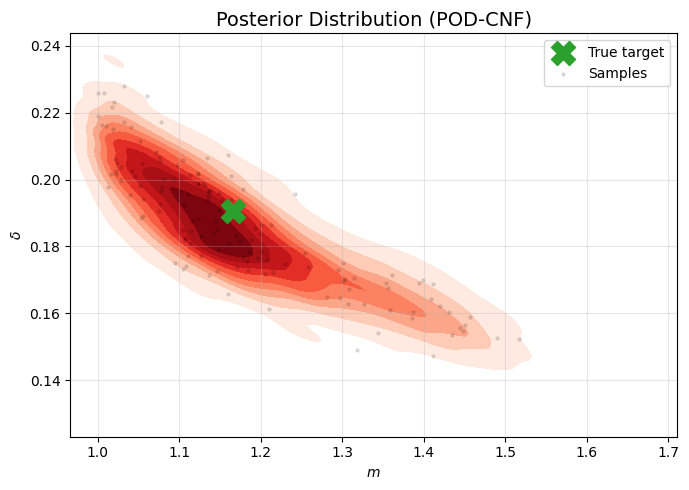} \\
    \end{tabular}
    \caption{Linear elasticity problem, \S\ref{sec:bivariate_case}. For both panels, the first row reports the MCMC traces and posterior distribution $\density(\bmu \mid \bu_{\text{sens}})$ from the high-fidelity numerical solver (FOM baseline), while the second row shows the corresponding results obtained via the POD-CNF surrogate.}
    \label{fig:inv_res}
\end{figure}

Figure \ref{fig:inv_res} shows the results obtained using this strategy in two different parametric scenarios. The pipeline has been implemented using an adaptive Metropolis-Hastings algorithm \cite{haario2001adaptive}, complemented by a two-phase strategy based on an initial posterior exploration (chain with higher \emph{temperature} and fixed proposal) followed by a subsequent refinement phase (lower temperature and adaptive proposal). During the iterations of the MCMC algorithm, the likelihood function has been estimated using the aforementioned KDE approximation, here implemented using a Gaussian kernel. For the latter, we used $N_{gen}=100$ and a bandwidth of $h=0.1$ during the exploratory phase of the MCMC, and $N_{gen}=200$ with $h=0.03$ during the refinement stage.

As show in Figure \ref{fig:inv_res} and Table \ref{inv_tbl1} the proposed strategy successfully estimates the loading parameters $m$ and $\delta$, returning a posterior concentrated near the ground truth values. The approximation is more accurate for $m$ (pointwise error below 2\%) and slightly less for $\delta$ (pointwise errors around 3\%-7\%, but with true values always falling in the proposed credibility intervals).

To corroborate the analysis, we repeated the pipeline using the FOM solver, but introducing suitable modifications that would keep the overall computational cost comparable to that of POD-CNF. Specifically, we reduced the number of iterations in the MCMC by an order of magnitude (1000 iterations for the exploratory phase and 3000 for the refinement stage), as well as the number of samples generated per iteration ($N_{gen}=10,25$ in the two phases). Even with these modifications, the FOM analysis remained substantially more intensive than the previous one (approximately one hour versus 1–6 minutes), but was at least feasible to run.

\begin{table}\centering
\begin{tabular}{lcccrccccr}
\hline
\textbf{True $m$} & \textbf{Est. $m$} & \textbf{95\% CI} & \textbf{Error} & & \textbf{True $\delta$} & \textbf{Est. $\delta$} & \textbf{95\% CI} & \textbf{Error} & \textbf{Time} \\
\midrule
1.7384 & 1.7412 & $[1.5438, 1.9414]$ & 0.16\% & & 0.2295 & 0.2127 & $[0.1869, 0.2385]$ & 7.34\% & 42s \\
1.6228 & 1.5195 & $[1.4951, 1.7095]$ & 6.37\% & &  0.1653 & 0.2088 & $[0.1449, 0.2270]$ & 26.35\% & 46s \\
1.1777 & 1.3440 & $[1.0607, 1.4985]$ & 14.12\% & &  0.2032 & 0.1665 & $[0.0904, 0.2408]$ & 18.08\% & 45s \\
1.8151 & 1.8463 & $[1.7462, 1.9728]$ & 1.72\% &  & 0.2354 & 0.2380 & $[0.2220, 0.2493]$ & 1.12\% & 6m 20s \\
1.2920 & 1.3550 & $[1.1414, 1.5936]$ & 4.88\% &  & 0.2278 & 0.2109 & $[0.1703, 0.2441]$ & 7.42\% & 6m 15s \\
1.0757 & 1.1165 & $[1.0077, 1.2689]$ & 3.79\% &  & 0.1513 & 0.1463 & $[0.1198, 0.1734]$ & 3.27\% & 42s \\
1.1656 & 1.1863 & $[1.0111, 1.4722]$ & 1.78\% &  & 0.1905 & 0.1838 & $[0.1509, 0.2177]$ & 3.52\% & 41s \\
1.2045 & 1.2320 & $[1.0326, 1.4742]$ & 2.29\% &  & 0.2188 & 0.2132 & $[0.1740, 0.2453]$ & 2.53\% & 6m 3s \\
\hline
\end{tabular}\caption{Summary of the inverse problem results for the POD-CNF surrogate model. CI = marginal credibility interval.}\label{inv_tbl1}
\end{table}

\begin{table}[tb]\centering
\begin{tabular}{lcccrccccr}
\hline
\textbf{True $m$} & \textbf{Est. $m$} & \textbf{95\% CI} & \textbf{Error} &  & \textbf{True $\delta$} & \textbf{Est. $\delta$} & \textbf{95\% CI} & \textbf{Error} & \textbf{Time} \\
\midrule
1.7384 & 1.7383 & $[1.6788, 1.8011]$ & 0.01\% & & 0.2295 & 0.2178 & $[0.2012, 0.2325]$ & 5.08\% & 63m 20s \\
1.6228 & 1.4549 & $[1.3883, 1.5074]$ & 10.35\% & & 0.1653 & 0.1835 & $[0.1676, 0.1994]$ & 11.01\% & 61m 40s \\
1.1777 & 1.3873 & $[1.3118, 1.4656]$ & 17.80\% & & 0.2032 & 0.1829 & $[0.1673, 0.1989]$ & 10.01\% & 65m 53s \\
1.8151 & 1.8221 & $[1.7673, 1.8757]$ & 0.39\% & & 0.2354 & 0.2414 & $[0.2338, 0.2496]$ & 2.57\% & 65m 0s \\
1.2920 & 1.8346 & $[1.7432, 1.9279]$ & 41.99\% & & 0.2278 & 0.1497 & $[0.1342, 0.1656]$ & 34.28\% & 60m 50s \\
1.0757 & 1.9002 & $[1.7527, 1.9978]$ & 76.64\% & & 0.1513 & 0.0846 & $[0.0678, 0.0993]$ & 44.06\% & 63m 54s \\
1.1656 & 1.0948 & $[1.0016, 1.2011]$ & 6.07\% & & 0.1905 & 0.1952 & $[0.1682, 0.2304]$ & 2.48\% & 57m 12s \\
1.2045 & 1.7189 & $[1.6497, 1.8027]$ & 42.71\% & & 0.2188 & 0.1501 & $[0.1345, 0.1654]$ & 31.40\% & 62m 30s \\
\hline
\end{tabular}\caption{Summary of the inverse problem results for the FOM baseline. CI = marginal credibility interval.}\label{inv_tbl2}
\end{table}

As shown in Figure \ref{fig:inv_res} the results are consistent with the previous analysis. In the first scenario, the two approaches yield posterior distributions localized around the same maximum value, although their spread differs substantially. Such difference could be caused by an intrinsic discrepancy between the two or by the different hyperparameters used for the MCMC algorithm. In the second scenario, instead, the two distributions identify a similar high probability region, presenting an elongated shape that extends towards the top-right corner of the parameter space.

Tables \ref{inv_tbl1}-\ref{inv_tbl2} provide a broader overview, extending the comparison to other parametric scenarios. When both models demonstrate strong predictive performance, their 95\% credibility intervals typically exhibit a significant overlap. In particular, the intervals generated by POD-CNF surrogate tend to intersect with the tighter ones produced by the FOM, indicating an agreement between the high-probability regions predicted by the two approaches.
At the same time, the pointwise estimate produced by the FOM are not always reliable, and the tightness of such intervals turns out to be unjustified: see, e.g., rows 5, 6 and 8 in Table \ref{inv_tbl2}. In these challenging cases, the FOM yields relative errors on $m$ and $\delta$ exceeding 40\% and 30\%, respectively. In contrast, the adaptive Metropolis-Hastings algorithm driven by the POD-CNF surrogate remains stable, successfully recovering true parameter values in most tested configurations at a fraction of the computational time.

The proposed analysis, however, is not free from limitations and only constitutes a preliminary exploration. Indeed, the necessity to approximate the likelihood function empirically and the subsequent the reliance on KDE currently pose major limitations, hindering the scalability of the approach to more complex inverse problems. Similarly, the comparison between the proposed approach and the FOM reference is not entirely exhaustive. These aspects could be addressed by consider more refined inference strategies that combine ideas from Bayesian statistics and Simulation Based Inference \cite{cranmer2020frontier}. Still, because Inverse UQ is an application of interest but not the main focus of this work, we leave these considerations for future explorations.

\section{Conclusions}
\label{sec:conclusions}

In this paper, we presented POD-CNF, a hybrid generative framework designed to tackle the 
challenges associated with hidden parameters in high-dimensional systems governed by stochastic PDEs, and the corresponding Uncertainty Quantification (UQ) tasks. By integrating linear dimensionality reduction via Proper Orthogonal Decomposition (POD) with the expressive density estimation capabilities of Conditional Normalizing Flows, we developed a tool capable of modeling the full conditional probability distribution of the physical field, $\density(\mathbf{u} | \bm{\mu})$, thus accounting for the aleatoric uncertainty induced by the unobserved latent variables.

The numerical experiments demonstrated the effectiveness of the proposed approach, further consolidating our theoretical analysis. A key indicator of success was the Wasserstein-1 distance ($\wasserstein$) between the generated and ground-truth distributions. In all test cases, the global $\wasserstein$ error was significantly lower than the one obtained by the benchmark models, 
including POD-GP and deterministic baselines, confirming 
the ability of POD-CNF to capture complex, non-Gaussian marginals rather than merely regressing a deterministic average. The data efficiency of the proposed strategy was assessed through a dedicated sensitivity analysis, showing robust performances even when reducing the number of training samples by a significant amount (-75\% with respect to the original sample size).


A major advantage of the POD-CNF framework lies in 
its generative capabilities.
These allowed us to seamlessly extend the method to tackle not only forward UQ tasks, but also inverse UQ problems related with parameter estimation. By embedding the pre-trained generative model within a Bayesian inference framework, we were able to approximate 
the posterior distribution of the physical parameters $\bm{\mu}$
through a specialized Markov chain Monte Carlo sampling.

Despite these promising results, the current methodology also presents some limitations:

\begin{itemize}
    \item 
    relying on POD implies a linear approximation of the solution manifold, which may struggle with, e.g., advection-dominated phenomena or moving discontinuities, especially if considering time-dependent problems. Future developments could replace POD with nonlinear dimensionality reduction techniques, such as variational autoencoders (VAEs) \cite{rezende2015}, which can compress these complex manifolds more effectively;
    

    \item the proposed method for constructing probabilistic surrogates, 
    and tackling UQ tasks, is purely data-driven; it learns the conditional distribution of the POD coefficients, yielding physically consistent outputs (as demonstrated by our numerical experiments), but without enforcing the underlying governing equations. Future work could incorporate physics-informed regularization or even hard constraints \cite{mucke2026physics}: for instance, by adding to the loss function a penalty term involving the PDE residual $\|\mathcal{R}(\hat{\mathbf{u}})\|^2$, the model could be encouraged to generate physically consistent solutions without strictly requiring a solver during inference;
    
    \item 
    our preliminary exploration concerning forward and inverse UQ tasks remains subject to limitations. In the first case, when estimating conditional expectations of given quantities of interest, increasing the number of Monte Carlo samples reduces the sampling error but it does not address surrogate model error, motivating the development of suitable multi-fidelity strategies for probabilistic surrogates. In the inverse UQ analysis, instead, the major limitation is given by the empirical approximation of the likelihood function via Kernel Density Estimation, which currently limits the scalability of the approach to more complex problems.
\end{itemize}

On top of this, a significant avenue for future research involves extending our framework to time-dependent problems. 
This would require modeling the temporal evolution of the latent coefficients $\mathbf{c}(t)$ by learning either their transition probabilities or the distribution of their trajectories, ultimately enabling the prediction of time-dependent quantities of interest under uncertainty.

\subsection*{Acknowledgments}
\noindent AV, NRF and AM acknowledge the Project "Reduced Order Modeling and Deep Learning for the real-time approximation of PDEs (DREAM)" (Starting Grant No. FIS00003154), funded by the Italian Science Fund (FIS) - Ministero dell’Università e della Ricerca (MUR). NRF and AM are members of GNCS (Gruppo Nazionale per il Calcolo Scientifico) and gratefully acknowledge its support. 
The present research is also part of the activities of project Dipartimento di Eccellenza
2023-2027, Department of Mathematics, Politecnico di Milano, funded by MUR.

\bibliographystyle{cas-model2-names}

\bibliography{references}



\appendix

\newtheorem{theoremA}{Theorem}[section]
\newtheorem{lemmaA}{Lemma}[section]
\newtheorem{propositionA}{Proposition}[section]
\newtheorem{remarkB}{Remark}[section]

\section{Proofs of Section~\ref{sec:theory}}
\label{appendix:proofs}
\;\\
\textbf{Proof of Lemma~\ref{lemma:decomposition}.}
Fix any $\bmu\in\pspace.$
Using the Kantorovich-Rubinstein duality \cite[Chapter 5]{villani2009optimal}, we can write the Wasserstein-1 distance between the distributions of $\bu_{\bmu}$ and $\bV\phi(\bmu,\bz)$ as
\begin{equation}
\label{eq:kantorovich}
\wasserstein(\bu_{\bmu}, \bV\phi(\bmu,\bz)) = \sup_{Q\in \mathscr{B}_1(N_h)} \left| \mathbb{E}_{\bu_{\bmu}}[Q(\bu_{\bmu})] - \mathbb{E}_{\bz}[Q(\bV\phi(\bmu, \bz))] \right|
\end{equation}
where $\mathscr{B}_1(m)$ denotes the set of all Lipschitz continuous maps $Q: \mathbb{R}^{m} \to \mathbb{R}$ with Lipschitz constant no larger than 1. Notice that the two expectations in \eqref{eq:kantorovich} are taken with respect to $\bu_{\bmu}$ and $\bz$, respectively. We now notice that, for every $Q\in \mathscr{B}_1(N_h)$, by the triangle inequality
\begin{equation}
\label{eq:qtriangle}
\left| \mathbb{E}_{\bu_{\bmu}}[Q(\bu_{\bmu})] - \mathbb{E}_{\bz}[Q(\bV\phi(\bmu,\bz))] \right|\le\left| \mathbb{E}_{\bu_{\bmu}}[Q(\bu_{\bmu})] - \mathbb{E}_{\bu_{\bmu}}[Q(\bV\bV^\top\bu_{\bmu})] \right| + \left| \mathbb{E}_{\bu_{\bmu}}[Q(\bV\bV^\top\bu_{\bmu})] - \mathbb{E}_{\bz}[Q(\bV\phi(\bmu,\bz))] \right|.
\end{equation}
The first term in the above can be bounded as
\begin{equation}
\label{eq:firstbound}
\left| \mathbb{E}_{\bu_{\bmu}}[Q(\bu_{\bmu})] - \mathbb{E}_{\bu_{\bmu}}[Q(\bV\bV^\top\bu_{\bmu})] \right| 
\le \mathbb{E}_{\bu_{\bmu}}\left[ \left| Q(\bu_{\bmu}) - Q(\bV\bV^\top\bu_{\bmu}) \right| \right]
\le \mathbb{E}_{\bu_{\bmu}}\| \bu_{\bmu} - \bV\bV^\top\bu_{\bmu} \|,
\end{equation}
since $Q\in \mathscr{B}_1(N_h)$. The second term in Eq. \eqref{eq:qtriangle}, instead, can be re-written as
\begin{equation*}\left| \mathbb{E}_{\bu_{\bmu}}[Q(\bV\bV^\top\bu_{\bmu})] - \mathbb{E}_{\bz}[Q(\bV\phi(\bmu,\bz))] \right|=\left| \mathbb{E}_{\bu_{\bmu}}[q(\bV^\top\bu_{\bmu})] - \mathbb{E}_{\bz}[q(\phi(\bmu,\bz))] \right|,\end{equation*}
where $q:\mathbb{R}^n\to\mathbb{R}$ is given as $q(\mathbf{c}):=Q(\bV\mathbf{c}).$ It is straightforward to see that $q\in \mathscr{B}_1(n)$. In fact, because $\bV$ is orthonormal, for every $\mathbf{c},\mathbf{c}'\in\mathbb{R}^n$ one has
$
|q(\mathbf{c}) - q(\mathbf{c}')| = |Q(\bV\mathbf{c}) - Q(\bV\mathbf{c}')| \leq \|\bV\mathbf{c} - \bV\mathbf{c}'\| = |\bV(\mathbf{c} - \mathbf{c}')\| = \|\mathbf{c} - \mathbf{c}'\|.
$
It follows that
\begin{equation}
\label{eq:secondbound}
\left| \mathbb{E}_{\bu_{\bmu}}[q(\bV^\top\bu_{\bmu})] - \mathbb{E}_{\bz}[q(\phi(\bmu,\bz))] \right|\le\sup_{\tilde{q}\in \mathscr{B}_1(n)}\left| \mathbb{E}_{\bu_{\bmu}}[\tilde{q}(\bV^\top\bu_{\bmu})] - \mathbb{E}_{\bz}[\tilde{q}(\phi(\bmu,\bz))] \right|=\wasserstein(\bV^\top\bu_{\bmu},\phi(\bmu,\bz)).\end{equation}
Therefore, combining \eqref{eq:qtriangle} with \eqref{eq:firstbound} and \eqref{eq:secondbound} we obtain
\begin{equation*}\left| \mathbb{E}_{\bu_{\bmu}}[Q(\bu_{\bmu})] - \mathbb{E}_{\bz}[Q(\bV\phi(\bmu,\bz))] \right| \le \mathbb{E}_{\bu_{\bmu}}\| \bu_{\bmu} - \bV\bV^\top\bu_{\bmu} \| + \wasserstein(\bV^\top\bu_{\bmu},\phi(\bmu,\bz)).\end{equation*}
Since $Q\in \mathscr{B}_1(N_h)$ was arbitrary, passing at the supremum yields the desired conclusion.\qed
\;\\\\
\textbf{Proof of Corollary~\ref{corollary:upperbound}.}
By integrating \eqref{eq:main inequality} with respect to $\bmu$, we can bound the left-hand-side in Eq. \eqref{eq:corollary} with $\int_{\pspace}\mathbb{E}\|\bu_{\bmu}-\bV\bV^\top\bu_{\bmu}\|\pdensity(\bmu)d\bmu + \int_{\pspace}\wasserstein(\bV^\top \bu_{\bmu}, \;\phi(\bmu, \bz))\pdensity(\bmu)d\bmu.$ By Jensen's inequality, the first term can be bounded further as \begin{equation*}\int_{\pspace}\mathbb{E}\|\bu_{\bmu}-\bV\bV^\top\bu_{\bmu}\|\pdensity(\bmu)d\bmu\le \left[\int_{\pspace}\mathbb{E}\|\bu_{\bmu}-\bV\bV^\top\bu_{\bmu}\|^2\pdensity(\bmu)d\bmu\right]^{1/2}.\end{equation*} However, by standard arguments based on the Eckart-Young-Schmidt Theorem, the POD matrix $\mathbf{V}$ satisfies $\int_{\pspace}\mathbb{E}\|\bu_{\bmu}-\bV\bV^\top\bu_{\bmu}\|^2\pdensity(\bmu)d\bmu=\sum_{j>n}\lambda_j$ (see, e.g. \cite[Proposition 6.3]{quarteroni2015reduced}), thus yielding the claimed inequality.\qed
\;\\\\
\textbf{Proof of Theorem~\ref{theorem:universality}.}
We first notice that, since $\mathbf{V}$ is orthonormal and the random vectors $\bu_{\bmu}$ are uniformly bounded, one has $\|\mathbf{V}^\top\bu_{\bmu}\|\le R$ almost surely for all $\bmu\in\pspace$. In particular, being $B(0,R)\subset \mathbb{R}^n$ the ball of radius $R$ centered at the origin, one has $\cup_{\bmu\in\pspace}\text{supp}(\prob_{n,\bmu})\subseteq B(0,R)$, meaning that the union of the supports of the reduced distributions is bounded. 
We then notice that the map $\bmu\to\prob_{n,\bmu}$ is continuous in the $\wasserstein$ metric. Indeed, since $Q=q\circ\mathbf{V}^\top$ is 1-Lipschitz whenever $q:\mathbb{R}^{n}\to\mathbb{R}$ is 1-Lipschitz, one has
\begin{equation*}\wasserstein(\mathbb{P}_{n,\bmu},\;\mathbb{P}_{n,\bmu'})=\sup_{q\in \mathscr{B}_1(n)}|\mathbb{E}[q(\mathbf{V}^\top\bu_{\bmu})]-\mathbb{E}[q(\mathbf{V}^\top\bu_{\bmu'})]|\le \sup_{Q\in \mathscr{B}_1(N_h)}|\mathbb{E}[Q(\bu_{\bmu})]-\mathbb{E}[Q(\bu_{\bmu'})]|=\wasserstein(\mathbb{P}_{\bmu},\;\mathbb{P}_{\bmu'})\end{equation*}
where, as before, $\mathscr{B}_1(m)$ denotes the set of all 1-Lipschitz maps from $\mathbb{R}^m\to\mathbb{R}.$ Hence, the continuity of $\bmu\mapsto\prob_{\bmu}$ transfers to $\bmu\mapsto\prob_{n,\bmu}.$ The conclusion now follows from Corollary~\ref{corollary:upperbound} and the universality of conditional normalizing flows, cf. Lemma~\ref{lemma:universality} in Appendix~\ref{appendix:universality}.\qed

\section{Universality of conditional normalizing flows}
\label{appendix:universality}
To be derived, Corollary 2 requires the universality of conditional normalizing flows. To the best of our knowledge, however, such result is not readily available in the existing literature. Prior works, in fact, either focus on the unconditional case, as in \cite{koehler2021representational, lee2021universal}, or discuss parametric problems in terms of diffeomorphism learning, as in \cite{lyu2022cflows}, rather than in the context of conditional distributions. To bridge this gap, we derive a simple universality result, sufficient for our purposes, which directly descends from \cite[Theorem 1]{koehler2021representational}.

\begin{lemmaA}
    \label{lemma:universality}
    Let $\pspace\subset\mathbb{R}^p$ be a compact set equipped with a suitable density $\pdensity=\pdensity(\bmu)$. Assume that each $\bmu\in\pspace$ is associated with some probability distribution $\prob_{\bmu}$ defined over $\mathbb{R}^n$. Further assume that each of those probability distributions is absolutely continuous with respect to the Lebesgue measure, and that the union of their supports is bounded. Finally, assume that the map $\bmu\mapsto\prob_{\bmu}$ is continuous with respect to the $\wasserstein$ metric. Let $\mathbb{Q}$ the standard Gaussian distribution in $\mathbb{R}^n$. Then, for every $\epsilon>0$, there exists a conditional normalizing flow $\phi$ with three alternating affine coupling layers, and maps $s,t$ represented by ReLU networks, such that
    \begin{equation*}\int_\pspace\wasserstein(\phi(\cdot,\bmu)\#\mathbb{Q},\;\prob_{\bmu})\pdensity(\bmu)d\bmu<\epsilon.\end{equation*}
\end{lemmaA}

\begin{proof}
Let $\epsilon'>0.$ By compactness and continuity, the map $\bmu\mapsto\prob_{\bmu}$ is uniformly continuous. Hence, there exists some $\delta=\delta(\epsilon')$ for which $\wasserstein(\prob_{\bmu}, \prob_{\bmu'})<\epsilon'$ whenever $\|\bmu-\bmu'\|<\delta.$ Leveraging compactness, let $\{\bmu_1,\dots,\bmu_m\}$ be a $\delta$-net for $\pspace$. Let $\prob_1,\dots,\prob_m$ be the corresponding probabilities, so that $\prob_j:=\prob_{\bmu_j}.$ In light of the universality of unconditional normalizing flows, for each of such probabilities there exists a Real NVP $\phi_j$ with three coupling layers, featuring ReLU networks only, such that $\wasserstein(\phi_j\#\mathbb{Q},\;\prob_{j})<\epsilon'$: see \cite[Theorem 1]{koehler2021representational} and Remark \ref{remark:universality}.

For each $j=1,\dots,m$, let $s_{j,1},s_{j,2},s_{j,3}$ and $t_{j,1},t_{j,2},t_{j,3}$ be the ReLU networks in the three affine coupling layers of $\phi_j.$
Define the parametric maps
\begin{equation}
\label{eq:exogeneous}
\tilde{s}_{k}(\bmu, \bx):=\sum_{j=1}^m\mathbf{1}_{B_j}(\bmu)s_{j,k}(\bx),\qquad\tilde{t}_{k}(\bmu, \bx):=\sum_{j=1}^m\mathbf{1}_{B_j}(\bmu)t_{j,k}(\bx)
\end{equation}
where $\mathscr{B}_1,\dots,B_m$ is a partition obtained by considering the $m$ balls of radius $\delta$ centered over the $\delta$-net $\bmu_1,\dots,\bmu_m$, namely $B_j:=\tilde{B}_j\setminus(\cup_{l=1}^{j-1}B_l)$ for $j\ge2$, $\mathscr{B}_1:=\tilde{B}_1$ and $\tilde{B}_j:=\{\bmu\in\mathcal{P}\;:\;\|\bmu-\bmu_j\|<\delta\}.$ 
Let $\tilde{\phi}$ be the \emph{exogeneous} conditional normalizing flow having $\tilde{s}_{k}$ and $\tilde{t}_{k}$ as maps for the alternating affine coupling layers. We call $\tilde{\phi}$ exogenoeus because its defining maps are not given by neural network architectures but rather more general functions. 

Given any $\bmu\in\pspace$, let $j(\bmu)$ be the unique index in $\{1,\dots,m\}$ for which $\bmu\in B_{j(\bmu)}$. Then, for every $\bmu\in\pspace$ we have
\begin{align}
\label{eq:delta-net}
\wasserstein(\tilde{\phi}(\cdot,\bmu)\#\mathbb{Q},\;\prob_{\bmu})&\le \wasserstein(\tilde{\phi}(\cdot,\bmu)\#\mathbb{Q},\;\prob_{j(\bmu)})+\wasserstein(\prob_{j(\bmu)},\prob_{\bmu})\\
&=
\wasserstein(\phi_{j(\bmu)}\#\mathbb{Q},\;\prob_{j(\bmu)})+\wasserstein(\prob_{\bmu_{j(\bmu)}},\prob_{\bmu})\le 2\epsilon'\nonumber.\end{align}


We now notice that, by assumption, there exists a compact set $\mathcal{S}$ such that $\mathcal{S}\supset\text{supp}(\mathbb{P}_{\bmu})$ for all $\bmu\in\pspace.$ Having fixed such $\mathcal{S}$, let $\eta>0$ be a parameter, whose value shall be chosen later on. Leveraging universality of ReLU networks, for every $k=1,\dots,3$, let $s^\eta_{k}$ and
$t^\eta_{k}$ be two ReLU networks, function of $\bmu$ and $\bx$, such that
\begin{equation}
\label{eq:density}
\int_{\pspace}\sup_{\bx\in\mathcal{S}}\|\tilde{s}_{k}(\bmu,\bx)-s^\eta_{k}(\bmu,\bx)\|\pdensity(\bmu)d\bmu< \eta\qquad\text{and}\qquad\int_{\pspace}\sup_{\bx\in\mathcal{S}}\|\tilde{t}_{k}(\bmu,\bx)-t^\eta_{k}(\bmu,\bx)\|\pdensity(\bmu)d\bmu<\eta.
\end{equation}
This is possible since the maps $\tilde{s}_{k}$ and $\tilde{t}_{k}$ are all continuous in $\bx$ and piecewise constant in $\bmu.$ Without loss of generality, we also require that $|s_k^\eta|,|t_k^\eta|\le M$ uniformly across $\pspace\times\mathcal{S}$, with $M$ independent of $\eta$ (but dependent on $\epsilon'$ through $\tilde{s}$ and $\tilde{t}$, which depend on $\phi_1,\dots,\phi_m$). Such boundness can be ensured in multiple ways. For instance, one may first approximate the indicator functions in \eqref{eq:exogeneous} using continuous maps (e.g., partitions of unity), obtaining suitable continuous approximants of $\tilde{s}_k$ and $\tilde{t}_k$, which can then be approximated uniformly on $\pspace\times\mathcal{S}$ by ReLU networks. 
\\\\
Let now $\phi_\eta$ be the conditional normalizing flow associated with $\{s_k^\eta,t_k^\eta\}_{k=1}^{3}.$ 
For every $\bmu\in\mathcal{P}$ let $\tilde{g}_{\bmu}$ be the density of $\tilde{\phi}(\cdot,\bmu)\#\mathbb{Q}$ and let $g^\eta_{\bmu}$ be that of $\phi_\eta(\cdot,\bmu)\#\mathbb{Q}.$ We first notice that, by \eqref{eq:density}, for almost every $\bmu\in\mathcal{P}$ we have 
$s_{k}^\eta(\bmu,\cdot)\to \tilde{s}_k(\bmu,\cdot)$ and $t_{k}^\eta(\bmu,\cdot)\to \tilde{t}_k(\bmu,\cdot)$ uniformly on $\mathcal{S}$ as $\eta\to0^+$. 
Because the density of a normalizing flow depends smoothly on the maps $s$ and $t$ in its coupling layers, it is straightforward to see that such convergence directly transfers to the conditional densities, meaning that $\sup_{\bx\in\mathcal{S}}\|g_{\bmu}^\eta(\bx)-\tilde{g}_{\bmu}(\bx)|\to0$ as $\eta\to0^+$ for almost every $\bmu\in\mathcal{P}$. However, because all densities are uniformly bounded by some constant $C=C(M)>0$, this immediately implies
\begin{equation*}\lim_{\eta\to0^+}\int_{\pspace}\sup_{\bx\in\mathcal{S}}\|\tilde{g}_{\bmu}(\bx)-g_{\bmu}^n(\bx)\|\pdensity(\bmu)d\bmu=0,\end{equation*}
by dominated convergence.
In particular, because $\mathcal{S}$ is bounded, we also have
$\int_{\pspace}\wasserstein(\tilde{\phi}(\cdot,\bmu)\#\mathbb{Q},\;\phi_\eta(\cdot,\bmu)\#\mathbb{Q})\pdensity(\bmu)d\bmu\to0$ as $\eta\to0^+$. Thus, we may choose $\eta=\eta_*(\epsilon')$ such that $\int_{\pspace}\wasserstein(\tilde{\phi}(\cdot,\bmu)\#\mathbb{Q},\;\phi_\eta(\cdot,\bmu)\#\mathbb{Q})\pdensity(\bmu)d\bmu<\epsilon'.$ Combining the latter with \eqref{eq:delta-net} we obtain
\begin{equation*}\int_{\pspace}\wasserstein(\phi_{\eta_*}(\cdot,\bmu)\#\mathbb{Q},\;\prob_{\bmu})\pdensity(\bmu)d\bmu\le\int_{\pspace}\wasserstein(\tilde{\phi}(\cdot,\bmu)\#\mathbb{Q},\;\phi_{\eta_*}(\cdot,\bmu)\#\mathbb{Q})\pdensity(\bmu)d\bmu+\int_{\pspace}\wasserstein(\tilde{\phi}(\cdot,\bmu)\#\mathbb{Q},\;\prob_{\bmu})\pdensity(\bmu)d\bmu<3\epsilon'.\end{equation*}
Setting $\epsilon':=\epsilon/3$ and letting $\phi:=\phi_{\eta_*(\epsilon/3)}$ yields the conclusion.
\end{proof}

\begin{remarkB}
    \label{remark:universality}
Our definition of Real NVPs differs slightly from that in \cite{koehler2021representational}, as we require a terminal activation $\rho_k$ in the scaling network $s$ of each coupling layer. Nevertheless, \cite[Theorem 1]{koehler2021representational} extends naturally to our setting whenever $\rho_k$ is continuous and unbounded below. Indeed, the universality proof in \cite{koehler2021representational} uses only constant negative scalings of arbitrarily large magnitude. It therefore suffices to show that there exists $M\in\mathbb{R}$ such that, for every $c<M$, our definition admits a scaling network $s_c$ satisfying $s_c(\boldsymbol{x})\equiv[c,\dots,c]$ for all inputs $\boldsymbol{x}$. To this end, set $M:=\rho_k(0)$. By continuity and unboundedness, it follows that $(-\infty,M]\subseteq\rho_k(\mathbb{R})$. Hence, for every $c<M$, there exists $\tau_c\in\mathbb{R}$ such that $\rho_k(\tau_c)=c$. The desired network can then obtained by setting all weights and biases to zero, except for the bias of the output layer, which can be constructed as $[\tau_c,\dots,\tau_c]$.

From a practical standpoint, the widespread use of bounded activations may pose some limitations in this direction, but their incompatibility with the above assumption is not as severe as it may sound. The universality discussed so far, in fact, covers a broad class of target distributions, including absolutely continuous distributions with disconnected supports. Such targets cannot be represented exactly by a diffeomorphic transformation of a Gaussian, since diffeomorphisms preserve connectedness of the support. For this reason, their approximation generally requires transformations whose Jacobians become increasingly singular. In the construction discussed above, this behavior arises as $c\to-\infty$, causing the corresponding scaling factors to approach zero. Bounded activations restrict access to this limiting regime, but they can still be very well-suited for approximating a broad class of different distributions.
\end{remarkB}

\end{document}